\input amstex
 \documentstyle{amsppt}
 \document

\def\RBS{{RBS}}

 \def\RR{{\Bbb R}}
\def\ZZ{{\Bbb Z}}
\def\QQ{{\Bbb Q}}
\def\CC{{\Bbb C}}
\def\CO{{\Cal O}}
\def\CA{{\Cal A}}
\def\CE{{\Cal E}}
\def\Ct{{\Cal C}}

\def\d{{\partial}}

\def\Dbar{{\overline D}}
\def\Xbar{{\overline X}}

\def\XSig{{X^{\text{tor}}}}
\def\XSigexc{{X^{\text{tor,exc}}}}

\def\Dbarred{{\overline D^{\text{red}}}}
\def\Ebarred{{\overline{\Cal E}^{\,\text{red}}}}
\def\Ebarexc{{\overline{\Cal E}^{\text{exc}}}}
\def\Xbarred{{\overline X^{\text{red}}}}
\def\Xbarexc{{\overline X^{\text{exc}}}}

\def\Pred{{P^{\text{red}}}}
\def\nabla{{\triangledown}}

\def\prp{{\prime\prime}}
\def\Int{{\roman{Int}}}

\def\RBS{{{RBS}}}

\def\GhP{{G_{h,P}}}

\def\GlP{{G_{\ell,P}}}
\def\GlQ{{G_{\ell,Q}}}
\def\Gam{{\Gamma}}

\def\Sig{{\Sigma}}
\def\CA{{\Cal A}}

\def\Xtorexc{{X^{\text{tor,exc}}}}

\def\e{{\text{exc}}}
\def\r{{\text{red}}}
\def\t{{\text{tor}}}
\def\te{{\text{tor,exc}}}
\def\im{{\,\roman{im}}}

\def\GCQ{{\roman{GCQ}}}
\def\LCM{{\roman{LCM}}}
\def\L{{\roman{LCM}}}
\def\Lb{{\roman{LCMb}}}
\def\sq{{$\quad\square$}}
\def\Cp{{\frak{Cp}}}
\def\PCp{{\frak{PCp}}}
\def\i{{\iota}}
\def\Y"{{Y^{\prime\prime}}}
\def\back{{\backslash}}

\hfuzz=5pt

\define\isoarrow{{\,\overset\sim\to\longrightarrow\,}}
\def\-{{-1}}

\pageheight{8.0in}
\voffset -0.3cm

 \NoBlackBoxes
\magnification\magstep1

\centerline{\bf On the reductive Borel-Serre compactification, II:}
\centerline{\bf Excentric quotients and least common modifications}
\bigskip

\centerline{Steven Zucker\footnote{Supported in part by the National Science
Foundation,
through Grant DMS9820958}}

\centerline{Department of Mathematics, Johns Hopkins University, Baltimore,
MD 21218
USA\footnote{e-mail address: zucker\@jhu.edu}}

\bigskip\bigskip

{\eightpoint {\it Abstract.} Let $X$ be a locally symmetric variety, i.e.,
the quotient of a bounded
symmetric domain by a (say) neat arithmetically-defined group of isometries.
Let $\Xbarexc$ and $X^\te$ denote its excentric Borel-Serre and toroidal
compactifications respectively.  We determine their least common
modification and use it to prove a conjecture of Goresky and Tai
concerning canonical extensions of homogeneous vector bundles. In the process,
we see that $\Xbarexc$ and $X^\te$ are homotopy equivalent.}

\topmatter
\toc
\widestnumber\subhead{2.A.}
\specialhead {} Introduction\endspecialhead
\head 1.  The category of compactifications of a space\endhead
\subhead 1.1. Fundamentals\endsubhead
\subhead 1.2. Example: Two compactifications of a simplicial cone\endsubhead
\subhead 1.3. Discrete quotients and diagonality\endsubhead
\head 2.  Locally-symmetric varieties and their compactifications\endhead
\subhead 2.1. Boundary components and the Baily-Borel compactification
\endsubhead
\subhead 2.2. Toroidal compactifications and their quotients\endsubhead
\subhead 2.3. The Borel-Serre compactification and its quotients\endsubhead
\subhead 2.A. Appendix to \S2: Hybrid compactifications\endsubhead
\head 3.  The least common modification of $\Xbarexc$ and $X^\te$\endhead
\subhead 3.1. Main result and elements of the proof\endsubhead
\subhead 3.2. Real and complex toroidal embeddings; dual compactifications
\endsubhead
\subhead 3.3. Adjustments to duality\endsubhead

% \subhead 3.4. Enter Jordan algebras\endsubhead
% \subhead 3.5. The role of Siegel sets\endsubhead

\subhead 3.4. Quotients by $\Gam(\GlP)$ and the fiber\endsubhead
\subhead 3.5. Fibers over the strata of $X^*$\endsubhead
\head 4.  Canonical extension of homogeneous vector bundles\endhead
\subhead 4.1. Standard notions\endsubhead
\subhead 4.2. Torus actions and torus embeddings, revisited\endsubhead
\subhead 4.3. Identification of canonical extensions\endsubhead
% \subhead 4.4. Compatibilities\endsubhead
% \subhead 4.5. Descent to $X^\te$\endsubhead
\subhead 4.4. The Goresky-Tai conjecture\endsubhead
\subhead 4.5. Proof of the excentric Goresky-Tai conjecture\endsubhead
\specialhead {} References\endspecialhead
\endtoc
\endtopmatter

\tenpoint\baselineskip=16pt

\centerline{\bf Introduction}
\medskip

This article is the continuation of [Z4] in the direction it left
open, namely the Goresky-Tai conjecture [GT:\,9.5] (Conjecture A
in [Z4]). The conjecture concerns two very different classes of
compactifications of a locally symmetric variety $X$, and the
corresponding notions of canonical extension of homogeneous vector
bundles $\CE$ on $X$.  One class, the good toroidal
compactifications $X^\t$ [AMRT], comes from $X$ as an algebraic
variety, with $X^\t$ a smooth complex projective variety.  The
other, the Borel-Serre compactification $\Xbar$ [BS] (a
manifold-with-corners) and its reductive quotient $\Xbarred$ (a real
stratified space; see [Z4:\S\S1,5]), being defined also for non-Hermitian
$X$, are more general. The Borel-Serre spaces have
boundary strata of odd real-codimension, so are quite far from
being complex-analytic spaces.  Nonetheless, one can view Conjecture A as
saying
that $\Xbarred$ is more fundamental than $X^\t$, at least
as far as homogeneous vector bundles are concerned: if $\CE^\t$
and $\Ebarred$ denote the respective bundle extensions of $\CE$ to
$X^\t$ and $\Xbarred$, then for certain continuous mappings
$h:X^\t\to \Xbarred$ (described below), one has $\CE^\t\simeq
h^*\Ebarred$.  In other words, $\CE^\t$ is determined as a topological vector
bundle by $\Ebarred$; their respective Chern classes are correspondingly
related.

A problem in dealing with Conjecture A is that $X^\t$ and $\Xbarred$ are
generally so
different.  Here, I am thinking {\it beyond} the result of Lizhen Ji [J] (in
essence,
the conjecture from [HZ2:(1.5.8)] revised),
which implies that the greatest common quotient of $X^\t$ and
$\Xbarred$ is the obvious one, namely the Baily-Borel
compactification $X^*$ (see (2.2.18) and (2.3.5)). (Each of $X^\t$ and $\Xbar$
can be viewed as a resolution of $X^*$, in the senses of complex algebraic
varieties and real stratified spaces, resp.)
In [HZ2:\S1], we introduced two auxiliary compactifications of $X$
which we now denote $X^\te$ and $\Xbarexc$, and call the {\it excentric}
toroidal and Borel-Serre compactifications respectively (see our Section 2).
They are straightforward quotients of $X^\t$ and $\Xbar$ respectively,
with $\Xbarexc$ mapping onto $\Xbarred$. The two excentric
compactifications could be seen to have much in common (see, e.g., (2.3.11)).
It was a naive
suspicion at the time that they might be isomorphic as
compactifications of $X$; it was important to realize that in
general they are not.

The goal in [GT] was to determine the least common modification
(language as in [HZ2:(1.1.4)]) of $\Xbarred$ and $X^\t$, denoted
$\L(\Xbarred,X^\t)$. They showed that the canonical mapping
$\L(\Xbarred,X^\t)\to X^\t$ is a homotopy equivalence.  Then
homotopy inverses are used to yield the mappings $h$ above, via
the composite
$$
X^\t\to\L(\Xbarred,X^\t)\to \Xbarred.
$$
In this article, we determine that $\L(\Xbarexc,X^\te)\to X^\te$
is a homotopy equivalence (Theorem 3.1.1), for more or less the same reason.
However, our argument shows that it is a natural consequence of
the calculation of the least common modification of dual
compactifications of a simplicial cone (1.2) (see also
[HZ1:\,2.3]). We thereby obtain a homotopy class of mappings $k:
X^\te\to\Xbarexc$. We recover the $\L$ determination from [GT] by
considerations of {\it LCM-basechange} (1.1.8), the name we give
to the determination for three compactifications $X_1$, $X_2$, and
$X_3$, with $X_3$ a quotient of $X_2$, whether the canonical
embedding
$$
\L(X_1,X_2)\hookrightarrow \L(X_1,X_3)\times_{X_3} X_2
$$
is an equality.  On the other hand, our theorem does not seem to follow
directly from [GT].

The preceding is made more significant by the fact that there is a
canonical extension $\CE^\te$ of the homogeneous vector bundle
$\CE$ to $X^\te$ that pulls back to $\CE^\t$ via $X^\t\to X^\te$
(4.3), and likewise $\Ebarexc$ on $\Xbarexc$ that is the pullback
of $\Ebarred$ via $\Xbarexc\to\Xbarred$.  We can then formulate
the analogue of Conjecture A (which we tentatively label
Conjecture A$'$: $k^*\Ebarexc\simeq \CE^\te$).  It is not hard to
see that Conjecture A is a consequence of Conjecture A$'$, but the
latter is talking about the pair of excentric compactifications,
which more resemble each other.  We give a proof of Conjecture
A$'$ in the last section of this article.

The organization of this article is as follows.  In Section 1, we
treat the compactifications of a space $X$ as a category in (1.1),
in which taking the $\L$ is a bifunctor.  Particularly important
is the discussion of when $\L$-basechange occurs.  In (1.2), we
make the critical calculation of the $\L$ of dual
compactifications of a simplicial cone.  This is followed in (1.3)
by a treatment of the somewhat elusive notion that we call {\it
diagonality}, where the commutation of $\L$ and discrete quotients
is analyzed.

Section 2 contains a description of the essential features of the
various classes of compactifications of a locally symmetric
variety: Baily-Borel (2.1), toroidal (2.2), and Borel-Serre (2.3).
The excentric quotients are recalled in (2.2.18) and (2.3.5).  In
an appendix, we introduce some additional Borel-Serre quotients
(called {\it hybrid compactifications}), lying between $\Xbarexc$
and $\Xbarred$, that are inspired by the determination of the
greatest common quotient of $\Xbar$ and $X^\t$ in [J].

Section 3 is devoted to the proof of our Theorem 3.1.1, which
asserts that $\L(\Xbarexc,X^\te)\to X^\te$ is a homotopy
equivalence.  We use the outcome of (1.2) to give in (3.2) an
approximation of the proof, revising it at the boundary in (3.3).
In Sections (3.4) and (3.5), we finish the proof and determine
consequences that follow by $\L$-basechange. We emphasize the
substantial conclusion given in Corollary 3.5.11: $\Xbarexc$ {\it
and $X^\te$ are homotopy equivalent}.

Finally, in Section 4, we carry out the toroidal construction of
[AMRT] on (the total space of) the vector bundle $\CE$, yielding a
vector bundle $\CE^\t$ on $X^\t$. We show that $\CE^\t$ is the
canonical extension of $\CE$ in the sense of [Mu] in (4.3).  (This
is analogous to what was done in [Z4] to yield $\Ebarred$.)  We
then descend this bundle to $\CE^\te\to X^\te$.  In (4.4), we
verify that Conjecture A$'$ implies the conjecture of Goresky-Tai
(Conjecture A). We complete the task by verifying Conjecture A$'$
in (4.5).

I want to thank Mark Goresky and Lizhen Ji for helpful
correspondence and discussions. The referee is to be commended for
his or her thorough job of reading the manuscript, and also for
wisely insisting on a substantial revamping of the exposition. To
my surprise, it was pointed out by the referee (correctly!) that
the verbal description in the first line of [HZ1:\,p.262] is a bit
garbled; the correct assertion was nearby, though, and it is
stated correctly (in slightly different notation) here in (2.2).
\medskip

\demo{Apology} It is only a little unnatural that this article is
appearing well after its sequel [Z5], which was written for a
special volume. The latter is largely independent, referring only
to Corollary 3.5.11 for a conditional assertion.  There is also
reference to Corollary 3.5.11 in [Z6:\,Prop.~2].
\enddemo
\demo{Comments on {\rm [Z4]}}  i) {\bf Erratum:} The quantity
$\delta$ in [Z4:(3.1.4)] should be described as, and taken to be,
the {\it sum} of the positive $\QQ$-roots, not the half-sum (the
``half" appears as $\tfrac 1 p$ when $p=2$). Subsequent statements
involving $\delta$ are correct as written.
\smallskip

ii) About the time [Z4] appeared, R.~Mazzeo asked a familiar
question: ``Why $L^p$-cohomology for $p\ne 2$?"  I think that the
article [Z4] provides a good answer. It is almost certain that I
first tried to calculate the $L^\infty$-cohomology quickly,
finding it to be infinite-dimensional.  Use of large finite $p$
offers a perturbation away from that difficulty, allowing for our
topological interpretation of the $L^p$-cohomology.
\enddemo

\demo{Changes in notation from {\rm [Z4]}}  The simultaneous
treatment of Borel-Serre and toroidal compactifications taxes
one's alphabetical resources, as was seen already in [HZ2].  Note
in particular the following changes of notation:
\smallskip

i) The unipotent radical of a parabolic subgroup $P$ (formerly $U_P$) is now
denoted $W_P$;
the symbol $U_P$ is used here for the {\it center} of the unipotent radical.
This is as in [AMRT].
\smallskip

ii) Symmetric spaces (formerly $X$) are now denoted $D$; the
symbol $X$ is now used for the arithmetic quotient $\Gam\back D$,
which was denoted $M$ in [Z4].
\smallskip

iii) The reductive Borel-Serre compactification of $X$ is denoted here
$\Xbarred$
(formerly $M^\RBS$ of $M$).  The ``bar" (overline) always indicates a
construction
of Borel-Serre type.
\enddemo
\bigskip

% The basic disparity between $\Xbarexc$ and $\Xtorexc$ is given, in essence,
% by [HZ2,(1.5.4)],
% which we paraphrase now.  If one places the closure of $|\widehat \Sigma_P|
% at infinity in
% $|\Sigma_P| \supseteq C_P$ as the locus of limits of geodesic orbits, the
% fibers
% of ${}^<\!Z^\te_P$ over $\Cal A_P$ are proper subsets of $\overline
% \Gamma_{\ell ,P}
% \backslash \widehat C_P$ unless the latter is compact (i.e., unless
% $(G_{\ell,P})$
% is $\Bbb Q$-anisotropic).  Moreover, the set of limit
% points of orbits of the geodesic action in the fiber is discrete.

\centerline{\bf 1. The category of compactifications of a space}
\medskip

The main purpose of this section is to recall (cf.~[HZ2,\S1]) the notion of
the least common
modification of two compactifications of the same topological space, and to
develop
properties that this notion has when a third compactification is invoked.
\medskip

{\bf (1.1)} {\bf Fundamentals.} We begin with some basic
terminology.  Let $X$ be a locally compact Hausdorff space
(non-compact) with a countable base for its topology.  One defines
a category $\Cp(X)$, the {\it compactifications of} $X$, as
follows.  The objects are pairs $(Y,\i)$, with $Y$ a compact
Hausdorff space, and $\i:X\to Y$ a dense open embedding. We will
write just $Y$ when $\i$ is understood, and view $X$ as a subset
of $Y$.  The {\it interior} of $Y$ is then understood to be $X$;
the {\it boundary} $\partial Y$ of $Y$ (qua compactification of
$X$) is the complement of $X$ in $Y$, a closed subset of $Y$. To
eliminate pathology, we will assume throughout this article that
we allow only spaces $Y$ (thus also $X$) whose topologies have
countable neighborhood bases. This holds in particular whenever
$Y$, hence also $X$, is metrizable.

A morphism in $\Cp(X)$, called a {\it morphism of compactifications} of $X$,
from $(Y_1,\i_1)$
to $(Y_2,\i_2)$, is a commutative triangle of continuous mappings:
$$\matrix
X & @> \i_1 >> & Y_1\\
 & \underset{\i_2\,\,}\to\searrow & @VVV\\
& & Y_2\endmatrix
$$
Note that there is at most one morphism from $Y_1$ to $Y_2$, as $X$ is dense
in $Y_1$,
and that morphisms in $\Cp(X)$ are surjective by compactness.  There is such a
morphism if
and only if $Y_2$ is a topological quotient of $Y_1$ at its boundary (so
$\partial Y_1$
maps onto $\partial Y_2$) . If so, one also says that $Y_1$ is a {\it
modification of} $Y_2$.

The set of compactifications of $X$ can be seen to form a
(non-distributive) lattice. Specifically, for $Y_1$ and $Y_2$ as
above, one takes $Y_1\lor Y_2$ to be the closure of $X$
(equivalently, any dense subset of $X$) under the diagonal
embedding in $Y_1\times Y_2$.  This is the {\it least common
modification} of the two compactifications, and is also denoted
$\LCM(Y_1,Y_2)$ or $\LCM_X(Y_1,Y_2)$; it admits canonical
morphisms to $Y_1$ and $Y_2$. The notion passes to the set of
homeomorphism classes of compactifications. Similarly, one takes
$Y_1\land Y_2$ to be the {\it greatest common quotient}
$\GCQ(Y_1,Y_2)$ of $Y_1$ and $Y_2$, which can be realized as the
inverse limit of all common quotients of $Y_1$ and $Y_2$ (the set
of such is non-empty, for the one-point compactification of $X$ is
always a quotient of both).  It receives canonical morphisms from
$Y_1$ and $Y_2$.

\proclaim{Lemma 1.1.1} In $\Cp(X)$, the following three statements are
equivalent:
\smallskip

i) There is a morphism $Y_1\to Y_2$.
\smallskip

ii) Via the natural projection, $\LCM(Y_1,Y_2)\simeq Y_1$.
\smallskip

iii) Via the quotient mapping, $\GCQ(Y_1,Y_2)\simeq Y_2$.\sq
\endproclaim

It is useful to introduce the category $\PCp(X)$ of partial
compactifications
of $X$, which contains $\Cp(X)$ as a full subcategory. A {\it partial
compactification}
of $X$ is a space $Y$ (not necessarily compact) containing $X$ as a dense open
subset.
In particular, $X$ itself is an object in $\PCp(X)$.  The notions of boundary,
morphism
and LCM can be
extended verbatim to $\PCp(X)$, though morphisms need not be surjective (thus
are not
necessarily quotient mappings).  Moreover, $\GCQ(Y_1,Y_2)$ need not be defined,
for $Y_1$ and
$Y_2$ may have no common quotients at all. Nonetheless, the following version
of Lemma 1.1.1
holds in $\PCp(X)$:
\proclaim{Lemma 1.1.1.1} Let $Y_1,Y_2\in\PCp(X)$.
\smallskip

i) There is a morphism $Y_1\to Y_2$ if and only if, via the
natural projection, $\LCM(Y_1,Y_2)\simeq Y_1$.
\smallskip

ii) There is a {\rm surjective} morphism $Y_1\to Y_2$ if and only if $\GCQ(Y_1,
Y_2)$ exists
and is isomorphic to $Y_2$ via the natural mapping,
$Y_2\to\GCQ(Y_1,Y_2)$.  \sq
\endproclaim
\demo{{\rm (1.1.1.2)} Example} Let $X$ be the interval $(0,1)$,
$Y_1 = (0,1]$ and $Y_2 = [0,1)$. One has $\LCM(Y_1,Y_2) \simeq X$.
In particular, the canonical morphisms, $\LCM(Y_1,Y_2)\to Y_i$
($i=1,2$), are not surjective.  Also, $\GCQ(Y_1,Y_2)$ is not
defined.
\enddemo

In the sequel, we will concern ourselves only with the $\LCM$.  Also, until
stated
otherwise, assertions are given for $\PCp(X)$. The following is elementary:

\proclaim{Lemma 1.1.2} In $\PCp(X)$, $\LCM(Y_1,Y_2)$ is equal to
$$
\{(y_1,y_2)\in Y_1\times Y_2\,|\,\,\exists\text{ a sequence
$\{x_j\}$ in $X$ with } \i_1(x_j)\to y_1,\,\i_2(x_j)\to
y_2\}.\quad\square
$$
\endproclaim

\demo{{\rm (1.1.3)} Remark} When we want to draw more attention to
the role of the boundaries $\d_1$ and $\d_2$ of $Y_1$ and $Y_2$
respectively, we renotate $\d\LCM(Y_1,Y_2)$ as $\Lb(\d_1,\d_2)$
(``b" as in ``boundary"); in general, whenever we write
$\Lb(\d_1,\d_2)$, $X$ and the way that $\d_1$ and $\d_2$ are
attached to $X$ are understood to be given information.
\enddemo

It is easy to see: \proclaim{Proposition 1.1.4} Let $Y_1$, $Y_2$
and $Z$ be partial compactifications of $X$, such that $Y_1$ and
$Y_2$ admit morphisms to $Z$.  Then $LCM(Y_1,Y_2)$ maps to $Z$,
and the following diagram commutes
$$\CD
\L(Y_1,Y_2) @>>> Y_2\\
@VVV      @VV\varphi_2 V\\
Y_1 @>\varphi_1>> Z \endCD
$$
In other words, $\LCM(Y_1,Y_2)\subseteq Y_1\times_Z Y_2\subset Y_1\times Y_2$.
\endproclaim
\demo{Proof} Let $\{x_j\}$ be a sequence as in Lemma 1.1.2. Then
by continuity, $\{\varphi_i(x_j)\}$ converges to
$\varphi_i(y_i)\in Z$ for both $i=1$ and $i=2$, and our assertion
follows. \sq
\enddemo
As such, we have
\proclaim{Corollary 1.1.5}  In the situation of {\rm Proposition 1.1.4}, let
$z\in Z$.  Let
$U$ be a closed neighborhood of $z$ in $Z$.  Let $Y_1(U)$ denote the subset of
$Y_1$ lying
over $U$, and do likewise for $Y_2$ and $\L(Y_1,Y_2)$.  Then
$$
\L(Y_1,Y_2)(U) = \L(Y_1(U),Y_2(U)).\quad\square
$$
\endproclaim

\noindent In particular, the fiber $\L(Y_1,Y_2)_z$ can be determined from
data
over any neighborhood of $z$.

Next, suppose there is a morphism, $\varphi:Y_1\to Y_2$, and let
$Y_3$ be any third element of $\PCp(X)$.  The induced mapping
$\varphi\times 1: Y_1\times Y_3\to Y_2 \times Y_3$ induces a
morphism $\L(Y_1,Y_3)\to\L(Y_2,Y_3)$. It is important to recognize
that there are no general pullback properties in this situation,
even in $\Cp(X)$. Under the above conditions, one always has
$$
\L(Y_1,Y_3)\subseteq Y_1\times_{Y_2}\L(Y_2,Y_3),\tag 1.1.6
$$
{\it but equality can fail.}  In other words, the commutative diagram
$$\CD
\L(Y_1,Y_3) @>>> \L(Y_2,Y_3)\\
@V\pi_1 VV          @VV\pi_2 V\\
 Y_1 @>\varphi >>    Y_2\endCD
$$
need not be Cartesian, even in $\Cp(X)$.
\demo{{\rm (1.1.7)} Example} A standard situation in which equality in (1.1.6)
fails is
given by taking a morphism $Y_1\to Y_2$ in $\Cp(X)$, and letting $Y_3=Y_1$.
Then $Y_1
\times_{Y_2}\L(Y_2,Y_3)=Y_1\times_{Y_2}Y_1$, but $\LCM(Y_1,Y_3)$ is just $Y_1$.
Indeed, equality
holds in (1.1.6) if and only if $Y_1\simeq Y_2$. (We add that it is easy to
find examples
in $\PCp(X)$ of the preceding sort where equality in (1.1.6) holds
without
$Y_1\simeq Y_2$.  Indeed, take $Y_1$ and $Y_2$ as in (1.1.1.2), and $Y_3=Y_2$.)
\enddemo
\demo{{\rm (1.1.8)} Definition} We will say that {\it
$\LCM$-basechange holds for $Y_3$ with respect to} $Y_1\to Y_2$
when we have equality in (1.1.6), i.e., $\LCM(Y_1,Y_3)=
Y_1\times_{Y_2}\LCM(Y_2,Y_3)$.  When we want to emphasize the role
of the boundaries (the notion is trivial on $X$) as in (1.1.3), we
will say that  {\it $\Lb$-basechange holds for $\d Y_3$ with
respect to} $\d Y_1\to \d Y_2$.
\enddemo
\noindent When this is the case, the fiber of $\pi_1$ at $y\in
Y_1$ is canonically homeomorphic to the fiber of $\pi_2$ at
$\varphi(y)$.
\demo{{\rm (1.1.9)} Remark} In terms of sequences,
as in (1.1.2), $\L$-basechange is equivalent to the following.
Given a sequence in $X$ with limit $y_2\in Y_2$ and limit $y_3\in
Y_3$, then for any $y_1\in Y_1$ that maps to $y_2\in Y_2$, there
is a sequence in $X$ converging to $y_1\in Y_1$ (so to $y_2\in
Y_2$) and to $y_3\in Y_3$.
\enddemo

We have the following complement to Corollary 1.1.5:
\proclaim{Proposition 1.1.10} In the situation of Proposition {\rm 1.1.4},
$$
\L(Y_1,Y_2)_z\subseteq (Y_1)_z\times(Y_2)_z,
$$
with equality for all $z\in Z$
if and only if $\L$-basechange holds for $Y_1$ with respect to $Y_2\to Z$.
\endproclaim
\demo{Proof} We note that by (1.1.6), $\L(Y_1,Y_2)\subseteq\L(Y_1,Z)\times_Z
Y_2 =
Y_1\times_Z Y_2$, with equality if and only if the $\L$-basechange assertion
holds.
The fiber of this over $z$ is the same as that in the statement of the
proposition. \sq
\enddemo

Example (1.1.7) and the failure of surjectivity in morphisms in $\PCp(X)$
suggest that
we return to $\Cp(X)$ as setting for the rest of this Section. We next assert:
\proclaim{Proposition 1.1.11 {\rm ($\LCM$-basechange in a tower)}} Let $\Y"\to
Y'\to Y$
be morphisms
of compactifications of $X$, and $Z$ a fourth compactification of $X$.  Then
$\LCM$-basechange
holds for $Z$ with respect to $\Y"\to Y$ if and only if $\LCM$-basechange
holds
for $Z$ with respect to both $\Y"\to Y'$ and $Y'\to Y$.
\endproclaim
\demo{Proof} From (1.1.6), we always have
$$\multline
\L(\Y",Z)\subseteq \L(Y',Z)\times_{Y'}\Y"\\
\subseteq (\L(Y,Z)\times_YY')\times_{Y'}\Y" = \L(Y,Z)\times_Y\Y".\endmultline
\tag 1.1.11.1
$$
We see that there is equality of the ends if and only if we have equality
at both
inclusion symbols.  This gives our assertion.\sq
\enddemo

One can also talk about {\it two-sided} LCM-basechange.  Let $Y'\to Y$
and $Z'
\to Z$ be morphisms of compactifications of $X$.  There is an embedding
$$
\L(Y',Z')\subseteq Y'\times_{Y}\L(Y,Z)\times_{Z}Z',\tag 1.1.12
$$
through which the embedding $\L(Y',Z')\subseteq Y'\times Z'$ factors.
\proclaim{Proposition 1.1.13} Under the conditions of {\rm (1.1.12)}, equality
holds in
{\rm (1.1.12)} if and only if $\L$-basechange holds for $Z$ and $Z'$ with
respect to
$Y'\to Y$ and for $Y$ and $Y'$ with respect to $Z'\to Z$.
\endproclaim
\demo{Proof} This follows immediately from an elementary fact about fiber
products:
if $A$ and $B$ map to $Z$, and $S$ is a proper subset of $A$, then $S\times_Z
B$ is a
proper subset of $A\times_Z B$.\sq
\enddemo

\demo{{\rm (1.1.13.1)} Remark} It is always the case that in the situation
of (1.1.12),
$$
\L(Y',Z')\simeq\L(\L(Y,Z'),\L(Y',Z)).
$$\enddemo

We give one instance of a simple and useful criterion for $\L$-basechange.
We return to
the situation
of {\rm (1.1.6).}
\proclaim{Proposition 1.1.14} Let $Y_1,Y_2,Y_3\in\Cp(X)$, with morphism
$Y_1\to Y_2$.
Suppose that the compact Lie group $H$ acts on
the pairs
$(Y_1,X)$ and $(Y_3,X)$, in such a way that $H\backslash\d Y_1\simeq\d Y_2$.
If $H$ acts
trivially on $\d Y_3$, then $\L$-basechange holds for $Y_3$ with respect to
$Y_1\to Y_2$.
\endproclaim
\demo{Proof} This is basically (1.1.9). We have that $\L(Y_1,Y_3)\to\L(Y_2,
Y_3)$
is surjective.  If a sequence $\{x_j\}$ in $X$ satisfies $x_j\to y_2$ in $Y_2$
and $x_j\to y_3$
in $Y_3$, a subsequence converges to some $y_1\in Y_1$, and then $y_1$ maps to
$y_2$.
Then for $h\in H$, $hx_j\to hy_1$ in $Y_1$; whereas in $Y_3$, $hx_j\to hy_3 =
y_3$.\sq
\enddemo

It is also possible for the opposite of $\L$-basechange to occur.  That
happens whenever
there is a morphism of compactifications $Y_3\to Y_1$ (cf.~(1.1.7)).  More
generally, one has:
\proclaim{Proposition 1.1.15}
In the situation of {\rm (1.1.6),} $\L(Y_1,Y_3)\to \L(Y_2,Y_3)$ is an
isomorphism
if and only if the canonical mapping $\L(Y_2,Y_3)\to Y_2$ factors through
$Y_1$.
\endproclaim
\demo{Proof} We are in the situation
$$\CD
\L(Y_1,Y_3) @>>>  \L(Y_2,Y_3)\\
@VVV   @VVV\\
Y_1 @>>> Y_2\endCD
$$
When the top arrow is an isomorphism, its inverse can be used to define the
factorization.
Conversely, if we have a morphism $\L(Y_2,Y_3)\to Y_1$, use that with the
canonical
mapping $\L(Y_2,Y_3)\to Y_3$, to give a morphism $\L(Y_2,Y_3)\to \L(Y_1,Y_3)$.
There is a
canonical morphism in the opposite direction.  This gives the isomorphism
we were seeking.\sq
\enddemo
Finally, we state a simple assertion that involves a second initial space:
\proclaim{Proposition 1.1.16}
Let $X'$ be a space with a proper surjection $f:X'\to X$.  Suppose that $Y_1$
and $Y_2$
are compactifications of $X$ and $Y'_1$ and $Y'_2$ compactifications of $X'$
for which $f$
extends to mappings $f_1:(Y'_1,\d Y'_1)\to (Y_1,\d Y_1)$ and $f_2:(Y'_2,\d
Y'_2)\to (Y_2,\d Y_2)$.
Then the induced mapping $\L(Y'_1,Y'_2)\to \L(Y_1,Y_2)$ is surjective.
\endproclaim
\demo{Proof} If $\{x_j\}$ is a sequence converging in both $Y_1$ and $Y_2$,
lift it to a
sequence $\{x'_j\}$ in $X'$.  The hypotheses implies that $\{x'_j\}$ has a
subsequence
that converges in both $Y'_1$ and $Y'_2$.  We are done by Lemma 1.1.2.\sq
\enddemo
\medskip

{\bf (1.2)} {\bf Example: Two compactifications of a simplicial
cone.} The following is an essential calculation, one that
underlies [HZ1,\,2.3], [HZ2,\,(1.5)], [GT,\S 7], and what is to
come in this article.

Let $\sigma$ be a closed simplicial cone that spans a real vector space of
dimension
$d$.  When convenient, we will understand that the origin has been removed.
Let
$\widehat\sigma$ denote the quotient of $\sigma$ by cone dilations, a simplex
that we
can identify with any cross-section of $\sigma\to\widehat\sigma$.  Let
$\{q_j\}$
be a set of $d$ linear functionals that are
non-negative
in $\sigma$ and whose zero-loci define the codimension-one faces of $\sigma$.
These will be
called
{\it linear coordinates} on $\sigma$; they induce barycentric coordinates on
$\widehat\sigma$.
A choice of cross-section enables us to write $\sigma\simeq\widehat\sigma
\times (0,\infty)$
in the usual way.  We compactify $\sigma$ to $\overline\sigma_1$
by taking
$\overline\sigma_1\simeq\widehat\sigma\times (0,\infty]$; thus, $\partial
\overline
\sigma_1\simeq\widehat\sigma$.

A second compactification $\overline\sigma_2$ of $\sigma$ is obtained by
converting to {\it
toroidal
coordinates}, by putting $t_j=\exp (-q_j)$.  This sets up a homeomorphism
$\sigma
\simeq (0,1]^d$, to which we attach the boundary that yields $\overline
\sigma_2
\simeq [0,1]^d$.  Then in toroidal coordinates
$$
\d\overline\sigma_2 \simeq \{ x\in [0,1]^d : t_j=0 \text{ for some $j$}\}\simeq
\bigcup\,\{\tau^
\vee(\widehat\sigma):\widehat\tau\text{ is a face of $\widehat\sigma$}\},
$$
where $\tau^\vee(\widehat\sigma)$ is the closed polyhedral face of $\overline
\sigma_2$
dual to $\widehat\tau$.
For $\tau$ a proper face of $\sigma$, we define $\overline\tau_2$ to be the
result of
the preceding construction when one considers $\tau$ as a subset of the linear
space it
spans.  This allows us to regard $\overline\tau_2$ naturally as a subset of
$\overline\sigma_2$.

It is our present goal to determine $\LCM(\overline\sigma_1,\overline
\sigma_2)$.
\proclaim{Proposition 1.2.1} The complement of $\sigma$ in $\LCM(\overline
\sigma_1,\overline\sigma_2)
\subset\overline\sigma_1\times\overline\sigma_2$ is given by
$$
\Lb(\widehat\sigma_1,\partial\overline\sigma_2)=\partial\LCM(\overline\sigma_1,
\overline\sigma_2)=\bigcup\,\{\widehat\tau
\times
\tau^\vee(\widehat\sigma): \widehat\tau\text{ is a face of }\widehat\sigma\}.
\tag 1.2.1.1
$$
\endproclaim
\demo{Proof}  We make use of Lemma 1.1.1, together with the simple
observation that $q_j$ goes to $\infty$ if and only if $t_j\to 0$.
For each subset $J$ of $\{1,2,..., d\}$, a face $\tau=\tau_J$ of
$\sigma$ is defined by the equations $q_j = 0$ for all $j\in J$.
Put $S=\sum_{1\le j\le d}q_j$, and note that for all $j$, the
ratio $\widehat q_j := q_j/S$ takes values in $[0,1]$. Suppose
that $\{\bold q_k\}$ is a sequence in the cone $\sigma$ that
converges to a boundary point $\bold q_\infty\in\widehat\sigma
=\d\overline \sigma_1$.  Then $\bold q_\infty$ lies in the
interior $\widehat\tau_J^\circ$ of $\widehat \tau_J$ for a
uniquely determined $J$. This implies for the sequence $\{\bold
q_k\}$ that $S\to\infty$, and thus the linear coordinates
$q_j\to\infty$ for all $j\notin J$. In toroidal coordinates, we
have that $t_j\to 0$ for all $j\notin J$.  As the face
$\tau^\vee_J$ in $\d\overline\sigma_2$ dual to $\widehat\tau_J$ is
defined by the equations $t_j = 0$ for all $j\notin J$, we see
that the right-hand side of (1.2.1.1) contains the boundary of
$\LCM(\overline\sigma_1,\overline\sigma_2)$.

To see that one gets every point of the latter space in the boundary, it
suffices to
exhibit a suitable family of curves in $\sigma$.  For each $J$, consider the
family of curves
$\bold q_J(s)$ given in linear coordinates by
$$
q_j(s) = \cases s\, q_{j,0} \text{ if } j\notin J,\\
            q_{j,0} \text{ if } j\in J,\endcases\tag 1.2.1.2
$$
with $q_{j,0}\in\RR^{\ge 0}$ for all $j$. The limit in $\overline\sigma_1$, as
$s\to\infty$,
is the point of $\widehat\tau_J\times\{\infty\}$ with barycentric coordinates
$\widehat q_j=0$ if $j\in J$, and $\widehat q_j=\widehat q_{j,0}$ if $j
\notin J$;
in $\overline\sigma_2$, it is the point of $\tau^\vee_J$ with toroidal
coordinates
$t_j = \exp (-q_{j,0})$ if $j\in J$, and $t_j = 0$ if $j\notin J$.  Thus, all
of
$\widehat\tau\times\tau^\vee(\widehat\sigma)$ is reached.\sq
\enddemo

We illustrate of the above when $\sigma$ is of dimension $d=2$.
Let $q_0$ and $q_1$ be the non-negative linear functionals that
define the edges of $\sigma$. We sketch four cone dilation orbits
and a curve $t_1=const$ (where $q_1$ is constant, so $\widehat
q_1\to 0$), which is of the sort given in (1.2.1.2), in both
$\overline\sigma_1$ and $\overline\sigma_2$:
% \roman{***FIGURE 1}
$$
\tag 1.2.2
$$

\demo{{\rm (1.2.3)} Remarks} i) We can write (1.2.1.1) as
$$
\Lb(\widehat\sigma_1,\partial\overline\sigma_2)=\bigcup\,\{\widehat\tau^\circ
\times
\tau^\vee(\widehat\sigma): \widehat\tau\text{ is a face of }\widehat\sigma\},
\tag 1.2.3.1
$$
where $\widehat\tau^\circ$ denotes the {\it interior} of $\widehat\tau$,
consisting
of those points in $\widehat\tau$ that do not belong to any proper face.
% Though (1.2.1.1) remains true if we take $\tau^\vee(\sigma)$ to be the
% {\it open} face, there are reasons not to write it that way (see (3.3)).
\smallskip

ii) It is not hard to see that the greatest common quotient of $\overline
\sigma_1$
and $\overline\sigma_2$ is the one-point compactification of $\sigma$.
\enddemo

The picture that follows has been used in [HZ1:\,2.3] and
elsewhere, and will be used in the present work.
\proclaim{Construction 1.2.4} There are quasi-canonical
homeomorphisms of $\d\overline\sigma_2$ and a (cubical) polyhedral
decomposition $P(\widehat \sigma)$ of $\widehat\sigma$ that
subdivides further to the barycentric subdivision of
$\widehat\sigma$.
\endproclaim
To see this, one
notes that
$\d\overline\sigma_2$ is union of the $d$ closed faces that pass through the
origin of
$[0,1]^d$.
% Let $F_j$ be the $j$-th (cubical) face, with equation $t_j=0$.
The barycentric
subdivision of $\widehat\sigma$ consists of taking all
simplices
spanned by the {\it barycenters of} a chain of faces of $\widehat\sigma$:
$$
\tau_1\subset\tau_2\subset ... \subset\tau_m;\tag 1.2.4.1
$$
For the polyhedral complex, one takes only those chains where the codimension
of each face
in the next in (1.2.4.1) is equal to one.  Removing the rest of the faces
yields
the polyhedral
decomposition.  The construction of a homeomorphism can be done recursively,
by doing it
first on all maximal faces, and then adjoining the origin of the cube and the
barycenter of
$\widehat\sigma$.  If one subdivides $\d\overline\sigma_2$ first, one can
treat the
$t_j$'s as ``linear" and make the homeomorphism simplicial
and equivariant
under permutation of the variables.  The specification of the functionals
$q_j$,
which are unique up to positive scalars, determine this mapping uniquely.

The picture looks like a small piece of a Voronoi decomposition (see
[N:\,p.98]).
For a 2-simplex in the boundary when $d=3$, the picture for $P(\widehat
\sigma)$ is:
$$
\tag 1.2.4.2
$$
% *** FIGURE 3
\medskip

{\bf (1.3)} {\bf Discrete quotients and diagonality.}  We next
consider the situation where there is a discrete group action, and
the effect of that on $\L$'s.

Let $D$ be a space on which a discrete group $\Gam$ acts.  Let $E$ be a
partial
compactification of $D$ to which the action of $\Gam$ extends continuously. We
suppose that
the action on $E$ (so also its restriction to $D$) is {\it separated and
discontinuous},
by which we mean that the following two conditions hold:
\smallskip

i) if $y\in E$ is not a $\Gam$-translate of $x$, then there are neighborhoods
$U_x$ of $x$
and $U_y$ of $y$ with $U_y\cap(\Gam\cdot U_x) =\emptyset$.
\smallskip

ii) if $x\in E$, there is a neighborhood $U_x$ of $x$ in $E$ such that if
$\gamma\in\Gam$
and $\gamma\cdot x\in U_x$ then $\gamma\cdot x = x$.
\smallskip

% ii$^\prime$) if $x\in E$, there is a neighborhood $U_x$ of $x$
% in $E$ such that if $1\ne\gamma\in\Gam$ and there is $y\in E$ with
% $y,\,\gamma\cdot y\in U_x$ then $\gamma\cdot x = x$.
% \smallskip

Suppose we have two such $E$---call them $E_1$ and $E_2$. Then $\Gam\times\Gam$
acts on
$E_1\times E_2$.  The diagonal $\Delta(\Gam)$, which we identify with $\Gam$,
is easily seen to
preserve the subset $\L(E_1,E_2)$. The projections of $\L(E_1,E_2)$ onto $E_1$
and $E_2$ are
equivariant for the $\Gam$-actions. Moreover, via either projection, one can
easily deduce
that the action on $\L(E_1,E_2)$ is separated and discontinuous.  We want a
little more.
\demo{{\rm (1.3.1)} Definition} i) We say that the actions of $\Gam$ on
$E_1$ and
$E_2$ are {\it diagonal} when the following holds: if $(e_1,e_2)\in\L(E_1,
E_2)$ and
$(\gamma_1\cdot e_1,\gamma_2\cdot e_2)\in\L(E_1,E_2)$ for some $\gamma_1,
\gamma_2\in\Gam$,
then there exists $\delta\in\Delta(\Gam)$
with $\delta\cdot e_1 = \gamma_1\cdot e_1$ and $\delta\cdot e_2=\gamma_2e_2$;
equivalently,
if $(e_1,e_2)\in\L(E_1,E_2)$ and $(\gamma\cdot e_1,e_2)\in\L(E_1,E_2)$,
there exists
$\delta\in\Delta(\Gam)$ with $e_1=(\delta\gamma)\cdot e_1$ and $e_2=\delta
\cdot e_2$,
i.e., $(\gamma\cdot e_1,e_2) = \delta\cdot(e_1,e_2)$.
\smallskip

ii) We say that the actions of $\Gam$ on $E_1$ and $E_2$ are {\it strongly
diagonal}
when the following holds: if $(e_1,e_2)\in\L(E_1,E_2)$ and $(\gamma\cdot e_1,
e_2)
\in\L(E_1,E_2)$,
then $\gamma\cdot e_1= e_1$ or $\gamma\cdot e_2=e_2$.
\enddemo
\noindent Strongly diagonal actions are diagonal.  We also say that
{\it diagonality holds} under the
conditions of (1.3.1). We prefer strong diagonality in the sequel.
It is easy to see that diagonality is not affected by switching
$E_1$ and $E_2$. We will give a rather simple, and quite
pertinent, example of non-diagonal actions in (1.3.6) below.

\demo{{\rm (1.3.2)} Remark} It is easy to see that if there is a
$\Gam$-equivariant morphism $E_1\to E_2$ in $\PCp(D)$, then
$\L(E_1,E_2)\simeq E_1$ is diagonal.
\enddemo

As the non-trivial issues lie at the boundaries of $E_1$
and $E_2$, we also speak of {\it boundary diagonality}, in the
spirit of the notion of $\Lb(E_1,E_2)$ from (1.1.3). A little
surprisingly, the useful properties of diagonality that we have
found are set-theoretical, i.e., not topological in nature. We
therefore make the following generalization of (1.3.1): \demo{{\rm
(1.3.3)} Definition} Let $E$ and $E'$ be sets given with actions
by a group $\Gam$, so that there is a natural action of
$\Gam\times\Gam$-action on $E\times E'$. Let $B$ be a
$\Delta(\Gam)$-invariant subset of $E\times E'$. We say that $B$
is {\it diagonal} if $(e,e'),\,(\gamma \cdot e,e')\in B$ implies
that $\gamma\cdot e=e$ or $\gamma\cdot e'=e'$. In other words,
{\it on $B$, $\Gam\times\Gam$-equivalence reduces to
$\Delta(\Gam)$-equivalence.}
\enddemo

The following is immediate:
\proclaim{Lemma 1.3.4} i)  Let
$B_1\subset B_2$ be $\Delta(\Gam)$-invariant subsets of $E_1\times
E_2$, with $B_2$ diagonal.  Then $B_1$ is also diagonal.
\smallskip

ii) Let $\Gam'\subset \Gam$. If $B$ is diagonal for $\Gam$, it is diagonal
for $\Gam'$.
\sq
\endproclaim

\noindent In the original setting of $E,\,E'\in\PCp(D)$, we
understand that $B = \L(E,E')$ in (1.3.3) unless stated otherwise.

\demo{{\rm (1.3.5)} Remark} We comment on diagonality in an important setting
related to
$\L$-basechange.  Let $\Gam$ be a group acting
on sets $E_1$,
$E_2$ and $E_3$.  Assume that a $\Gam$-equivariant mapping $E_3\to E_1$ is
specified.
{\it We now take $E_1\times E_2$ to have the diagonal $\Gam$-action.}  Consider
the mappings
$$
E_3\times_{E_1}\!(E_1\times E_2)\subset E_3\times(E_1\times E_2)\to (E_3\times
E_2)\to (E_1\times E_2).\tag 1.3.5.1
$$
The projection $E_3\times (E_1\times E_2)\to E_3\times E_2$ is $(\Gam\times
\Gam)$-equivariant;
however, the subset $E_3\times_{E_1}(E_1\times E_2)$ of the domain is only
$\Delta(\Gam)$-invariant.
It is not true, for instance, that $B$ diagonal in $E_1\times E_2$ implies
$E_3\times_{E_1}\!B$
diagonal in $E_3\times E_2$.
\enddemo
We give next an example of a pair of actions that are non-diagonal
in the sense of (1.3.1), one that fits well with the
considerations of (3.2):

\demo{{\rm (1.3.6)} Example} Let $E_1$ be the real line,
partitioned at the integer points.  This gives a one-dimensional
simplicial complex with 1-simplices $[n,n+1]$ for all $n\in\ZZ$.
The group $\Gam=\ZZ$ acts on $E_1$ in the usual way, by
translation.  Let $E_2$ be the dual cell complex, which is $\RR$
with 1-cells $[n-\tfrac 1 2, n+\tfrac 1 2]$, with the same
$\Gam$-action.  Let
$$
B=\bigcup\,\{\tau\times\tau^\vee:\,\tau=[n,n+1],\,n\in\ZZ\}
$$
(cf.~(1.2.3.1)).  Then the points $(0,\tfrac 12)$ and $(1,\tfrac 12)$ are
in $B$,
though 1 is a
$\Gam$-translate of 0, and $\Gam$ has no fixed points in either $E_1$ or
$E_2$. Note,
however, that $B$ {\it is} diagonal for the action of, say, $2\ZZ$.
\enddemo

The next assertion covers a rather common situation.
\proclaim{Proposition 1.3.7 {\rm (Diagonality and products)}} Let
$D =D_1 \times D_2$, the product of two spaces. Suppose that the
discrete group $\Gam=\Gam_1\rtimes\Gam_2$ (semi-direct product,
$\Gam_1$ normal in $\Gam$) acts factorwise on $D$.  Let
$E=E_1\times E_2$ and $E'=E'_1\times E'_2$, for
$\Gam_1$-equivariant $E_1,\,E'_1\in\PCp(D_1)$, and
$\Gam_2$-equivariant $E_2,\,E'_2\in\PCp(D_2)$. Assume that
diagonality holds for the actions of $\Gam_1$ on $E_1$ and $E'_1$,
and for the actions of $\Gam_2$ on $E_2$ and $E'_2$. Then the
actions of $\Gam$ on $E$ and $E'$ are diagonal, and $\L(\Gam \back
E,\Gam\back E')\simeq\Gam\back \L(E,E')$ is a $\Gam_2\back\L(E_2,
E'_2)$-fibration over $\Gam_1\backslash\L(E_1,E'_1)$.
\endproclaim
\demo{Proof} We begin by taking the $\Gam_2$-quotients.  That gives at once
(one can use Lemma 1.1.2):
$$\multline
\L(E_1\times(\Gam_2\back E_2),E'_1\times(\Gam_2\back E'_2))\simeq
\L(E_1,E'_1)\times\L((\Gam_2\back E_2),(\Gam_2\back E'_2))\\
\simeq\Gam_2\back\L(E,E').\endmultline
$$
Taking the quotient by $\Gam/\Gam_2\simeq \Gam_1$ (isomorphism of groups),
we obtain the result.\sq
\enddemo

The following extension of Proposition 1.1.16 combines the above
notions, and is tailored to our later needs in \S3: \proclaim
{Proposition 1.3.8}  Let $D$ be a space with a separated
discontinuous action of a discrete group $\Gam$, and let
$X=\Gamma\backslash D$.  Suppose that $E_1$ and $E_2$ are partial
compactifications of $D$, yielding spaces for which the $\Gamma$
actions are diagonal. Assume that there is a subset $S$ of $D$
that maps onto $X$, and whose respective closures $S_1$, $S_2$ in
$E_1$ and $E_2$ are compact; thus $Y_1 = \Gamma \backslash E_1$
and $Y_2 = \Gamma\backslash E_2$ are compactifications of $X$.
Then the canonical mapping
$$
\Gamma\backslash\L(E_1,E_2)\to \L(Y_1,Y_2)
$$
(where $\Gam$ acts on $\L(E_1,E_2)$ as $\Delta(\Gam)$) is a homeomorphism.
\endproclaim
\demo{Proof} Let $\{x_j\}$ be a sequence in $X$ that converges in both $Y_1$
and $Y_2$.
The existence of $S$ allows the lifting of $\{x_j\}$ to $S$, such that a
subsequence
converges in $S_1$ and $S_2$, a fortiori in $E_1$ and $E_2$. Thus $\L(E_1,E_2)
\to\L(Y_1,Y_2)$
is surjective.  We wish to show that the fibers are $\Delta(\Gam)$-equivalence
classes.
So, suppose that $(e_1,e_2)$ and $(e'_1,e'_2)$ have the same image in
$\L(Y_1,Y_2)$.
Then $e'_1=\gamma_1\cdot e_1$ and $e'_2=\gamma_2\cdot e_2$ for $\gamma_1,
\gamma_2\in\Gam$.
By diagonality, we can take $\gamma_1=\gamma_2$, and we are done.\sq
\enddemo
\bigskip

\centerline{\bf 2. Locally-symmetric varieties and their compactifications}
\medskip

In this section we present key elements of the structure of the relevant
compactifications
of locally symmetric varieties.  (A list of notational changes from [Z4] is
provided
at the end of the introduction.)
\medskip

{\bf (2.1)} {\bf Boundary components and the Baily-Borel
compactification.} Let $G$ be a semi-simple algebraic group
defined over the rational number field $\QQ$. For $P$ a parabolic
subgroup of $G/\QQ$ ($\QQ$-{\it parabolic subgroup}; we will
henceforth suppress the ``$\QQ$''), let $W_P$ denote the unipotent
radical of $P$.  The {\it Levi quotient} of $P$ is $\Pred =
P/W_P$, and the projection $P\to \Pred$ is split by a choice of
Levi subgroup $L_P\subset P$.

Let $D$ be the symmetric space of non-compact type
% (possibly with Euclidean factors)
associated to $G(\RR)$.  For each parabolic subgroup $P$, its real points act
transitively on $D$.
Unless specified otherwise, {\it we assume throughout that $G(\RR)$ is of
Hermitian
type}, i.e., that $D$ has a $G(\RR)$-invariant complex structure.

We further assume that $G$ is (almost) simple over $\QQ$. Then for
each maximal parabolic $P$,
$L_P$ is an almost-direct product of two groups, commonly denoted $\GhP$ and
$\GlP$ after [AMRT:\,p.209].
The rational boundary component $D_P$ of $D$ that is normalized by $P$ is
biholomorphic to
the Hermitian
symmetric space associated to $\GhP$. On the other hand, the group $\GlP$
(rarely of
Hermitian type) acts trivially on $D_P$, as does $W_P$.  The set
$$
D^* = D\sqcup\bigsqcup\{D_P: P \text{ is maximal parabolic}\}\tag 2.1.1
$$
has a canonical $G(\QQ)$-action.  It is given the $G(\QQ)$-equivariant {\it
Satake topology}
[S,\S 2] (see [Z3:(3.9)]).

The $\QQ$-root system of $G$ is of classification type $BC$ or $C$ [BB:2.9],
whose Dynkin
diagram is a linear graph with distinguished end.  The set $\Delta$ of simple
roots is
thus totally ordered by specifying the root at the distinguished end to be
the minimal
element.  In a standard lattice of parabolic subgroups, the specification of a
maximal
parabolic subgroup $P$ is equivalent to selecting one simple
root $\beta_P\in\Delta$;
the correspondence is: the root space for a root $\beta$ is contained in $W_P$
if and only
if $\beta_P$ occurs with positive coefficient in the expansion of $\beta$. One
says that $Q\prec P$ whenever $\beta_{Q}\prec
\beta_P$.
Moreover, the Dynkin diagram for $\GhP$ is the segment $\Delta_{h,P}$ of type
(B)C below
$\beta_P$
(in the total ordering), and that of $\GlP$ is the segment $\Delta_{\ell,P}$
of type A
above $\beta_P$.  Thus, one can assert:
\proclaim{Proposition 2.1.2}
When $G$ is simple over $\QQ$, the following are equivalent for maximal
parabolic
subgroups: (i) $Q\prec P$,\, (ii) $D_{Q}$ is a boundary component of $D_P$,\,
(iii) $G_{h,Q}\subset \GhP$\,,\, (iv)
$G_{\ell,Q}\supset \GlP$\,.\sq
\endproclaim

Let $\Gam \subset G(\QQ)$ be an arithmetic subgroup. We assume throughout that
$\Gam$ is
neat, i.e., that for every algebraic subquotient group $H$ of $G$, the induced
subquotient
of $\Gam$, denoted $\Gam (H)$ (instead of the ``$\Gam_H$" that appears
frequently
in the literature), is torsion-free. The latter for $H=G$ already gives that
the quotient space $X =\Gam
\backslash D$,
called a {\it locally-symmetric} (or {\it arithmetic}) {\it
variety}, is non-singular. We are interested in some of its
compactifications.

The {\it Baily-Borel Satake compactification} of $X$ is
$$
X^*=\Gam\backslash D^*,\tag 2.1.3
$$
homeomorphic in the sense of (1.1) to the Satake compactification of $X$ that
is parametrized
by the root at the distinguished end (see [Z3:(3.9)]).  It is shown in [BB]
that $X^*$
is a projective algebraic variety (of which $X$ is the regular locus), whence
the word ``variety" in the name for
$X$.
The decomposition (2.1.1) induces a stratification of (2.1.3) as follows.
For $P$ as
above, a boundary stratum
$$
M_P=\Gamma (P)\backslash D_P\tag 2.1.4
$$
of $X^*$ is induced by the inclusion $D_P\subset D^*$ in (2.1.1).  For
$Q\prec P$,
$M_{Q}$ lies in the
closure of $M_P$ in $X^*$,
the latter giving the Baily-Borel compactification $M^*_P$ of $M_P$.  The
quotient by
$\Gam$ identifies the boundary components of $\Gam$-conjugate $P$'s, so
boundary
strata in $X^*$ are parametrized by {\it $\Gam$-conjugacy classes} of maximal
parabolic
subgroups, and these are finite in number.
\bigskip

{\bf (2.2)} {\bf Toroidal compactifications and their quotients}.
The toroidal compactifications $\XSig$ of [AMRT] are predicated on
the notion that they will map to the Baily-Borel compactification
$X^*$.  As the construction is rather complicated, we give here
only a brief description, referring to the literature for details.

For $P$ maximal parabolic, the center $U_P$ of $W_P$ plays a key role in
the construction;
the quotient $V_P = W_P/U_P$ is commutative and even-dimensional, as follows
from the
root structure.  It is elementary but fundamental that
$U_Q
\supset U_P$ whenever $Q\prec P$ (there is no corresponding assertion for the
whole unipotent
radical).

For each $P$, there is a tower of quotient mappings of groups
$$
P'=G_{h,P}(\RR)W_P(\RR)U_P(\CC)\to G_{h,P}(\RR)V_P(\RR)\to G_{h,P}(\RR).
\tag 2.2.1
$$
These act homogeneously on a tower of spaces $\check D(P)\to D^A_P\to D_P$,
with
common isotropy subgroup
$K_{h,P}$,
which is maximal compact in $G_{h,P}(\RR)$. From (2.2.1) and the Cayley
transform
associated to $P$ comes the Siegel domain picture of $D$ relative to $D_P$ (see
[AMRT:III,\S 4],
and also our (4.1)).
Then,
dividing out the respective actions of $\Gamma(G_{h,P}W_P)$, $\Gamma(G_{h,P}
V_P)$, and
$\Gamma(G_{h,P})$, one gets the basic tower (essentially of mixed Shimura
varieties;
see [HZ1:1.6]):
$$
M'_P @>\pi_2 >> \Cal A_P @>\pi_1 >> M_P.\tag 2.2.2
$$
(The use of ``prime" in $M'$ and $P'$ differs from convention that begins in
(2.3).)

We underscore the absence of $\GlP$ in (2.2.1).  (It disappears through
inclusion
in the isotropy subgroup for a suitable basepoint at infinity.)  The action of
$\Gam(\GlP)$ on
(2.2.2) is
induced by conjugation on the groups in (2.2.1).  It is trivial on $M_P$ and
free on $M'_P$.
As $\pi_1$ is a fibration by abelian varieties,
in particular
proper, the $\Gam (G_{\ell,P})$-orbits in $\Cal A_P$ are rather ugly,
except
in the case that $\Cal A_P = M_P$ (i.e., when $W_P$ is commutative).  It is an
essential
observation that $\pi_2$ is a principal torus bundle with fiber
$$
T_P = \Gam (U_P)\back U_P(\CC).\tag 2.2.3
$$

Through (2.2.3), one uses torus embeddings to attach
a substantial boundary to $M'_P$ for each $P$, in a fashion that is both
equivariant
for the action of
$\Gam(G_{\ell,P})$ and compatible with order relations and
conjugacy among maximal parabolic subgroups [AMRT:\,p.117].
One then uses reduction theory to do the
same for $X=\Gam\backslash D$, and this produces the
desired compactification of $X$ [AMRT:\,III,\S\S5,6].

We want to specify, and comment on, some of the details.  The
space $\XSig$ depends on a collection $\Sigma$ of simplicial cone
complexes ({\it fans}), $\Sigma_P\subset U_P(\RR)$ (as $P$
varies), such that for ${}^\gamma\! P=\gamma P \gamma^{-1}$ (with
$\gamma \in\Gam$),
$\Sigma_{{}^\gamma\!P}=\gamma(\Sigma_P)\gamma^\-$; also, whenever
$Q\prec P$, $\Sigma_Q$ contains $\Sigma_P$ as a closed boundary
stratum (in a precise sense; see below), and $\Gam(G_{\ell,P})$
acts separated and discontinuously on $\Sigma_P$, compatibly with
the order relation.  We always assume that $\Sigma_P$ is
sufficiently fine so as to have {\it full boundary}, in the sense
of [HZ1:2.2.6], to avert combinatorial anomalies.

The use of the term {\it cone complex} in the above means that a fan $\Upsilon$
is
subject to the following two axioms:
\smallskip

A1. {\it If $\sigma\in\Upsilon$, and $\tau$ is a face of $\sigma$, then
$\tau\in
\Upsilon$.}

A2. {\it If $\sigma,\sigma'\in\Upsilon$, and $\sigma\cap\sigma'\ne\{0\}$, then
$\sigma\cap\sigma'$
is a face of both $\sigma$ and $\sigma'$.}

\noindent Any fan $\Upsilon$ in
$U_P(\RR)$ determines in a standard way a torus embedding $T_P\subset T_{
P,\Upsilon}$,
on which $T_P$ acts (see [HZ1:\,1.3] and references cited therein).  Briefly,
this goes
as follows; we present it for an arbitrary torus $T$.
Each cone
$\sigma$ determines an affine torus embedding $T=T_{\{0\}}\subset T_\sigma$
($T$-equivariant
partial compactification),
such that
$T_\tau\subset T_\sigma$ whenever $\tau$ is a face of $\sigma$.  The axioms
allow for
gluing, viz., if $\tau = \sigma\cap\sigma'$,
$$
T_{\{\sigma,\sigma'\}} = T_\sigma\cup_{T_\tau}T_\sigma'.\tag 2.2.4
$$
One takes
$T_{\Upsilon} = \,\bigcup\,\{T_\sigma\, :\, \sigma\in\Upsilon\}$,
with identifications as given by (2.2.4).

The notation in (2.2.4) contains the incidental statement that $T_{\{\sigma,
\sigma'\}}$
is determined only by $\sigma$ and $\sigma'$ (and their intersection).
Moreover,
nothing
new is obtained when $\sigma'$ is a face of $\sigma$; $\sigma$ and the fan it
generates
(i.e., the one consisting of $\sigma$ and all of its faces) produce the same
torus embedding.
Thus, axiom A1, which provides for gluing, can be relaxed.  We say that
a collection of
cones $\Upsilon$ is a {\it loose fan} if it satisfies the following weakened
version
of the axioms of a fan:
\smallskip

A$2'$. {\it If $\sigma,\sigma'\in\Upsilon$, then $\sigma\cap\sigma'$ is in
$\Upsilon$
(so is a face of both $\sigma$ and $\sigma'$).}
\smallskip

\noindent A loose fan $\Upsilon$ generates a fan $\overline\Upsilon$, given as
the set
of all elements of $\Upsilon$ and their faces. It also determines a torus
embedding
$T_\Upsilon$, by the same procedure that is specified for fans.  It is clear
that
$T_\Upsilon = T_{\overline\Upsilon}$.

It is convenient to remove the origin and contract out the cone
dilations in a fan, yielding from $\Upsilon$ a simplicial complex
$\widehat\Upsilon$; likewise, we will take off our hats to get the
fans, i.e., pass routinely from $\widehat\Upsilon$ to $\Upsilon$
when $\widehat\Upsilon$ is mentioned first.  Moreover, we standardize
notation by using $\sigma$, $\tau$, etc.~for simplicial cones, and
$\widehat\sigma$, $\widehat\tau$, etc.~for simplices. Each cone
$\sigma\in\Upsilon$ corresponds to a $T$-orbit, which we denote
$T(\sigma)$, with $\roman{codim}_ \CC\, T(
\sigma)=\dim_\RR\sigma$, and then
$$
T_\sigma = \bigsqcup\,\{ T(\tau): \tau\text{ is a face of
$\sigma$}\}; \tag 2.2.4.1
$$
here we allow the improper face ($\tau = \sigma$).  We note that $\widehat
\Upsilon$
is isomorphic to the nerve of the set of closed (i.e., closures of) $T$-orbits
in $T_\Upsilon$,
via the correspondence $\widehat\tau\mapsto\tau^\vee$ (cf.~(1.2); see also
(3.2.4.1) below).

With $T = T_P$ and $\Upsilon = \Sigma_P$, the above construction
determines a partial compactification $M'_{P,\Upsilon}$ of $M'_P$,
with
$$
M'_{P,\Upsilon} = M'_P\times^{T_P}T_{P,\Upsilon}\to\Cal A_P\to
M_P.\tag 2.2.5.1
$$
The mapping onto $\Cal A_P$ extends $\pi_1$ in (2.2.2); the
restriction of (2.2.5.1) to the boundary is
$$
\d M'_{P,\Upsilon}=:\widetilde Z_{P,\Upsilon}=M'_P\times^{T_P}\d
T_{P,\Upsilon}\to\Cal A_P\to M_P.\tag 2.2.5.2
$$
We refer to the preceding as the {\it toroidal construction from
$\Upsilon$} (and $T_P$ and $M'_P$). Then $\widehat \Upsilon$ is
isomorphic to the nerve of the set of all irreducible components
of $\widetilde Z_{P,\Upsilon}$, as it is already true for
$T_\Upsilon$.  Because we are taking only simplicial fans $\d
T_{P,\Upsilon}$, hence also $\widetilde Z_{P,\Upsilon}$, is a
union of smooth divisors with normal crossings.

Given a subset
$\Sigma$
of the fan $\Upsilon$, we define a subset $T(\Sigma)$ of $T_\Upsilon$ by:
$$
T(\Sigma) = \bigsqcup\,\{T(\sigma):\sigma\in\Sigma\};\tag 2.2.6
$$
note that with this notation, $T(\Upsilon)= T_\Upsilon$. The space $T(\Sigma)$
is invariant
under $T$, and it therefore determines a subset of $\widetilde Z_{P,\Upsilon}$:
$$
\widetilde Z(\Sigma) = M'_P\times^T \d T(\Sigma).\tag 2.2.7
$$
We allow ourselves to write $T(\widehat\Sigma)$ and $\widetilde Z(
\widehat\Sigma)$, and the stratification describe

instead of $T(\Sigma)$ and $\widetilde Z(\Sigma)$ resp.

\demo{{\rm (2.2.8)} Remark} There is a familiar and simple characterization of
when
$\d T(\widehat\Sigma)$ is closed in $\d T(\widehat\Upsilon)$.  It involves the
{\it star}
of $\widehat\Sigma$ in $\widehat\Upsilon$, given as
$$
\roman{Star}\,(\widehat\Sigma)=\{\tau\in\widehat\Upsilon:
\text{ a
vertex of $\widehat\tau$ is contained in }\widehat\Sigma\}.
$$
Then the closure of $\d T(\widehat\Sigma)$ in $\d T(\widehat\Upsilon)$ is $\d
T(\roman{Star}
(\widehat\Sigma))$; thus, $\d T(\widehat\Sigma)$ is closed if and only if the
inclusion
$\widehat\Sigma\subseteq\roman{Star}\,(\widehat\Sigma)$ is an equality. In
particular,
when $\Sigma$ is a proper subfan of $\Upsilon$, $\d T(\widehat\Sigma)$ is
never
closed. (See also (2.2.10) below.)
\enddemo

The interior of $\Sigma_P$ is, in fact, $\Gam(\GlP)$-equivariantly homeomorphic
to the
homogeneous cone $C_P$ occurring in the Siegel domain mentioned in conjunction
with (2.2.2),
an orbit of the adjoint representation of $\GlP$ on $U_P(\RR)$,
which is how it enters the construction.
Let $D_{\ell,P}$ be the space of type $S-\QQ$ (in the sense of [BS:\,2.3]) for
$\widehat
G_{\ell,P}=G_{\ell,P}/A_P$, with $A_P$ as in (2.3) below ($D_{\ell,P}$
need not
be a symmetric space, as $\widehat G_{\ell,P}$ may contain central anisotropic
tori).
Put $X_{\ell,P}=\Gamma(G_{\ell,P})\backslash D_{\ell,P}$ and $\widehat
\Sigma'_P =
\Gamma(G_{\ell,P})\backslash \widehat\Sigma_P$.
The complex
$\widehat\Sigma'_P$ is actually a stratified triangulation of a Satake
compactification
of $X_{\ell,P}$ (in the sense of [Z3:\,\S3] and [HZ2:\,(2.1)]).  As such,
$\widehat\Sigma'_P$
is a non-Hermitian analogue of $M_P^*$, which is itself a Satake
compactification
of $M_P$ (see [Z3:(3.9)]).

We mention two good choices for the fan $\Upsilon$: $\Sigma_P$
and its
subfan $\Sigma_P^c$; between them lies the set $\Sigma_P^\circ$. All three are
closed under
the action of $\Gam(\GlP)$ on $\Sigma_P$.  We give the definition of the
latter two:
$\Sigma_P^c$ is the fan spanned by the interior edges of $\Sigma_P$,
and
$\Sigma_P^\circ$ consists of those cones in $\Sigma_P$ that contain at least
one edge
from $\Sigma_P^c$. Thus, $\Sigma_P^\circ$ is usually not even a loose fan.  In
terms of simplices,
$$
\widehat\Sigma_P^\circ=\widehat\Sigma_P-\bigcup\,\{\widehat\Sig_Q:Q\succ P\},
\quad\text{so}\quad
\widehat\Sigma_P=\bigsqcup\,\{\widehat\Sigma_Q^\circ:Q\succeq P\}.\tag 2.2.9
$$
\demo{{\rm (2.2.10)} Remark}
We illustrate $\widehat\Sigma_P^c$ and $\widehat\Sigma_P^\circ$ by describing
them in
the case where $\widehat\Sigma_P$ is replaced by the fan generated by a
one-simplex,
which consists of two vertices [0] and [1] spanning the simplex
[0,1].
We declare that [1] is interior and [0] is boundary (if this were contained in
$\widehat\Sigma_P$,
then there would be $Q\succ P$ with $[0]\in\widehat\Sigma_Q$). Then:
$$
\widehat\Sigma_P^c = \{[1]\},\quad \widehat\Sigma_P^\circ =\{[1],[0,1]\},\quad
\text{and}\quad\widehat\Sigma_P=\{[0],[1],[0,1]\}.
$$
\enddemo

The boundary of $X^\t$ stratifies naturally in terms of $\Sigma$.  We admit
with some
apology that the notion of {\it boundary stratum} changed during the
progression
from [HZ1] to [HZ2].  Also, the term ``stratum" is being used
loosely,
in that our strata need not be smooth.

The {\it Baily-Borel-type $P$-stratum}
in $\XSig$,
denoted ${}^<\!Z_P$ (cf.~[HZ1:\,1.5]), is the $\Gam(G_{\ell,P})$-quotient of
${}^<\!
\widetilde Z_P  :=\widetilde Z(\Sigma_P^\circ)$; this in fact equals the
closure
of $\widetilde Z(\Sigma_P^c)$ in $\widetilde Z(\Sigma_P)$.
It follows
directly from the construction that
% , as a set,
$$
\widetilde Z(\Sigma_P) = \bigsqcup\, \{{}^<\!\widetilde Z_Q:Q\succeq P\}
$$
(for such $Q$, $\Sigma_Q$ is a boundary component of $\Sigma_P$).  In fact,
${}^<\!\widetilde Z_Q$
is produced
in the toroidal construction from $\widehat\Sigma_Q$ (and $T_Q$ and
$M'_Q$),
and is partially compactified in the toroidal construction from $\widehat
\Sigma_P
\supset\widehat\Sigma_Q$.  We
let $Z_P$
denote the closure of ${}^<\!Z_P$ in $\XSig$, and refer to it as the {\it
closed
$P$-stratum}.  The complement of ${}^<\!Z_P$ in $Z_P$ is the (disjoint) union
of all
$ {}^<\!Z_Q\cap Z_P$
with $Q\prec P$---now with $\widehat\Sigma_P$ a boundary component of $\widehat
\Sigma_Q$.
Indeed, the combinatorial data for constructing $Z_P$
is contained
in the union of the {\it open} regular neighborhoods of $\widehat\Sigma_P^c$
in $\widehat
\Sigma_Q$, as $Q\preceq P$ varies (see~[HZ2:(2.5.1,\,i)]); when $Q=P$, this
neighborhood
coincides with
$\widehat\Sigma_P^\circ$.  Under the natural mapping $\pi:\XSig\to X^*$
[AMRT:\,p.254],
$\pi^{-1}(M_P) = {}^<\!Z_P$, and then $\d X^\t =\bigsqcup_P{}^<\!Z_P$, with
the union taken
over $\Gam$-conjugacy classes.

Next, let $R$ be an arbitrary parabolic subgroup (possibly maximal) of $G$.
We want to
specify the (cubical) {\it $R$-stratum} $\overset\circ\to Z_R$ in $\XSig$. One
can write
$R$ uniquely as an intersection of maximal parabolic subgroups in a standard
lattice:
$$
R = \bigcap\,\{Q: Q\in\Cal S\},\quad\text{ where}\tag 2.2.11
$$
$$
\Cal S=\Cal S_R =\{Q\text{ maximal}: Q\supseteq R\}\leftrightarrow \{\beta\in
\Delta
\,: \beta=\beta_Q\text{ for some maximal }Q\supseteq R\}.
$$
(Note that this usage of subsets of $\Delta$ is complementary to the
usual one
for indexing parabolic subgroups; a minimal parabolic subgroup corresponds to
$\Delta$ here,
not $\emptyset$.)  Let $P$ be the smallest element in $\Cal S$.
As in
[HZ2:\,(2.2)], we then say that $R$ is {\it subordinate to} $P$, and write
$P=\Pi(R)$.
The $R$ that are subordinate to a given $P$ are in canonical one-to-one
correspondence
with the standard parabolic
subgroups $G_{\ell,R}$ of $G_{\ell,P}$ relative to $\Delta_{\ell,P}$,
including the
improper one, by taking
$$
G_{\ell,R}= R\cap G_{\ell,P},\quad\text{so then}\quad R = G_{\ell,R}\cdot\GhP
\cdot W_P.\tag 2.2.12
$$

{\it If} there were a good way to remove the $\Gam$-quotient from $X^\t$ to
yield a
space we might call $\widetilde X^\t$, the closed strata $Z_P$ of $X^\t$
would be
induced from their inverse images $\widetilde Z_P$ in $\widetilde
X^\t$. The
stratification would be {\it cubical} (or {\it corner-like}), in that the
closed
strata $\{\widetilde Z_R\}$ of $\widetilde X^\t$ would satisfy
$$
\widetilde Z_R = \bigcap\{\widetilde Z_Q:Q\in \Cal S\}\tag 2.2.13
$$
(compare (2.3.3) below), so the strata would be deducible from (2.2.13).
However,
the preceding is an oversimplification.

Here is one way to proceed.  Define for each $R$ with $\Pi(R)=P$ a subset
of $\widehat\Sigma_P$:
$$
\widehat\Sigma_R^c := \{\widehat\tau\in\widehat\Sigma_P\,|\,\widehat
\tau\text{ has
all vertices in $\bigcup_{Q\in\Cal S}\,\widehat\Sigma^c_Q$ and a
vertex in each $\widehat\Sigma^c_Q$}\}.\tag 2.2.14
$$
% $$
% \widehat\Sigma_R^\circ := \{\widehat\tau\in\widehat\Sigma_P\,|\,\widehat\tau
% \text{ has at least one
% vertex in }\widehat\Sigma^c_Q\text{ for each }Q\in\Cal S\}.\tag 2.2.15
% $$
This agrees with what we gave earlier when $R$ is maximal. As $R$ varies
(subordinate
to $P$), the $\widehat\Sigma_R^c$'s are disjoint.
\demo{{\rm (2.2.15)} Definition}
Let $\check Z_R = \widetilde Z(\widehat\Sigma_R^c)$. The {\it $R$-stratum} of
$\d X^\t$ is
$\overset\circ\to Z_R\simeq\Gamma(G_{\ell,R})\backslash\check Z_R$.
\enddemo
\noindent Note that replacing $R$ by any $\Gam$-conjugate of $R$ gives the
same stratum;
they are otherwise disjoint.
We denote by $Z_R$ the closure of $\overset\circ\to Z_R$ in $\d X^\t$ (the
closed $R$-stratum).
Because $Z_R$ contains points from more than one standard lattice, it is
somewhat
inconvenient to express it in terms of torus orbits.

We note that one has for the Baily-Borel-type $P$-stratum
$$
{}^<\!Z_P =\,\bigsqcup\,\{\overset\circ\to Z_{R'}:\Pi(R')=P\}.\tag 2.2.16.1
$$
This suggests that one might also put, for general $R$,
$$
{}^<\!Z_R =\,\bigsqcup\,\{\overset\circ\to Z_{R'}:R'\subseteq R,\,\Pi(R')=P\},
\tag 2.2.16.2
$$
which gives the portion of $Z_R$ that is created in the toroidal construction
from
$\widehat\Sigma_P$.

%%  Note that
%%  $$
%%  \widehat\Sigma_R^\circ=\bigsqcup\,\{\widehat\Sigma_{R'}^c: R'\subseteq R\}.
%%  $$

% On the other hand, in terms of
% $T_P$-orbits,
% we have ${}^<\!\widetilde Z_R = Z(\widehat\Sigma_R^\circ)$, where
% $$
% \widehat\Sigma_R^\circ := \{\widehat\tau\in\widehat\Sigma_P\,|\,\widehat\tau
% \text{ has at least one
% vertex in }\widehat\Sigma^c_Q\text{ for each }Q\in\Cal S\}.\tag 2.2.12
% $$
% (This agrees with what we gave earlier if $R$ is maximal.)
% Let ${}^<\!Z_R$ denote the image of ${}^<\!\widetilde Z_R$ in $\XSig$,
% and  $Z_R$ the
% closure of ${}^<\!Z_R$ in $\XSig$.
% The latter is the {\it closed} $R$-stratum of $\XSig$.  The set ${}^\circ\!
% \widetilde Z_R$
% is now given by $Z(\widehat\Sigma_R^c)$, where
% $$
% \widehat\Sigma_R^c := \{\widehat\tau\in\widehat\Sigma_R^\circ\,|\,\widehat
% \tau
% \text{ has no vertices
% outside }\bigcup_Q\,\{\widehat\Sigma^c_Q : Q\in\Cal S\}\}.\tag
% $$

%% Note that
%% ${}^<\!Z_R = {}^<\!Z_P\cap Z_R$,
%% and ${}^<\!Z_R\simeq\Gamma(G_{\ell,R})\backslash {}^<\!\widetilde Z_R$.

\demo{{\rm (2.2.17)} Remark} $Z_R$ is generally only a connected component of
$\bigcap\,
\{Z_Q: Q\in \Cal S\}$.  Because of the $\Gam$-quotient, this intersection
usually has more
than one
connected component, only one of which one wants to associate with $R$. The
others
correspond to certain $G(\QQ)$-conjugates of $R$ (see [HZ2:(3.5,\,App.)]).
\enddemo

In [HZ2:(1.4,\,(d))], we introduced a canonical (topological)
quotient of $\XSig$, called its {\it real boundary quotient}.
Specifically, one collapses ${}^<\!Z_P$ (2.2.16.1) of $\XSig$ to
${}^<\!Z_P/T_P^c$, which inherits the fibrations (2.2.5.2) over
$\Cal A_P$ and $M_P$; here, $T_P^c$ denotes the compact factor of
$T_P$, viz., $\Gam (U_P)\back U_P(\RR)$. Because its construction
is parallel to that of the excentric Borel-Serre compactification
[HZ2:(1.4,\,(b))] (see (2.3) below), we rename it an {\it excentric
toroidal compactification}
of $X$, and denote it $\XSigexc$. We get a tower of
compactifications (for each $\Sig$):
$$
\XSig\to\XSigexc\to X^*.\tag 2.2.18
$$
Likewise, we write ${}^<\!Z_P^\e$ for ${}^<\!Z_P/T_P^c$ and
$\overset\circ\to Z_R^{\text{exc}}$ for the $R$-stratum of
$\XSigexc$, viz.~$\overset\circ\to Z_R/T_P^c$; define $\check
Z_R^{\text{exc}}$ analogously as $\check Z_R/T_P^c$.  We can see
that the topological structure of $\check Z_R^{\text{ exc}}$ is
rather simple.  First, we make a basic observation:
\proclaim{Proposition 2.2.19 {\rm [HZ1:\,2.1]}}
Let $T\simeq
(\CC^*)^n$ be a torus, so that $T/T^c \simeq (\RR^+)^n$.  Given
any simplicial fan $\Upsilon$ for $T$, let $T_\Upsilon$ denote the
corresponding torus embedding. Then $T_\Upsilon/T^c$ is a
contractible $n$-dimensional manifold-with-corners. \sq
\endproclaim
From this, we see immediately that
\proclaim{Corollary 2.2.20}
The natural mapping $\check Z_R^{\text{exc}}\to\Cal A_P$ is a homotopy
equivalence.\sq
\endproclaim
\bigskip

{\bf (2.3)} {\bf The Borel-Serre compactification and its
quotients.} {\sl The spaces $\Xbar$ and $\Xbarred$ below are
defined the same way in the absence of Hermitian structure.}

The reader may assume again, for simplicity only, that $G$ is simple over
$\QQ$,
for the general group is almost a product of these.  Let $R$ be any
parabolic
subgroup (not necessarily maximal), and $S_R$ the maximal $\QQ$-split torus
of
the center of $R^{\text{red}}$.  Then $A_R$, the identity component of
$S_R(\RR)$, acts
{\it geodesically}
on $D$ [BS:\S 3]. (In the case where $G=SL(2)$, $D$ is the upper half-plane,
and $P$
is the group of upper-triangular matrices, $A_P$ is the the group of diagonal
matrices,
and the $A_P$-orbits are the vertical lines.) When we write $R$ as in (2.2.11),
we have
$$
A_R = \prod_{P\in \Cal S}A_P.\tag 2.3.1
$$
The simple roots in $\Cal S$ set up an isomorphism $A_R\simeq (0,\infty)^{
\Cal S}$,
which one uses to define $\overline A_R$ as $(0,\infty]^{\Cal S}$, and then
the corner
$D(R) = D\times^{A_R} \overline A_R$.  The Borel-Serre boundary face
associated to $R$ is
$$
e(R)\simeq D\times^{A_R}\{(\infty,...,\infty)\}\subset D(R);
$$
it decomposes (as a manifold) as
$$
e(R)\simeq D/A_R\simeq e(R)^{\text{red}}\times W_R(\RR);\quad\text{also}\quad
D\simeq  e(R)\times A_R.
\tag 2.3.2
$$
Here $e(R)^{\text{red}}$ is the ``symmetric space" (with Euclidean factors
allowed)
of $L_R/A_R$ (cf.~$D_{\ell,P}$ in (2.2)).
It is easy to see that $D(R)$ contains $e(Q)$ for all $Q\supset R$.  Then
with a natural topology,
$$
\overline D = \bigcup\,\{D(R): R\text{ is parabolic in } G\} =
D\sqcup\bigsqcup\,\{e(R): R\text{ is parabolic in } G\}
$$
is a manifold-with-corners with the $e(R)$'s as its open faces
[BS:\S 7].
The closure $\overline{e(R)}$ of $e(R)$ in $\overline D$ adjoins to $e(R)$
exactly those
$e(Q)$ with $Q\subseteq R$; indeed, the Borel-Serre construction actually
applies to
$e(R)$, and then $\overline{e(R)}$ is {\it its} Borel-Serre compactification.
We have
$$
\overline{e(R)} = \bigcap\, \{\overline{e(P)}:P\in\Cal S\}.\tag 2.3.3
$$

There is an evident $G(\QQ)$-action on $\overline D$, such that for $g\in
G(\QQ)$,
$g\cdot e(R) = e(gRg^\-)$.  Given any neat
arithmetic
subgroup $\Gam$ of $G(\QQ)$, one puts $\overline X = \Gam\backslash\overline
D$; this
is the {\it Borel-Serre compactification} of $X$.  It is a compact
manifold-with-corners
with open faces of the form
$$
e'(R)\simeq \Gam(R)\backslash e(R),
$$
and these are parametrized
by $\Gam$-conjugacy
classes of parabolic subgroups $R$ of $G$.  From (2.3.2), one
sees that
$e'(R)$ is a fiber bundle over $e'(R)^{\text{red}}=: \Gam(R)\backslash e(R)
^{\text{red}}$
with fiber $\Gam(W_R)\backslash W_R(\RR)$, a compact nilmanifold.  As was
the case
for the toroidal boundary in (2.2),
$$
\bigcap\, \{\overline{e'(P)}: P\in\Cal S\}
$$
has the same finite number of connected components (arising for the same
reason as (2.2.17)), one of which is
$\overline{e'(R)}$.

We will be working with two natural quotients of $\Xbar$,
$\Xbar^{\text{exc}}$ and $\Xbarred$, the {\it excentric} and {\it
reductive} Borel-Serre compactifications of $X$ resp. (The space $\Xbarred$
played a central role in our previous article [Z4], where it was
denoted $M^{RBS}$.) They are determined by the respective
quotients of the boundary strata:
$$
e(R)\to U_P(\RR)\backslash e(R)=: e(R)^{\text{exc}}\to e(R)^{\text{red}},
\tag 2.3.4
$$
when $R$ is subordinate to the maximal parabolic $P$.  As such, all three
spaces
have the same number of boundary strata, indeed the same number as $X^\t$.
They
fit into a tower of compactifications:
$$
\Xbar\to\Xbar^{\text{exc}}\to\Xbarred\to X^*.\tag 2.3.5
$$

That the Baily-Borel compactification $X^*$ is a quotient of $\Xbarred$ is
shown in
[Z3:(3.11)] and [GT:2.6.3]. The quotient mapping $\Xbarred\to X^*$ is given
by the following.
Suppose that $R$ is subordinate to the maximal parabolic $P$.  Then there
is a product decomposition
$$
R^{\text{red}} = (G_{\ell,R})^{\text{red}}\times\GhP,\tag 2.3.6
$$
with $G_{\ell,R}$ as in (2.2.12).  This induces a projection
$$
e(R)^{\text{red}}\to D_P, \qquad \text {and thereby} \qquad e'(R)^{\text{red}}
\to M_P.\tag 2.3.7
$$
Letting $R$ vary over parabolics subordinate to $P$, one gets (cf.~[Z2:\,
p.~351]):

\proclaim {Proposition 2.3.8}
Let $D_{\ell,P}$ be the space of type $S-\QQ$ associated to $\widehat
G_{\ell,P}$
(introduced after {\rm (2.2.8)}), and $X_{\ell,P}$
the quotient of
$D_{\ell,P}$ by $\Gam(G_{\ell,P})$.  Let $\overline\pi^\r : \Xbarred\to X^*$
be the
natural mapping.  For $x\in M_P\subset X^*$,
$$
(\overline\pi^\r)^\-(x)\simeq (\overline {X_{\ell,P}})^{\text{red}},
$$
the  reductive Borel-Serre compactification of $X_{\ell,P}$.\sq
\endproclaim

In analogy with (2.3.6), we write
$$
R^{\text{exc}} = (R_{\ell,P})\times\GhP\times V_P,\tag 2.3.9
$$
with $V_P=W_P/U_P$ again.  It is easy to see that the following analogue
of Proposition 2.3.8 holds:
\proclaim {Proposition 2.3.10} With notation as in {\rm Proposition 2.3.8,}
let $\overline\pi^\e:
\Xbarexc\to X^*$ be the natural mapping.  For $x\in M_P\subset X^*$,
$(\overline\pi^\e)^\-(x)$
is a $\Gam(V_P)$-fibration over $\overline {X_{\ell,P}}$.\sq
\endproclaim

It is the central issue in this article that, except in easy cases (viz.,
when $G$ is of
$\QQ$-rank zero or $\RR$-rank one), the tower (2.3.5) is incompatible with
the tower
(2.2.18), in that there are no morphisms of compactifications of $X$, in
either
direction, between a Borel-Serre and a toroidal compactification of any sort
(plain,
excentric, or reductive).  It was shown by L.~Ji that $\text{GCQ}(\Xbarred,
X^\t) = X^*$ [J]. In [GT],
continuous
mappings $\XSig\to\Xbarred$ are constructed that are homotopy equivalent to
a morphism
of compactifications, and we will be elaborating on that in \S 3.

Suppose that a discrete group $\Theta$ acts continuously on a space $Y$.
Recall that Borel
defined a homotopy class of spaces $\langle\Theta,Y\rangle$ that is realized
by $\Theta
\backslash Z$ for any free $\Theta$-space $Z$ mapping $\Theta$-equivariantly
to $Y$.
This is known as the {\it Borel construction}.  It is at the heart of the
notion of
equivariant cohomology.  From (2.3.9) and Corollary 2.2.20, we deduce:
\proclaim{Proposition 2.3.11}
The excentric boundary strata $\overset\circ\to Z_R^\e$ and $e'(R)^\e$ are homotopy
equivalent,
both giving models for the Borel construction $\langle\Gam(G_{\ell,R}),
\Cal A_P\rangle$
for the action of $\Gam(G_{\ell,R})$ on $\Cal A_P$.\sq
\endproclaim
\demo{{\rm (2.3.12)} Remark} The analogous statement can be seen to hold for
the closed strata,
i.e., $Z_R^\e$, $\overline{e'(R)^\e}$, and $\langle\Gam(G_{\ell,R}),
\overline{\Cal A_P}^{\,\e}\rangle$.
Such facts are suggestive of (3.5.11) below.
\enddemo

We conclude with a curious assertion that holds in the non-Hermitian case as
well.  Any
{\it Satake
compactification} $X^{Sa}$ of $X$ (in the sense of [Z3]) is a quotient of
$\Xbarred$
[Z3:(3.11)].
An instance of this is $X^{Sa}=X^*$, as in (2.3.5).  We have morphisms of
compactifications of $X$:
$$
\Xbar\to\Xbarred\to X^{Sa},\tag 2.3.13
$$
which is the $\Gam$-quotient of partial compactifications of $D$:
$$
\Dbar\to\Dbarred\to D^{Sa}.\tag 2.3.14
$$
If we assume that $\Gam$ is neat, we have that $\Gam$ acts freely on
$\Dbar$. From
[Z3], one sees that the fibers of the mappings in (2.3.14) are contractible,
so these
mappings are homotopy equivalences (indeed, as $\Dbar$ is contractible, so
are the other
two spaces).  We need say no more (compare [HZ1:(3.9.1)]):
\proclaim{Proposition 2.3.15}  Let $D$ be a symmetric space of non-compact
type, and
$\Gam$ a neat arithmetically-defined group of isometries of $D$.  Then
$\Xbar$ is a model
for the Borel constructions $\langle\Gam,\Dbarred\rangle\approx\langle\Gam,
D^{Sa}\rangle$.\sq
\endproclaim
\medskip

{\bf (2.A)} {\it Appendix to {\rm \S 2}:} {\bf Hybrid
compactifications.} We wish to bring work of L.~Ji [J] into the
current context.  When $X$ is a locally symmetric variety, the
spaces $\XSig$ and $\Xbar$ are quite different in general, as we
mentioned at the end of (2.3).  Ji showed that their greatest
common quotient is a topological space he called the {\it
intermediate compactification}, which we will notate as $X^J$,
that differs from $X^*$ only in its highest-dimensional boundary
stratum.  An easier, but less detailed version of that is the one
intended in [HZ2:\,Conj.\,(1.5.8)]: the greatest common quotient
of $\XSig$ and $\Xbarred$ is $X^*$.

We first present an auxiliary construction, namely the determination of
certain
``group-theoretic" compactifications between $\Xbarexc$ and $\Xbarred$.  This
gets decided on
$\overline D^{\text{exc}}$, so we consider the quotient mappings $e(R)^{
\text{exc}}\to
e(R)^{\text{red}}$ for all $R$.  We seek to
determine whether one gets a Hausdorff space by contracting the fibers for
only certain $R$;
to give them a name,
we call such spaces {\it hybrid compactifications}.
Taking the question back to $\overline D$,
we will use the criterion implicit in [Z1]:
\proclaim {(2.A.1) Condition}
For each $R$, let $\Pi (R)=P$, and take $\widehat W_R$ to be one of $W_R$ or
$U_P$.  Assume that
this is done compatibly
with $G(\QQ)$-conjugacy: $\widehat W_{gRg^\-} = g(\widehat W_R)g^\-$.  Then
the following
condition must be satisfied: $\widehat W_R\subseteq\widehat W_{R'}$
whenever $R'\subset R$.
\endproclaim

\noindent {\it Remark.} The above is satisfied, of course, for the choice
$\widehat W_R
=W_R$ for all $R$, or when $\widehat W_R=U_P$ for all $R$.
\medskip

We make use of the following simple lemma, which we state without proof:
\proclaim {(2.A.2) Lemma}
Let $R'\subset R$ be parabolic, with $\Pi(R)=P$ and $\Pi(R')=P'$.  Then
\smallskip

i) $U_P\subseteq W_{R'}$ always,
\smallskip

ii) $W_R\subseteq U_{P'}$ only if both $R=P'=P$ and $W_P=U_P$.\sq
\endproclaim

From the preceding lemma, it follows that aside from $\Xbarred$ there are
only $r$
distinct hybrid
compactifications of the above sort: for each maximal parabolic type $Q$, put
$\widehat W_R=W_R$ if and only if
$\Pi(R)\prec Q$, and call the space $X^{\prec Q}$.
(Note that we get $\Xbarexc$ when $Q$ is smallest with respect to ``$\prec$".)
These fit
into a tower, with a mapping $X^{\prec Q}\to X^{\prec P}$ whenever $Q\prec P$.
We have been leading up to:
\proclaim {Proposition 2.A.3}
When $Q$ is the largest maximal parabolic type, there is a Cartesian square
$$\CD
X^{\prec Q} & @>>> & X^J\\
@VVV && @VVV\\
\Xbarred & @>>> & X^*
\endCD
$$\sq
\endproclaim
\tenpoint\baselineskip=16pt

\bigskip

\centerline{\bf 3. The least common modification of $\Xbarexc$ and $\XSigexc$}
\medskip

    In this Section, we will give an elaboration on a theorem of [GT],
which asserts that the natural mapping
$\L(\Xbarred,\XSig)\to\XSig$ has contractible fibers (so is a
homotopy equivalence). We make systematic use of the language of
\S 1.
\medskip

{\bf (3.1)} {\bf Main result and elements of the proof.} Our
present goal is to prove: \proclaim{Theorem 3.1.1} The canonical
morphism $\L(\Xbarexc,\XSigexc)\to \XSigexc$, having contractible
fibers, is a homotopy equivalence.
\endproclaim
\noindent We will then see (Corollary 3.5.13, (ii)) that the
corresponding assertion (Goresky-Tai) for $\L(\Xbarred,
\XSig)\to\XSig$ is deducible from Theorem 3.1.1 via
$\L$-basechange (1.1.8).

To begin, we have from \S2 that the spaces $\Xbarexc$ and $\XSigexc$ are
stratified, and have $X^*$
as a
(stratified) common quotient.  By Proposition 1.1.5, the assertion of
Theorem 3.1.1
can get decided locally over $X^*$.  In both cases, the fibers of the natural
mapping
to $X^*$ over the stratum $M_P$ involve $V_P$ and the homogeneous cone $C_P$
from (2.2),
through the boundaries of partial compactifications of $C'_P :=\Gam(\GlP)
\backslash C_P$ (see (3.5) below).

% In the case of $\Xbarexc$, it is the Borel-Serre compactification
% $\Xbar_{\ell,P}$ at infinity; in the case of $\XSigexc$, it is
% $\Gam(\GlP)\backslash |\widehat\Sigma_P^c|^\vee$ (a compact subset of
% $X_{\ell,P}$ itself); it is a deformation retract of the
% manifold-with-corners $\Xbar_{\ell,P}$.

Let $\check D(P)$ denote the open subset $U_P(\CC)\cdot D$ of the
compact dual of $D$ (cf.~(2.2.1)). In
[AMRT:\,III,\,pp.\,235--236,\,250] we find the following, which
allows us to examine the excentric Borel-Serre and toroidal cases
simultaneously: \proclaim{Proposition 3.1.2} With $P$ acting on
$C_P$ and $D_P$ through $P^\r$, the canonical $P$-equivariant
decomposition
$$
U_P(\RR)\backslash D\simeq C_P\times D_P\times V_P(\RR)\tag 3.1.2.1
$$
extends to a $P$-equivariant decomposition
$$
U_P(\RR)\backslash\check D(P)\simeq U_P(\RR)^\im\times [U_P(\CC)\backslash
\check D(P)]
\simeq U_P(\RR)^\im\times D_P\times V_P(\RR),\tag 3.1.2.2
$$
where ``$\im$" reminds us that this copy of $U_P(\RR)$ is
represented by the imaginary part of $U_P(\CC)$.\sq
\endproclaim
\noindent It follows that the excentric quotients on $\Xbar$ and $X^\t$ are
induced by
constructions on (3.1.2.2), after which one must undo the $U_P(\RR)$-quotients
in the interior.

Underlying Proposition 3.1.2 are the facts that $\GhP$ centralizes
$U_P$,
and that the
change of basepoint for the Cayley transform can be effected by translation
by an
element of $U_P(\RR)^\im$ (which is a consequence of the same for $G=SL_2$).
In effect,
when working over the Baily-Borel boundary stratum $M_P$, after acknowledging
the role of $V_P(\RR)$, we can ``replace"
$X$ by $C_P$
and determine the least common boundary modification (see (1.1.3)) of the
two partial
compactifications of the cone.

We will also make use of the fact mentioned in (2.3) that
$\Gam(\GlP)\backslash \widehat \Sigma_P$ is a Satake
compactification of $X_{\ell,P}$, so is in particular a quotient
of $\Xbar^\r_{\ell,P}$ [Z3:(3.11)].
\medskip

{\bf (3.2)} {\bf Real and complex toroidal embeddings; dual
compactifications.} A spanning simplicial cone $\sigma$, as in
(1.2), determines a smooth (complex) affine torus embedding:
$$
T\simeq (\CC^*)^d\subset \CC^d\simeq T_\sigma\tag 3.2.1
$$
in the usual manner (see [AMRT:\,I,\S1]).  As described in [HZ1,\S 2],
$\sigma$ also determines
a real torus
embedding $T_{\RR;\sigma}$.  The latter admits two descriptions: it is the
quotient
of $T_\sigma$ by the maximal compact torus of $T$ (viz., $T^c\simeq (S^1)^d$),
and it is
represented by the closure of $(\RR^+)^d$ in $T_\sigma$.  Then $\sigma\subset
T_{\RR;\sigma}
\simeq (\RR^{\ge 0})^d$ is given by $(0,1]^d$ (see (1.2)).

Let $\Sigma$ be a fan, as in (2.2).  This determines a complex
torus embedding $T_\Sigma$ by gluing the $T_\sigma$'s, and
likewise a real torus embedding $T_{\RR;\Sigma}$.  Give $|\Sigma|$
the largest topology for which the inclusion mappings of its cones
are continuous. We then define two partial compactifications,
$\overline \Sigma_1$ and $\overline\Sigma_2$, by attaching to
$\Sigma$ boundaries induced by the constructions in (1.2) on each
simplex; with $\Sig\simeq\widehat\Sig\times (0,\infty)$,
$$
\overline\Sigma_1 = \bigcup\,\{\overline\sigma_1:\sigma\in\Sigma\}\qquad
\text{and}\qquad
\overline\Sigma_2 = \bigcup\,\{\overline\sigma_2:\sigma\in\Sigma\},
\tag 3.2.2
$$
likewise with the corresponding topology for the inclusions $\overline\sigma_i
\hookrightarrow
\overline\Sigma_i$. Also, put
$$
\widehat\Sigma_1 = \d\overline\Sigma_1 \qquad\text{and}\qquad
\widehat\Sigma_2 = \d\overline\Sigma_2.\tag 3.2.3
$$
Then both $\widehat\Sigma_1$ (plainly) and $\widehat\Sigma_2$ (cf.~(1.2.4))
are
homeomorphic to $\widehat\Sigma$.

\proclaim{Lemma 3.2.4} For any $\widehat\tau\in\widehat\Sigma$, let $\tau^\vee(
\widehat
\sigma)$ be as in {\rm Proposition 1.2.1}, and put
$$
\tau^\vee=\bigcup\,\{\tau^\vee(\widehat\sigma)\,|\,\widehat\sigma \text{ is a
top-dimensional
simplex in
$\widehat\Sigma$ that contains $\widehat\tau$}\}.\tag 3.2.4.1
$$
Then $\tau^\vee$ is contractible.
\endproclaim
\demo{Proof} In the notation of (1.2), we have that for $\tau=\tau_J
\subseteq
\sigma$, $\tau$ is given by the equations $t_j=1$ for all $j\in J$, and
$\tau^\vee(
\widehat\sigma)$ is given by the equations $t_j=0$ for all $j\notin J$.  This
contracts to
the point that has further equations $t_j=1$ for all $j\in J$, viz., $\widehat
\tau\cap
\tau^\vee(\widehat\sigma)$. This point is independent of the
top-dimensional
cone $\sigma$ (cf.~[HZ1:\,2.3]), so we can do the contractions simultaneously
for all $\sigma$.\sq
\enddemo
\noindent (3.2.5) {\it NB---}The union in (3.2.4.1) can be
infinite if $\widehat\tau$ is in the boundary of $\widehat\Sigma$.
%
% \noindent A partial picture of (3.2.4.1) when $d=2$ is:
% \smallskip
%
% \noindent (3.2.5) ***FIGURE 2

\smallskip

Recall the notion of $\Lb$ that was defined in (1.1.3).  We have
the following extension of Proposition 1.2.1, which we write for
$\Sig=\Sig_P$:
\proclaim{Proposition 3.2.6 {\rm (Canonical
duality)}}
$\Lb(\widehat\Sigma_1,\widehat\Sigma_2)=\bigcup\,\{\widehat\tau\times\tau^\vee:\,
\widehat\tau\in\widehat\Sigma_P\}$.
\endproclaim
\demo{{\rm (3.2.7)} Remark} The formula in Proposition 3.2.6 can
be rewritten as
$$
\Lb(\widehat\Sigma_1,\widehat\Sigma_2)=\bigsqcup\,\{\widehat\tau^\circ\times
\tau^\vee: \widehat\tau\in\widehat\Sigma_P\},\tag 3.2.7.1
$$
for $\tau^\vee\subset\omega^\vee$ when $\widehat\omega$ is a face of $\widehat
\tau$.
Also, we note that the interior of $\widehat\Sigma_1$ can be written as
$$
\widehat C_P =
\bigsqcup\{\widehat\tau^\circ:\widehat\tau\in\widehat\Sigma_P^
\circ\}.\tag 3.2.7.2
$$
In the above, $\widehat\tau^\circ$ denotes the interior of the simplex
$\widehat\tau$.
We can switch the roles in (3.2.7.1) and write:
$$
\Lb(\widehat\Sigma_1,\widehat\Sigma_2)=\bigcup\,\{\widehat\tau\times
(\tau^\vee)^\circ: \widehat\tau\in\widehat\Sigma_P\}.\tag 3.2.7.3
$$

\demo{Proof of Proposition {\rm 3.2.6}}
Let $\{x_j\}$ be a sequence in $\Sig$ that converges in both
$\overline\Sigma_1$
and $\overline\Sigma_2$.  Suppose that the limit in $\widehat\Sigma_1$ lies
in the
interior $\widehat\tau^\circ$ of the simplex $\widehat\tau$. We claim that
there is a
top-dimensional cone $\sigma$ containing $\tau$ (as a face) such that $\sigma$
contains a
subsequence of $\{x_j\}$ (because there can be infinitely many such $\sigma$,
there is
something to check).  If not, each such $\sigma$
contains
only finitely many $x_j$'s.  Thus, there is a neighborhood $N_\tau(\sigma)$
of $\tau$ in
$\sigma$ that contains no $x_j$'s. By the definition of the topology, $\bigcup_
\sigma N_
\tau(\sigma)$ is open in $\Sigma$, and it contains no $x_j$'s. This contradicts
the convergence.
It follows that we may assume that the sequence is in a single $\sigma$, and
we are
reduced to Proposition 1.2.1: the second limit lies in $\tau^\vee(\widehat
\sigma)$.
Since this happens for every $\sigma\supset\tau$, we use
(3.2.4.1) and
our formula for $\Lb(\widehat\Sigma_1,\widehat\Sigma_2)$ follows.\sq
\enddemo

Suppose next that a discrete group $\Gam_\ell$ acts linearly on a
cone complex $\Sigma$ in a separated and discontinuous manner,
such that there are only finitely many $\Gam_\ell$\,-equivalence
classes of cones.  We use ``prime" to indicate a
$\Gam_\ell$-quotient, as in $\Sigma'= \Gam_\ell \backslash\Sigma$,
and do likewise in (3.2.2) to obtain
$$
\overline\Sigma'_1=\Gam_\ell\backslash\overline\Sigma_1\qquad
\text{and} \qquad\overline\Sigma'_2 =
\Gam_\ell\backslash\overline\Sigma_2.\tag 3.2.8
$$
The spaces in (3.2.8) are viewed as compactifications of $\Sig'$.
Let $\widehat\Sigma$ be the corresponding simplicial complex, as
in (2.2).  Then we have $\widehat\Sigma'=\Gam_\ell \backslash
\widehat\Sigma$, a compact space.  From (3.2.3) and (3.2.6) we get
a pair of polyhedral models of $\widehat\Sigma'$ having finitely
many faces:
$$
\widehat\Sigma'_1=\Gam_\ell\backslash\widehat\Sigma_1\qquad\text{and}\qquad
\widehat \Sigma'_2=\Gam_\ell\backslash\widehat\Sigma_2,\tag 3.2.9
$$
with the latter given as in (1.2.4).

Since $\Gam_\ell$ acts on $\overline\Sigma_1$ and
$\overline\Sigma_2$, it acts of their $\L$, i.e. the $\Lb$ in
Proposition 3.2.6.  The expected conclusion holds to a large
extent:

\proclaim{Proposition 3.2.10} If $\Gam_\ell$ is neat, the actions
of $\Gam_\ell$ on $\overline\Sigma_1$ and $\overline\Sigma_2$ are
diagonal, and
$$
\Lb(\widehat\Sigma'_1,\widehat\Sigma'_2)=\Gam_ \ell\backslash
\bigcup\,\{\widehat\tau\times\tau^\vee:\widehat\tau\in\widehat\Sigma\}.
\tag 3.2.10.1
$$
\endproclaim
\demo{Proof} By hypothesis, there is a subset $S$ of $\Sig$, the
union of finitely many closed cones, such that the mapping
$S\to\Sig\to\Sig'$ is surjective. This gives the same for the
respective closures of $S$ in $\overline\Sigma_1$ and
$\overline\Sigma_2$:
$$
\overline S_1\to\overline\Sigma_1\qquad\text{and}\qquad\overline
S_2\to \overline\Sigma_2.
$$
Moreover, $\overline S_1$ and $\overline S_2$ are compact.  Thus,
we are in the setting of Proposition 1.3.8, provided we verify
that the actions of $\Gam_\ell$ on $\overline\Sigma_1$ and
$\overline\Sigma_2$ are diagonal.

Let $\widehat\tau\in\widehat\Sigma_1=\d\overline\Sigma_1$, $e_1\in
\widehat\tau$ and $e_2\in(\tau^\vee)^\circ$. As $\widehat\tau$
varies, this gives all points $(e_1,e_2)$ of $\Lb(\widehat
\Sigma_1, \widehat\Sigma_2)$, by (3.2.7.3).  Suppose that
$(\gamma\cdot e_1, e_2)\in\Lb(\widehat\Sigma_1, \widehat\Sigma_2)$
for some $1\ne\gamma\in\Gam_\ell$.  A priori, there are two cases.
If $\gamma\cdot \widehat\tau= \widehat\tau$, then $\gamma$ acts as
a permutation of the vertices of $\widehat\tau$, so it fixes the
barycenter of $\widehat\tau$. Since $\Gam_\ell$ us assumed neat,
$\gamma$ must fix the whole boundary component containing
$\tau^\circ$, hence also its closure. In particular, $\gamma\cdot
e_1 = e_1$. On the other hand, suppose that
$\gamma\cdot\widehat\tau\ne\widehat\tau$. The only way that both
$(e_1,e_2)$ and $(\gamma\cdot e_1, e_2)$ could be in the $\Lb$ is
if $e_2\in\tau^\vee \cap (\gamma\cdot\tau)^\vee$. However, this
set would lie in the {\it boundary} of $\tau^\vee$, a
contradiction.  The conclusion is, then, that $\gamma\cdot e_1 =
e_1$, so the actions are diagonal.

% If there are such simplices in $\widehat\Sigma_1$, there will not be
% any in $\widehat\Sigma_1^{(1)}$, the barycentric subdivision of
% $\widehat\Sigma_1$ unless $\gamma=1$.  Said another way, if
% $\widehat\Sigma_1$ is the barycentric subdivision of another
% $\Gam_\ell$-invariant simplicial complex,---take this as the
% meaning of ``sufficiently fine" in the statement of this
% proposition---we must have $\gamma\cdot e_1=e_1$, so strong
% diagonality holds.

This gives that the canonical mapping
$$
\Gam_\ell\backslash\Lb(\widehat\Sigma_1,\widehat\Sigma_2)\to\Lb(\widehat
\Sigma'_1, \widehat\Sigma'_2)
$$
is an isomorphism. We use Proposition 3.2.6 to obtain formula
(3.2.10.1).\sq
\enddemo
\demo {{\rm (3.2.10.2)} Remarks} i) One might try to prove
Proposition 3.2.10 by using (3.2.7.1) instead of (3.2.7.3).  The
reader is invited to investigate why we took the other route.
\smallskip

ii) If $\widehat\tau$ above is just the vertex $\{e_1\}$, the
situation is very much like that in (1.3.6), but here we keep
$\Gam_\ell$ fixed.
\enddemo

\medskip

{\bf (3.3)} {\bf Adjustments to duality.} For a pair of reasons,
$\L(\Xbarexc, \XSigexc)$ cannot be determined from Proposition
3.2.10 alone. First, the projection (from Proposition 3.1.2) onto
$\widehat C'_P$ of the $P$-stratum of $\XSigexc$ is not dense in
$\widehat\Sigma'_2$; it is only the (compact) $\Gam(
\GlP)$-quotient of the subcomplex
$$
(\widehat\Sigma_P^c)^\vee = \bigcup\,\{\tau^\vee:\widehat\tau\in\widehat
\Sigma_P^c\},\tag 3.3.1
$$
which is contained in the interior of $\widehat\Sigma_2$. Then
$(\widehat\Sigma_P^c)^\vee$ can be identified, via (1.2.4), with
$\widehat\Sigma_P^{(1)}$.
Here, the
barycentric subdivision $\Sigma_P^{(1)}$ of $\Sigma_P$ also gives a compatible
collection of fans as $P$ varies. In Corollary 3.3.5 below, we
show that
$$
(\widehat\Sigma_P^c)^\vee\isoarrow\,(\widehat\Sigma_P^\circ)^\vee,\tag 3.3.3
$$
where $\widehat\Sigma_P^\circ$ denotes the $\Gam(\GlP)$-invariant subset of
$\widehat
\Sigma_P$ given by (2.2.9).

%%% = \Gam(\GlP)\backslash(\widehat\Sigma_P^{(1)})

Looking from the inside, we have that going to infinity in $\Xbar$
occurs not just by letting the dilation variable of the cone,
relative to a cross-section of $\Sigma_P$, go to $\infty$. Rather,
it must go to infinity at least as fast as a so-called {\it core}
for the cone,---see [AMRT:\,II,\,\S 5] for the notion of a
core---which gives a cross-section of $C_P\subset\Sigma_P$ that
blows up at the boundary cones of $\Sigma_P$ (the precise meaning of a core is
irrelevant for this article). One gets something
essentially the same by using the ``$\GlP$ side" of the
description in [Z2:\,3.18].

We take a brief excursion into the calculus of joins. Let
$\widehat\omega$ be the face of a simplex $\widehat\sigma$ that is
opposite the face $\widehat\tau$ of $\widehat\sigma$, i.e., $\widehat\omega$
is spanned by the vertices of $\widehat\sigma$ that are {\it not} in
$\widehat\tau$. Then $\widehat \sigma$ is the {\it join}
$\widehat\omega*\widehat\tau$, namely the quotient of
$\widehat\omega \times [0,1]\times\widehat\tau$ obtained by
separately collapsing $\widehat \omega\times\{1\}$ to a point and
$\{0\}\times\widehat\tau$ to a point.  The face $\widehat\omega$
of $\widehat\sigma$ lies at $s=0$ and $\widehat\tau$ is at $s=1$,
where $s$ denotes the variable of [0,1]. One can collapse
$\widehat\sigma-\widehat\omega$, and also its subset $\tau^
\vee(\sigma)$, onto $\widehat\tau$ by the following deformation
retraction: $\{h_s\}$ is given by mappings on the product
$\widehat\omega\times(0,1]\times\widehat\tau$ that depend only on
the interval variable: $f_s(t) = \min\{t+s,1\}$;
$h_s(w,t,x)=(w,f_s(t),x)$. Note that $\tau^ \vee(\sigma)$ is
closed under the flow.

The setting above is for disjoint faces of one simplex in a complex. However,
the notion can easily
be extended to intersecting faces.  Let $\widehat\tau$ and $\widehat\omega$
be simplices
with intersection $\widehat\alpha$.  Write $\widehat\tau = \widehat\alpha *
\widehat\tau'$ and
$\widehat\omega = \widehat\alpha * \widehat\omega'$. Then $\widehat\tau *
\widehat\omega =
\widehat\tau * \widehat\omega' =\widehat\tau' *\widehat\omega=\widehat\alpha *
\widehat\tau'
*\widehat\omega'$, the smallest
simplex containing both $\widehat\tau$ and $\widehat\omega$. We are interested
in the dual, in
the sense of (3.2.4.1), of a join.
% We will treat only the join of disjoint faces, for the general case
% follows easily:
\proclaim {Proposition 3.3.4} i) Let $\widehat\tau$ and $\widehat\omega$ be
faces of
a top-dimensional simplex $\widehat\sigma$ in a simplicial complex $\widehat
\Sigma$. Then
$$
(\widehat\tau * \widehat\omega)^\vee(\widehat\sigma)=\tau^\vee(\widehat\sigma)
\cap\omega^\vee
(\widehat\sigma).\tag 3.3.4.1
$$

ii) Conversely, if $\widehat\tau$ and $\widehat\omega$ are simplices in
$\widehat\Sigma$
such that $\tau^\vee\cap\omega^\vee\ne\emptyset$, then $\widehat\tau * \widehat
\omega$ is
defined and
$$
(\widehat\tau * \widehat\omega)^\vee=\tau^\vee\cap\omega^\vee.\tag 3.3.4.2
$$
\endproclaim
\demo{Proof} We use the notation from the proof of Proposition 1.2.1. If we
write
$\widehat\tau = \widehat\tau_J$ and $\widehat\omega = \widehat\tau_K$, then
$\widehat\tau *
\widehat\omega =\widehat\tau_{(J\cap K)}$. Similarly, we have $\tau^\vee=\tau^
\vee_J$, with
defining equations parametrized by the {\it complement} of $J$, and likewise
for $K$, so
$\tau^\vee_J\cap\tau^\vee_K = \tau^\vee_{(J\cap K)}$. This gives i).

As for ii), because more than one top-dimensional simplex is involved in the
dual, we
must be a little careful.  From the definition, we have
$$
\tau^\vee\cap\omega^\vee = \bigcup\,\{\tau^\vee(\widehat\sigma)\cap\omega^\vee(
\widehat\sigma'):
\widehat\sigma,\,\widehat\sigma'\text{ top-dimensional simplices in $\widehat
\Sig$}
\}.\tag 3.3.4.3
$$
There is a natural notion of $\tau^\vee(\widehat\upsilon)$ for any simplex
$\widehat\upsilon
\supset\widehat\tau$: for any top-dimensional $\widehat\sigma\supset\widehat
\upsilon$,
take $\tau^\vee(\widehat\sigma)\cap\widehat\upsilon$, and check that this is
independent
of the choice of $\widehat\sigma$. Then $\tau^\vee(\widehat\sigma)\cap\omega^
\vee(\widehat
\sigma')=(\tau^\vee\cap\omega^\vee)(\widehat\sigma\cap\widehat\sigma')$.
We see
that we can then rewrite (3.3.4.3) as
$$
\tau^\vee\cap\omega^\vee =
\bigcup\,\{\tau^\vee(\widehat\sigma)\cap\omega^\vee(
\widehat\sigma):\widehat\sigma\text{ is a top-dimensional simplex
in $\widehat \Sig$}\}.\tag 3.3.4.4
$$
If the term in the right-hand side of (3.3.4.4) coming from $\widehat\sigma$
is non-empty,
then both $\widehat\tau$ and $\widehat\omega$ are faces of $\widehat\sigma$.
% (Of course,
% there can be more than one, even infinitely many, $\widehat\sigma$ with this
% property.)
We can now proceed as in i), term by term.\sq

We apply Proposition 3.3.4 to $\widehat\Sigma_P$:
\proclaim{Corollary 3.3.5}
Let $\widehat\tau\in\widehat\Sigma_P^\circ$.  Then $\tau^\vee\in(\widehat
\Sigma_P^c)^\vee$.
\endproclaim
\demo{Proof} Since $\widehat\tau$ has a vertex in
$\widehat\Sigma_P^c$, we can write it in the form
$\widehat\alpha*\widehat\omega$, with
$\widehat\alpha\in\widehat\Sigma_P^c$. It follows from
(3.3.4,\,ii) that $\tau^\vee\subseteq\alpha^\vee$.\sq
\enddemo

\noindent We remember this fact as the equivalent (3.3.3).

With that said, the calculations in (3.2) get used as follows. We
now regard $\overline\Sigma_1$ as an element of $\PCp(C_P)$. Let
$(\Sigma_P)^\te$ be the partial compactification of $C_P$
contained in $\overline\Sigma_2$ with $\d (\Sigma_P)^\te
=(\widehat \Sigma_P^c)^\vee \subset \widehat\Sigma_2$. We consider
$$
\Lb(\widehat\Sigma_1,(\widehat\Sigma_P^c)^\vee)=\d\L(\overline\Sigma_1,(\Sigma_
P)^\te)\subset \Lb(\widehat\Sigma_1,\widehat\Sigma_2),\tag 3.3.5.1
$$
in which the rightmost space is given by (3.2.7.3).  If
$\widehat\tau\in \widehat \Sigma_P^\circ$, then $\tau^\vee$ is
interior to $\d\overline\Sig_2$ by Corollary 3.3.5, and thus the
contributions of $\widehat\tau$ to
$\Lb(\widehat\Sigma_1,(\widehat\Sigma_P^c)^\vee)$ and
$\Lb(\widehat\Sigma_1,\widehat\Sigma_2)$ are the same. At issue,
then, are the simplices $\widehat\tau$ contained in the boundary
of $\widehat\Sig_P$.

We see from (3.2.7.1) and (3.2.7.2) that the subset of
$\Lb(\widehat\Sigma_1,(\widehat\Sigma_P^c)^\vee)$ that maps to
$\widehat C_P$ under the projection onto $\widehat\Sigma_1$ is
given by
$$ \Cal T :=
\bigsqcup\,\{\widehat\tau^\circ\times\tau^\vee:\widehat\tau\in
\widehat\Sigma_P^\circ\}.\tag 3.3.6
$$
$\Cal T$ visibly maps {\it onto} $(\widehat\Sigma_P^c)^\vee$ under
the projection
$\Lb(\widehat\Sigma_1,(\widehat\Sigma_P^c)^\vee)\to(\widehat\Sigma_P^c)^\vee$.
Therefore, $\Lb(\widehat\Sigma_1,(\widehat\Sigma_P^c)^\vee)$
necessarily contains the closure $\widetilde{\Cal T}$ of $\Cal T$
in $\widehat\Sigma_1\times (\widehat\Sigma_P^c)^\vee$:
$$
\Lb(\widehat\Sigma_1,(\widehat\Sigma_P^c)^\vee)\supseteq\widetilde{
\Cal T} =
\bigcup\,\{\widehat\tau\times\tau^\vee:\widehat\tau\in\widehat\Sigma_P^\circ\}.\tag
3.3.7
$$
The following determination is fundamental: \proclaim {Proposition
3.3.8} $\Lb(\widehat\Sigma_1,(\widehat\Sigma_P^c)^\vee)$ is equal
to $\widetilde{\Cal T}$.
\endproclaim
\demo{Proof} One must keep in mind that $\widehat\Sigma_P^\circ$
is not a complex.  Let $\widehat \omega$ be a simplex in
$\widehat\Sigma_P$. (We are thinking that $\widehat \omega$ is in
the boundary of $\widehat C_P$, but it does not matter.) In
canonical duality this contributes, in (3.2.7.1),
$\widehat\omega^\circ\times\omega^\vee$ to $\Lb( \widehat\Sig_1,
\widehat\Sig_2)$; in $\Lb(\widehat\Sig_1,(\Sig_P^c)^\vee)$ the
contribution is
$$
\widehat\omega^\circ\times(\omega^\vee\cap(\Sig_P^c)^\vee).\tag
3.3.8.1
$$
Suppose that $\widehat\tau\in\widehat\Sig_P^c$ and
$\omega^\vee\cap\tau^\vee\ne \emptyset$. By Proposition
3.3.4,\,ii), $\widehat\omega*\widehat\tau$ is defined and its dual
is $(\widehat\omega*\widehat\tau)^\vee =
\tau^\vee\cap\omega^\vee$. Since $\widehat
\omega*\widehat\tau\in\widehat\Sig_P^\circ$, we get that
$$
\widehat\omega^\circ\times(\omega^\vee\cap\tau^\vee)\subset(\widehat\omega
*\widehat\tau)\times(\widehat\omega*\widehat\tau)^\vee.\tag
3.3.8.2
$$
As $(\omega^\vee\cap\tau^\vee)\in (\widehat\Sig_P^c)^\vee$ by
Corollary 3.3.5, we are done.\sq
\enddemo

\noindent Thus, the first adjustment to canonical duality is to
replace $\widehat\Sig_P$ by $\widehat\Sig_P^\circ$ in the formula
in Proposition 3.2.6 (see (3.3.11.2) below.)

The second one can now be handled rather quickly. Given the cone $C=C_P$, let
$s:\widehat
C\to C$ be any cross-section of the cone.  Here we do not assume that $s$
extends to
a cross-section of $\Sig_P$.  One uses $s$ to write
$$
C\simeq\widehat{C}\times (0,\infty).\tag 3.3.9
$$
This again determines variables $(r,\widehat x)$ on $C$, such that
$x=r\cdot s(\widehat x)$ for all $x\in C$.  We can define
$\overline C(s)$ in terms of (3.3.9) as $\widehat{C}\times
(0,\infty]$ (compare the definition of $\overline\sigma_1$ in
(1.2)).  Note that if $s_2
> s_1$, then
convergence to infinity in $\overline C(s_2)$ implies the same in $\overline
C(s_1)$,
and the limits coincide in $\widehat C\times\{\infty\}$.

Take $s$ to be a section determined by a $\Gam(\GlP)$-invariant core of $C$,
which always
exists [AMRT:\,p.123].  The stronger notion of convergence in $\overline C(s)$
than in
$\overline\Sigma_1$ implies that
$$
\L(\overline
C(s),(\Sigma_P)^\te)\subseteq\L(\overline\Sigma_1,(\Sigma_P)^\te)\tag 3.3.10
$$
But the set of pairs of limits we get, using $\widehat\Sig_1$, is
proved in Proposition 3.3.8 to be the least possible, so it must
be the same set here. It follows that we have equality in
(3.3.10); the fact that $s$ blows up at the boundary of $C$ is
irrelevant for the determination of the $\L$. In other words, the
second adjustment to canonical duality is vacuous.  In sum,
\proclaim{Proposition 3.3.11} Let $s$ be a $\Gam(\GlP)$-invariant
cross-section of $C$ that blows up at the boundary of $\Sig_P$.
Then
$$
\Lb(\widehat C(s),(\Sig_P^c)^\vee) =
\Lb(\widehat\Sigma_1,(\Sig_P^c)^\vee),\tag 3.3.11.1
$$
with the right-hand side given by {\rm Proposition 3.3.8} and {\rm
(3.3.7)}, i.e.:
$$
\Lb(\widehat\Sigma_1,(\Sig_P^c)^\vee)=\bigcup\,\{\widehat\tau\times
\tau^\vee:\widehat\tau\in\widehat\Sigma_P^\circ\}.\tag 3.3.11.2
$$
In particular, this hold when $s$ defines a core of $C_P$.\sq
\endproclaim
\medskip

{\bf (3.4)} {\bf Fibers of mappings of $\Gam(\GlP)$-quotients.} We
next treat the role of taking $\Gam(\GlP)$-quotients in
determining the fiber in mappings related to Theorem 3.1.1. {\it
We assume throughout this section only that $\Gam(\GlP)$ acts
freely on} $\widehat C_P$.

First, we make a general and elementary observation:
\proclaim
{Proposition 3.4.1} Suppose that $\Gam$ acts separated and
discontinuously on spaces $Z_1$ and $Z_2$.  Let $\psi:Z_1\to Z_2$
be a $\Gam$-equivariant mapping. Let $\widetilde F$ be the fiber
of $\psi$ at $z\in Z_2$, and $F$ the fiber of the induced mapping
$\Gam\backslash Z_1\to \Gam\backslash Z_2$ at the corresponding
point $[z]\in\Gam\backslash Z_2$. Then $F\simeq
\Gam_z\backslash\widetilde F$, where $\Gam_z$ is the isotropy
subgroup of $\Gam$ at $z$.
\endproclaim
\demo{Proof} It is immediate that $F\simeq\Gam\backslash(\Gam\!\cdot\!
\widetilde F)$.
The inclusion of $\widetilde F$ in $\Gam\!\cdot\!\widetilde F$ induces
$\Gam_z
\backslash\widetilde F\simeq F$.\sq
\enddemo
\proclaim{Corollary 3.4.2} Under the conditions of {\rm Proposition 3.4.1,}
suppose
that $\Gam$ acts freely on $Z_2$.  Then $F\simeq\widetilde F$.\sq
\endproclaim

% \demo{{\rm (3.2.11)} Remark} ***The above argument for diagonality
% applies to the $\L$ when  $\overline\Sig_1$ is replaced by a
% partial compactification of $C_P$ that maps onto it.
%\enddemo

% With $\d\L(\overline\sigma_1,\overline\sigma_2)$ given by (1.2.1.1), and the
% interior of
% $\overline\Sigma'_1$ and $\overline\Sigma'_2$ understood to be $\Sigma'$, we
% use (3.2.4.1) to obtain:? ***
% \demo{{\rm (3.2.8)} Remark}  We will use this for $\widehat\Sigma=\widehat
% \Sigma_P$ and $\Gam_\ell=\Gam(\GlP)$.
% Although it was not essential for (3.2.7), we would like to have
% $$
% \Lb(\widehat\Sigma_1,\widehat\Sigma_2)=\bigcup\,\{\tau\times\tau^\vee|
% \tau\in\widehat\Sigma\}
% $$
% ($\Gam$-quotients are not taken here).  The difficulty seems to be pushing a
% biconvergent
% sequence into a single cone.  The problem lies with boundary simplices with
% infinitely
% many elements, as in (3.2.4.1).

% {\it provided} we change the topology of $\widehat\Sigma_1$ to the {\it
% Satake topology} [S].  This topology is coarser than the given topology.
% That adds the condition: for $y\in\widehat\Sigma_1$ with stabilizer
% $\Gam_y$ in $\Gam_\ell$, there is a neighborhood basis consisting of
% $\Gam_y$-invariant open sets.

We use a ``prime" as a standard way to denote a
$\Gam(\GlP)$-quotient, with the one exception $\widehat\Sig^\prp_2
= \Gam(\GlP)\back\d(\Sig_P)^\te$. Because $\Gam(\GlP)$ acts freely
on $\widehat C_P$, the action on $(\widehat \Sig_P^c)^ \vee$,
which is a subset of $\widehat C_P$, is also free.  The next fact
follows immediately from Proposition 3.2.10 when $\Gam(\GlP)$ is
neat.  We provide a small modification of the argument to obtain
some improvement of the result:

\proclaim{Proposition 3.4.3} i) For
$\widehat\tau\in\widehat\Sig_P^\circ$, the fiber of
$\Lb(\widehat\Sigma_1,(\widehat \Sig_P^c)^\vee)\to(\widehat
\Sig_P^c)^\vee$, at any point in the interior $(\tau^\vee)^ \circ$
of $\tau^\vee$, is $\widehat\tau$. In particular, the fiber is
contractible.  Thus
$$
\Lb(\widehat\Sigma_1,(\widehat\Sig_P^c)^\vee)=
\bigcup\,\{\widehat\tau\times (\tau^\vee)^\circ:
\widehat\tau\in\widehat\Sigma_P^\circ\}=\bigcup\,
\{\widehat\tau\times\tau^\vee:
\widehat\tau\in\widehat\Sigma_P^\circ\}.\tag 3.4.3.1
$$

ii)  Under the mild assumption that $\Gam(\GlP)$ acts freely on
$\widehat C_P$, the actions of $\Gam(\GlP)$ on $\overline\Sigma_1$
and $(\Sigma_P)^\te$ are diagonal, so
$$
\Lb(\widehat\Sigma'_1,\widehat\Sigma^\prp_2)=\Gam(\GlP)\backslash
\bigcup\,\{\widehat\tau\times(\tau^\vee)^\circ:\widehat\tau\in
\widehat\Sigma_P^\circ\}.\tag 3.4.3.2
$$
The fiber of
$\Lb(\widehat\Sigma'_1,\widehat\Sigma^{\prime\prime}_2)\to
\widehat \Sig_2^{\prime\prime}$, at any point represented by a
point of $(\tau^\vee)^\circ$, is $\widehat\tau$, hence
contractible.
\endproclaim
\demo{Proof}
%  ***We use Corollary 3.4.2 together with (3.3.8)***
i) We note that $(\widehat\Sig_P^c)^\vee$ is a polyhedral complex,
whose strata are the sets of the form $(\tau^\vee)^\circ$, with
$\widehat\tau\in\widehat\Sig_P^\circ$. The analogue of (3.2.7.2),
namely
$$
(\widehat\Sig_P^c)^\vee =
\bigcup\,\{(\tau^\vee)^\circ:\widehat\tau\in\widehat\Sig_
P^\circ\},\tag 3.4.3.3
$$
holds.  The fiber of the projection
$\Lb(\widehat\Sigma_1,(\widehat\Sig_P^c)^ \vee)\to
(\widehat\Sig_P^c)^\vee$ over $(\tau^\vee)^\circ$ is
$\widehat\tau$ by canonical duality.  (One cannot switch the roles
of $\widehat\tau$ and $\tau^\vee$ here.) Thus:
$$
\Lb(\widehat\Sigma_1,(\widehat\Sig_P^c)^\vee)=
\bigcup\,\{\widehat\tau\times (\tau^\vee)^\circ:
\widehat\tau\in\widehat\Sigma_P^\circ\};\tag 3.4.3.4
$$
taking the closure of $(\tau^\vee)^\circ$ in (3.4.3.4) effects no
change.

ii) Again, there is a subset $S$ of $\Sig_P$, the union of
finitely many closed cones, such that the mapping
$S\to\Sig_P\to\Sig'_P$ is surjective. This gives the same for the
respective closures of $S$ in $\overline\Sigma_1$ and
$(\Sigma_P)^\te$:
$$
\overline S_1\to\overline\Sigma_1\qquad\text{and}\qquad S^\te\to
(\Sigma_P)^\te.
$$
As $\overline S_1$ and $S^\te$ are compact, we can then apply
Proposition 1.3.8, provided we verify that diagonality holds.

Let $\widehat\tau\in\widehat\Sigma_1^\circ$, $e_1\in\widehat\tau$,
and $e_2\in(\tau^\vee)^\circ$. As $\widehat\tau$ varies, this
gives all points $(e_1,e_2)$ of $\Lb(\widehat \Sigma_1,
(\widehat\Sigma_P^c)^\vee)$, by (3.4.3.3).  Suppose that
$(\gamma\cdot e_1, e_2)\in\Lb(\widehat\Sigma_1,
(\widehat\Sigma_P^c)^\vee)$ for some $1\ne\gamma\in\Gam(\GlP)$.
The discussion continues almost exactly as in the proof of
Proposition 3.2.10: if $\gamma\cdot \widehat\tau= \widehat\tau$,
then $\gamma$ fixes the barycenter of $\widehat\tau$. Here,
though, the barycenter of $\widehat\tau$ is a point of $\widehat
C_P$ because $\widehat\tau\in \widehat\Sigma_1^\circ$, in
contradiction to the freeness hypothesis. Thus,
$\gamma\cdot\widehat\tau\ne\widehat\tau$. This is not possible
either, as in the aforementioned proof.  The conclusion is that
the situation described in the beginning of this paragraph does
not occur, so the actions are diagonal.

This gives that the canonical mapping
$$
\Gam(\GlP)\backslash\Lb(\widehat\Sigma_1,(\widehat\Sigma_P^c)^\vee)\to\Lb(
\widehat \Sigma'_1, \widehat\Sigma^\prp_2)
$$
is an isomorphism. We use (3.3.11.2) to obtain formula (3.4.3.2).
By Corollary 3.4.2, the taking of $\Gam(\GlP)$-quotients does not
change the fiber, as $\Gam(\GlP)$ acts freely on
$(\widehat\Sig_P^c)^\vee$.  Our assertion follows.\sq
\enddemo

\demo{{\rm (3.4.3.5)} Remarks} i) Given (3.4.3.3), Proposition
3.4.3 asserts that {\it all} fibers of that mapping are
contractible.
\smallskip

ii) Item i) of the proposition is to be contrasted with
Proposition 3.4.12 below.
\smallskip

iii) The outcome of the check for diagonality in the above proof
is important enough that we restate it for emphasis: {\it if both
$(e_1,e_2)$ and $(\gamma\cdot e_1, e_2)$ ($\gamma\in\Gam(\GlP)$)
are in $\Lb(\widehat\Sigma_1, (\widehat\Sigma_P^c)^\vee)$, then
$\gamma = 1$.}
\smallskip

iv) An assertion from the proof of Proposition 3.2.10, tacitly
used in the above, is that $\tau^\vee
\cap(\gamma\cdot\tau)^\vee\ne\emptyset$. By (3.3.4.2), there must
be a simplex in $\widehat\Sigma_1$ (not necessarily unique)
containing both $\widehat \tau$ and $\gamma\cdot\widehat\tau$ as
faces, with $\tau^\vee \cap (\gamma\cdot\widehat\tau)^\vee =
(\widehat\tau*(\gamma\cdot\widehat\tau))^\vee.$ However, this
would not improve the argument.
\enddemo

We use Proposition 3.4.3 in conjunction with Proposition 3.3.11 to
yield:

\proclaim{Corollary 3.4.4} $\Lb(\widehat C'(s),\Sig^\prp_2)=
\Lb(\widehat\Sigma'_1,\Sig^\prp_2)$.\sq
\endproclaim

We did not really use much about $\widehat\Sigma_1$ in making the
argument for the key Proposition 3.3.8, which underlies our
calculations.  The main point was that $\widetilde{\Cal T}$ gave
the part of the boundary contained in $\widehat C_P\times
(\widehat \Sig_P^c)^\vee$.  As such, we can try to apply it
likewise to other $\Gam(\GlP)$-equivariant partial
compactifications of $C_P$. Most interesting are ones that have
familiar elements $Y\in\PCp(\widehat C_P)$ at the boundary;
writing $\widehat C_P$ as $D_{\ell,P}$, we have in mind
$Y=\Dbar_{\ell,P}$\,, $Y=\overline D^\r_{\ell,P}$\,, or $Y$ any
Satake partial compactification $D^{Sa}_{\ell,P}$ of $D_{\ell,P}$.

We need to be specific about the corresponding elements of
$\PCp(C_P)$.  With $\Sig$ written as $\widehat\Sig\times
(0,\infty)$, inducing $C_P\simeq\widehat C_P\times(0,\infty)$, we
take $Y\times (0,\infty]$ with $\d Y\times (0,\infty)$ removed; it is 
stated this way because a cross-section for one $Y$ need not
work for another (compare the end of (3.3)).
There is a ``version" of $\widetilde{\Cal T}$, i.e., the closure
of (3.3.6), for each case. These we denote respectively
$\widetilde{\Cal T}^{BS}$, $\widetilde{\Cal T}^\r$, and
$\widetilde{\Cal T}^{Sa}$, by which we mean
$$
\widetilde{\Cal T}^{BS}
=\bigcup\,\{\widehat\tau^{BS}\times\tau^\vee:\widehat\tau \in
\widehat\Sigma_P^\circ\},\quad\widetilde{\Cal T}^\r
=\bigcup\,\{\widehat\tau^\r\times\tau^\vee:\widehat\tau \in
\widehat\Sigma_P^\circ\},\,\,\,\,\text{etc.,}\tag 3.4.5
$$
where $\widehat\tau^{BS}$ is the closure of $\widehat\tau^\circ$
in $\overline C_P$, $\widehat\tau^\r$ is the closure of $\widehat\tau^\circ$
in $\overline C_P^\r$, etc.

\demo {{\rm (3.4.6)} Remark and Convention} We henceforth
understand, when we write $D^{Sa}_{\ell,P}$ below, that the latter
maps to $\widehat\Sig_P$, which is itself a minimal Satake
compactification of $\widehat C_P$ (see, e.g., [Z3:(2.10)] about
morphisms of Satake partial compactifications).  It would not
surprise me if this turned out to be an unnecessary assumption.
\enddemo

One sees directly that the diagram
$$\matrix
\widetilde{\Cal T}^{BS} & @>>> & \widetilde{\Cal T}^\r & @>>> &
\widetilde{ \Cal T}^{Sa} & @>>> &\widetilde{\Cal T}  \\
@VVV  @VVV  @VVV @VVV \\
\Dbar_{\ell,P} & @>>> & \overline D^\r_{\ell,P} & @>>> &
D^{Sa}_{\ell,P} & @>>> & \widehat\Sig_P \endmatrix\tag 3.4.7
$$
is Cartesian, with fiber given by Proposition 3.4.12 below.  The
quotients by $\Gam(\GlP)$ of the bottom row fit into a tower:
$$
\Xbar_{\ell,P}\to\Xbar^\r_{\ell,P}\to X^{Sa}_{\ell,P}\to
\widehat\Sig'_P,\tag 3.4.7.1
$$
though fixed points of $\Gam(\GlP)$ at the boundary of the latter
two prevents the $\Gam(\GlP)$-quotient of (3.4.7) from being
Cartesian.

We determine next: \proclaim{Proposition 3.4.8} The natural
mappings associated to {\rm (3.4.7)},
$$\multline
\Lb(\Dbar_{\ell,P},(\widehat\Sig_P^c)^\vee)\to\Lb(\overline D^\r_{
\ell,P},(\widehat\Sig_P^c)^\vee)\to\\
\Lb(D^{Sa}_{\ell,P},(\widehat\Sig_P^c)^\vee)\to\Lb(\widehat\Sig_1,
(\widehat\Sig_P^c)^\vee)\to(\widehat\Sig_P^c)^\vee,
\endmultline\tag 3.4.8.1
$$
\noindent satisfy $\Lb$-basechange {\rm (1.1.8)} for
$(\widehat\Sig_P^c)^\vee$. Moreover, all are homotopy equivalences
with contractible fibers, as are the projections onto
$(\widehat\Sig_P^c)^\vee$. The same holds for the mappings on
$\Lb(\,\,\cdot\,\,,\,(\widehat\Sig_P^c)^ \vee)$ induced by
morphisms of Satake compactifications.
\endproclaim
\demo{Proof} We show that
$$
\d\L(\overline C_P,(C_P)^\te)=
\Lb(\Dbar_{\ell,P},(\widehat\Sig_P^c)^\vee)\simeq\widetilde{\Cal
T}^{BS}, \tag 3.4.8.2
$$
the other cases being essentially the same.  As in the proof of
Proposition 3.3.8, one has the inclusion
$$
\Lb(\Dbar_{\ell,P},(\widehat\Sig_P^c)^\vee)\supseteq\widetilde{\Cal
T}^{BS}.\tag 3.4.8.3
$$
On the other hand, we have the tautology
$(\widehat\tau^{BS}\times\tau^\vee)
=\widehat\tau^{BS}\times_{\widehat
\Sig_P}(\widehat\tau\times\tau^\vee)$ for $\widehat\tau\in
\widehat\Sig_P$, and therefore
$$
\widetilde{\Cal
T}^{BS}\simeq\Dbar_{\ell,P}\times_{\widehat\Sig_P}\widetilde{ \Cal
T}=
\Dbar_{\ell,P}\times_{\widehat\Sig_P}\Lb(\widehat\Sig_P,(\widehat
\Sig_P^c)^\vee);\tag 3.4.8.4
$$
by (1.1.6),
$$
\Lb(\Dbar_{\ell,P},( \widehat\Sig_P^c)^\vee)\subseteq
\Dbar_{\ell,P}\times_{\widehat\Sig_P}\Lb(\widehat\Sig_P,(\widehat
\Sig_P^c)^\vee).\tag 3.4.8.5
$$
We see from (3.4.8.3), (3.4.8.4) and (3.4.8.5) that (3.4.8.2)
holds, and and then that the inclusion (3.4.8.5) is an equality.
By parallel arguments, we get the same for $\Dbar_{\ell,P}^\r$ and
$D_{\ell,P}^{Sa}$.  This yields $\Lb$-basechange for
$(\widehat\Sig_P^c)^\vee$ with respect to the mappings in the
bottom row of (3.4.7).

The fibers in (3.4.8.1), over the interior of $\tau^\vee$ in
$(\widehat \Sigma_P^c)^\vee$ are the respective closures
$\widehat\tau^{BS}$, $\widehat\tau^\r$, $\widehat\tau^{Sa}$ of
$\widehat\tau^\circ$ in $\Dbar_{\ell,P}$, $\Dbar_{\ell,P}^\r$ and
$D^{Sa}_{\ell,P}$. These closures are contractible, and the
mappings between them have contractible fibers, as can be seen
from [Z3:(3.8)]; that also covers the case of a morphism of Satake
compactifications.\sq
\enddemo
\demo{{\rm (3.4.8.6)} Remark} Lest it be forgotten, if one has
contractible fibers for the $\Lb$, one has contractible fibers for
the whole $\L$, for trivial reasons.
\enddemo

We next deduce the associated result for $\Gam(\GlP)$-quotients:

\proclaim{Proposition 3.4.9} The natural mappings induced by {\rm
(3.4.7.1)}:
$$\multline
\Lb(\Xbar_{\ell,P},\widehat\Sig^{\prime\prime}_2)\to\Lb(\overline
X^\r_{ \ell,P},\widehat \Sig^{\prime\prime}_2)\to  \\
\Lb(X^{Sa}_{\ell,P},\widehat\Sig^{\prime \prime}_2)
\to\Lb(\widehat\Sig'_1,\widehat\Sig^{\prp}_2)\to\widehat\Sig^{\prp}_2
\endmultline\tag 3.4.9.1
$$
are all homotopy equivalences with contractible fibers.
The same holds for the mappings on $\Lb(\,\,\cdot\,\,,\widehat\Sig^{\prime
\prime}_2)$
induced by morphisms of Satake compactifications.
\endproclaim
\demo{Proof} Because $\Gam(\GlP)$ acts freely on
$(\widehat\Sigma_P^c)^\vee$, Corollary 3.4.2 would again apply. We
see that again the issue in deducing our assertion from
Proposition 3.4.8 is diagonality. The argument we present for
$\Xbar_{\ell,P}$ is easily seen to apply to the other cases.

Let $\widehat\tau\in\widehat\Sigma_1^\circ$,
$e_1\in\widehat\tau^{BS}$, and $e_2\in(\tau^\vee)^\circ$. As
$\widehat\tau$ varies, this gives all points $(e_1,e_2)$ of
$\Lb(\Dbar_{\ell,P}, (\widehat\Sigma_P^c)^\vee)$, by the equality
in (3.4.7.5). Suppose that $(\gamma\cdot e_1,
e_2)\in\Lb(\widehat\Sigma_1, (\widehat\Sigma_P^c)^\vee)$ for some
$\gamma\in\Gam(\GlP)$. Let $\overline e_1$ be the projection of
$e_1$ onto $\widehat\Sig_1$.  Then $\overline e_1\in\widehat\tau$,
$e_2\in(\tau^\vee)^\circ$, and $(\gamma\cdot\overline e_1,e_2)\in
\Lb(\widehat\Sigma_1,(\widehat\Sigma_P^c)^\vee)$. By
(3.4.3.5,\,iii), $\gamma = 1$, so the actions are diagonal.  This
proves our assertion.\sq
\enddemo

% We next consider the $\Gam(\GlP)$-quotient mapping involving (3.3.6), viz.,
% $\widetilde
% {\Cal T}\to\widetilde{\Cal T'}$, and examine $\widetilde{\Cal T}$ and
% related spaces in some detail.

Having considered the mapping
$\Lb(\widehat\Sigma_1,(\widehat\Sigma_P^c)^\vee)\to(\widehat\Sigma^c_P)^\vee$
for Proposition 3.4.3, we look now at the other projection:
$\Lb(\widehat\Sigma_1,(\widehat\Sigma_P^c)^\vee)\to\widehat\Sigma_1\simeq
\widehat\Sigma_P$. For $\widehat \tau\in\widehat\Sigma_P$, put
$$\alignat 1
B(\widehat\tau) & =\bigcup\,\,\{\upsilon^\vee:\,
\widehat\upsilon\in\widehat\Sigma_P^ \circ
\text{ and $\widehat\tau$ is a face of $\widehat\upsilon$}\}\tag 3.4.10\\
 & =\bigcup\,\,\{\upsilon^\vee:\,\widehat\upsilon\in\widehat\Sigma_P^\circ
\text{ contains
$\widehat\tau$ as a codimension-one face}\}.\endalignat
$$
We invoke the treatment of duals and joins from (3.3).  The
condition that $\widehat\upsilon \in\widehat\Sig_P^\circ$ contains
$\widehat\tau$ as a codimension one face can be rewritten as
$\widehat\upsilon = \widehat\tau\ast\widehat\alpha$, where
$\widehat \alpha$ is a vertex of $\widehat\Sig_P^c$.  This allows
us to write
$$
B(\widehat\tau)=\bigcup\,\{(\tau^\vee\cap\alpha^\vee):\,\alpha\text{
as above}\} = \tau^\vee\cap (\Sig_P^c)^\vee,\tag 3.4.10.1
$$
which provides another way of looking at (3.3.8.1). When $\widehat
\tau\in\widehat\Sig_P^\circ$, $B(\widehat\tau)$ is just
$\tau^\vee$; for $\widehat\tau\subset\widehat\sigma$,
$B(\widehat\tau)\supseteq B(\widehat\sigma)$.  Also, one should
not forget (3.2.5).

% (3.4.10.2)  FIGURE ***

\proclaim{Lemma 3.4.11} Let $\widehat\tau\in\widehat\Sig_1$. Then
$B(\widehat\tau)$ is contractible.
\endproclaim

\noindent {(3.4.11.1) {\it NB---}For general triangulated spaces,
the set $B(\widehat\tau)$ need not be contractible.  Indeed, if
$\widehat\tau$ were the vertex (declared to be ``boundary") of the
cone on a circle, $B(\widehat\tau)$ would be a circle.
\smallskip

% has the homotopy type of the nerve of $\widehat\Sig_P^\circ(
% \widehat\tau)$, the latter denoting the set
% of $\widehat\upsilon$ in (the first line of) (3.4.10).

\demo{Proof of Lemma {\rm 3.4.11}}  This is easy if
$\widehat\tau\in\widehat\Sig_P^\circ$:
$B(\widehat\tau)=\tau^\vee$, which is contractible by Lemma 3.2.4.
If $\widehat\tau$ is in $\d\widehat\Sig_1$, $B(\widehat\tau)$ is
the boundary of $\tau^\vee$ in $(\widehat\Sig_P^c)^\vee$ by
(3.4.10.1). To see it is contractible, we will determine the link
of the boundary component of $\widehat C_P$ that contains
$\widehat\tau^\circ$.

Recall that the $\QQ$-root system of $\GlP$ is of type A, with
simple roots of the form $\Delta_\ell :=\{\beta_j: j<k\}$ inside
the root system of type BC or C from (2.1).  In particular, it too
has a linear Dynkin diagram, and the picture is analogous to that
of (2.1). The boundary components of $\widehat C_P$ are normalized
by maximal $\QQ$-parabolic subgroups $Q_{\ell,P}$ of $\GlP$. The
standard ones are determined by deleting a root from
$\Delta_\ell$. The parameter (a non-empty set of simple roots) for
this Satake compactification is $\{\beta_1\}$ (see [HZ2:\S 2.1]),
on the opposite end of the Dynkin diagram for $G$ from that for
(2.1.4). The selection and omission of $\beta_j\in\Delta_\ell$
splits $\Delta_\ell$ into two disjoint pieces of type A, which we
write as $\Delta_\ell^- \sqcup\Delta_\ell^+$, with $\Delta_\ell^-$
giving simple $\QQ$-roots for the automorphism group, denoted
$G^-_{\ell,P}$, of the boundary component $D_\ell^-$ of $\widehat
C_P$.  The corresponding Levi subgroup of $Q_{\ell,P}$ is an
almost-direct product of the form $G^-_{\ell,P}\cdot
G^+_{\ell,P}$, with the latter factor having $\Delta_\ell^+$ as
simple $\QQ$-roots.

We can now see that the link of $D_\ell^-$ is contractible. As allowed by 
the opening
paragraph of section (2.3), we may take $D$ there to be the non-Hermitian
$D_{\ell,P}$.
 We can identify the link on $\Dbar_{\ell,P}$. Thus, let $R=Q_{\ell,P}$ in
(2.3.2). As $Q_{\ell,P}$ is maximal, we have that $\dim A_R = 1$.
The link of $D_\ell^-$ is then given by an embedded copy of
$D_\ell^+\times W(\RR)\subset D_{\ell,P}$, where the
$A_R$-component in (2.3.2) is held constant; $W$ is the unipotent
radical of $Q_{\ell,P}$ and $D_\ell^+$ is the symmetric
space---taken in the sense of (2.1)---of $G^+_{\ell,P}$. The set
$B(\widehat\tau)$ is of the same homotopy type as the link, so it
too contractible.\sq
\enddemo

There are partial analogues of Propositions 3.4.8 and 3.4.9 ``from
the other side.''
\proclaim{Proposition 3.4.12} {\rm i)} For
$\widehat\tau\in\widehat\Sig_1$, the fiber of
$$
\Lb(\widehat\Sigma_1,(\widehat\Sigma_P^c)^\vee)\to\widehat\Sigma_1\tag
3.4.12.1
$$
over $\widehat\tau^ \circ$ is $B(\widehat\tau)$.  In particular,
$$
\Lb(\widehat\Sigma_1,(\Sig_P^c)^\vee)=\bigsqcup\,\{\widehat\tau^\circ\times
B(\widehat\tau):\widehat\tau\in\widehat\Sigma_P\}=\bigcup\,\{\widehat\tau\times
B(\widehat\tau):\widehat\tau\in\widehat\Sigma_P\},\tag 3.4.12.2
$$
and the fibers in {\rm (3.4.12.1)} are contractible.

{\rm ii)} It is analogous for the vertical mappings in {\rm
(3.4.7)}:
$$\Lb(\Dbar_{\ell,P},(\widehat\Sig_P^c)^\vee)=\bigsqcup\,\{(\widehat\tau^{BS})^
\circ\times
B(\widehat\tau):\widehat\tau\in\widehat\Sigma_P\}=\bigcup\,\{\widehat\tau^{BS}
\times
B(\widehat\tau):\widehat\tau\in\widehat\Sigma_P\},\tag 3.4.12.3
$$
where $(\widehat\tau^{BS})^\circ$ denotes the inverse image in
$\Dbar_{\ell,P}$ of $\widehat\tau^\circ$; likewise for the other
cases ($\Dbar_{\ell,P}^\r$ and $D_{\ell,P}^{Sa}$).
\endproclaim
\demo{{\rm (3.4.12.4)} Remark} For $\tau\in\widehat\Sig_1$, the
mapping $(\widehat\tau^{BS})^\circ\to\widehat\tau^\circ$ is a
homeomorphism if and only if $\widehat\tau\in\widehat
\Sig_P^\circ$. The same holds in the other cases.
\enddemo
\demo{Proof of Proposition {\rm 3.4.12}} i) That the fiber of
(3.4.12.1) over $\widehat\tau^ \circ$ is $B(\widehat\tau)$ is
contained in (3.3.8.1).  This implies the first equality in
(3.4.12.2).  The second equality follows from the fact that $B$ is
order-reversing on simplices, mentioned before Lemma 3.4.11. The
fibers are contractible by Lemma 3.4.11.

ii) From (3.4.7), there is $\Lb$-basechange for
$(\widehat\Sig_P^c)^\vee$, so the spaces occurring as fibers in
$\widetilde{\Cal T}^{BS}\to\Dbar_{\ell,P}$ are the same as those
in $\widetilde{\Cal T}^{Sa}\to D^{Sa}_{\ell,P}$. Taking
$D^{Sa}_{\ell,P}$ to be $\widehat\Sigma_P$, we have from i) that
(3.4.12.3) holds and the fibers are contractible.\sq
\enddemo

Proposition 3.4.12 has an interesting consequence:
\proclaim{Corollary 3.4.13} The mapping
$$
\Lb(\Xbar_{\ell,P},\widehat\Sig^{\prime\prime}_2)\to\Xbar_{
\ell,P}
$$ has contractible fibers, so is a homotopy equivalence.
\endproclaim
\demo{Proof} Our running hypothesis is that $\Gam(\GlP)$ acts
freely on $D_{\ell,P}$.  The same holds for its action on
$\Dbar_{\ell,P}$ (a basic feature of the Borel-Serre construction;
see [BS:9.5]).  Moreover, the actions on $\Dbar_{\ell,P}$ and
$(\widehat\Sig_P^c)^\vee$ are diagonal, as determined in the proof
of Proposition 3.4.9.  From Corollary 3.4.2, the fiber is
determined prior to taking the $\Gam(\GlP)$-quotient. From
Proposition 3.4.12, ii), the fiber is $B(\widehat\tau)$, which is
contractible by Lemma 3.4.11.\sq
\enddemo
\medskip

{\bf (3.5)} {\bf Fibers over the strata of $X^*$.}  We next
proceed to show how Theorem 3.1.1 follows from a mild variant of
Proposition 3.4.9.  The assertions in (3.4) are all about
boundaries attached to the homogeneous cones $C_P$, and these must
be brought to bear upon $D$ itself.

From \S2, all of the compactifications of $X$ under consideration
admit a morphism onto $X^*$. We restrict our attention to the
portion of these that maps to a neighborhood of the $P$-stratum of
$X^*$, and then determine the $\L$'s. This is legitimate by
Corollary 1.1.5.

For the Borel-Serre spaces, one starts with the decomposition associated to
$P$:
$$
D=D_P\times D_{\ell,P}\times A_P\times W_P\qquad (A_P\simeq\RR_{\ge 0}),\tag
3.5.1
$$
and adjoins accordingly (see (2.3)):
$$\matrix
\text{for } \Dbar\,\,\,\,\,\,\,: & \qquad D_P\times \Dbar_{\ell,P}\times\{\infty\}\times W_P;\\
\text{for } \Dbar^\e: & \qquad D_P\times \Dbar_{\ell,P}\times\{\infty\}\times V_P;\\
\text{for } \Dbar^\r: & \qquad D_P\times \Dbar_{\ell,P}^\r\times\{\infty\}\times\{1\};
\\
\text{for } D^*\,\,\,\,: & \qquad D_P\times \{pt\}\times\{\infty\}\times\{1\}.
\endmatrix\tag 3.5.2
$$
% As before, we fix a maximal parabolic subgroup $P$.  We have the
% decomposition of $D$, viewed as a homogeneous space for $P(\RR)$:
% $$
% D\simeq C_P\times D_P\times W_P(\RR).\tag 3.5.1
% $$
% In the construction of $\Dbar$, one attached the boundary face
% $$
% e(P) \simeq D_{\ell,P}\times D_P\times W_P(\RR)
% $$
% as the limit of $A_P$-orbits (cf.~(2.3.2)).  For parabolic
% subgroups $R$ that are subordinate to $P$ (see (2.2.12)), one can write
% accordingly
% $$
% e(R) \simeq\widetilde e(R_{\ell,P})\times D_P\times W_P(\RR),
% $$
% where ``$\widetilde e$'' indicates a Borel-Serre face for
% $\Dbar_{\ell,P}$. Taking into account all such $R$, we have the
% stratification
% $$
% \overline D_{\ell,P}\simeq\bigsqcup_R \,\widetilde e(R_{\ell,P}).
% $$
We recall the basic determinations (see Propositions 2.3.8 and
2.3.10): \proclaim{Proposition 3.5.3} i) The fiber of $\Dbar\to
D^*$ over $D_P$ is $\Dbar_{\ell,P}\times W_P(\RR)$.
\smallskip

ii) The fiber of $\Dbar^\e\to D^*$ over $D_P$ is $\Dbar_{\ell,P}\times
V_P(\RR)$.
\smallskip

iii) The fiber of $\Dbar^\r\to D^*$ over $D_P$ is $\Dbar^\r_{\ell,P}$.\sq
\endproclaim

We put $W'_P=\Gam(W_P)\backslash W_P(\RR)$ and $V'_P=\Gam(V_P)\backslash V_P(
\RR)$.
Passing to arithmetic quotients in Proposition 3.5.3, one gets:

\proclaim{Corollary 3.5.4} i) The fiber of $\Xbar\to X^*$ over $M_P$ is a
$W'_P$-fibration over $\Xbar_{\ell,P}$.
\smallskip

ii) The fiber of $\Xbar^\e\to X^*$ over $M_P$ is a $V'_P$-fibration over
$\Xbar_{\ell,P}$.
\smallskip

iii) The fiber of $\Xbar^\r\to X^*$ over $M_P$ is $\Xbar^\r_{\ell,P}$.\sq
\endproclaim

We will look at the above in relation to the fiber of $\Xtorexc\to
X^*$ over $M_P$. Recall from Proposition 3.1.2 that we can use
(3.5.1) in the form $U_P(\RR)\backslash D\simeq C_P\times
D_P\times V_P(\RR)$ for that purpose. The portion of the excentric
toroidal boundary over $M_P$, like that of the excentric
Borel-Serre, is adjoined by a construction on
$$
A_P\times D_{\ell,P}\simeq C_P
$$
(Again, one must not forget that the $U_P(\RR)$-quotient is
actually taken only at the boundary, and not in the interior as
the preceding may suggest.)  The Baily-Borel-type $P$-stratum of
the excentric toroidal boundary is an arithmetic quotient of
$\d(U_P(\RR)^\im)_{\Sigma_P^ \circ}\times D_P\times V_P(\RR)$.  It
has a canonical mapping onto $M_P$ (see (2.2.5.2), then take the
$(T_P^c)$-quotient, as for (2.2.18)) that is induced by projection
of that product onto $D_P$.  Since $\d U_P(\RR)_{ \Sigma^\circ_P}
\simeq(\widehat\Sigma_P^c)^\vee $, $\Gam( \GlP)\backslash\d
U_P(\RR)_{\Sigma^\circ_P}$ coincides with $\widehat
\Sig^{\prime\prime}_2$ from (3.3), the toroidal analogue of
Corollary 3.5.4 ii) takes the following form.  Put $\Gam' =
\Gam(\GlP\cdot V_P)$.

\proclaim{Proposition 3.5.5} i) The action of \, $\Gam'$ on $(\widehat
\Sigma_P^c)^
\vee \times V_P(\RR)$ is free.
\smallskip

$\quad$ ii) The fiber of $\Xtorexc\to X^*$ over $M_P$ is
$\Gam'\backslash((\widehat\Sigma_P^c)^\vee \times V_P(\RR))$. It
is a $V'_P$-fibration over $\widehat\Sig^{\prime\prime}_2$. \sq
\endproclaim

We can now proceed to determine the least common modifications. That entails
adapting
Proposition 3.4.9 to include the role of $W_P$. Let $\widetilde\Gam' = \Gam(
\GlP\cdot W_P)$.

\proclaim{Proposition 3.5.6} i) The fibers over $M_P\subset X^*$ of the natural
projections

{\eightpoint
$$
\L(\Xbar,X^\te)\to X^\te,\quad\L(\Xbarexc,X^\te)\to X^\te,\quad\L(\Xbarred,
X^\te)\to X^\te
$$}
are given respectively by the rows of the commutative diagram:

{\eightpoint $$\matrix
\Lb(\widetilde\Gam'\backslash(\Dbar_{\ell,P}\times W_P(\RR)),\Gam'\backslash(
(\widehat\Sig_P^c)^\vee
\times V_P(\RR))) & @>>> &\Gam'\backslash(\widehat\Sig_P^c)^\vee
\times V_P(\RR))\\
@VVV  @VV = V\\
\Lb(\Gam'\backslash(\Dbar_{\ell,P}\times V_P(\RR)),\Gam'\backslash(
(\widehat\Sig_P^c)^\vee
\times
V_P(\RR))) & @>>> &
\Gam'\backslash((\widehat\Sig_P^c)^\vee\times V_P(\RR))\\
@VVV  @VV = V\\
\Lb(\overline X^\r_{\ell,P},\Gam'\backslash((\widehat\Sig_P^c)^\vee\times
V_P(\RR))) & @>>>&
\Gam'\backslash((\widehat\Sig_P^c)^\vee\times V_P(\RR))
\endmatrix\tag 3.5.6.1
$$}

ii) The rows in {\rm (3.5.6.1) are surjective and have contractible fibers.
\endproclaim
\demo{Proof} The process is familiar by now: we consider the situation
before arithmetic
quotients are taken, viz.,

{\eightpoint $$\matrix
\Lb((\Dbar_{\ell,P}\times W_P(\RR)),((\widehat\Sig_P^c)^\vee
\times V_P(\RR))) & @>>> &(\widehat\Sig_P^c)^\vee
\times V_P(\RR))\\
@VVV  @VV = V\\
\Lb((\Dbar_{\ell,P}\times V_P(\RR)),((\widehat\Sig_P^c)^\vee
\times
V_P(\RR))) & @>>> &
(\widehat\Sig_P^c)^\vee\times V_P(\RR))\\
@VVV  @VV = V\\
\Lb(\Dbar^\r_{\ell,P}\times\{1\},((\widehat\Sig_P^c)^\vee\times
V_P(\RR))) & @>>>&
((\widehat\Sig_P^c)^\vee\times V_P(\RR))
\endmatrix\tag 3.5.6.2
$$}

\noindent and then verify diagonality for ii).

To do the analogue of i) over $D_P$, we may set aside the common
factor of $D_P$ in (3.5.2), and view the issue as one of
determining the corresponding $\L$'s and projections in
$\PCp(D_{\ell,P}\times A_P\times W_P(\RR))$. These are given in
$\Lb$ format in (3.5.6.2).
% We use the fact about products of $\L$'s (Proposition 1.1.***), as applied
% to $(D_{\ell,P}\times A_P)$ and $W_P(\RR)$.
We have the following inputs:
$$
\Lb(\Dbar_{\ell,P},(\widehat\Sig_P^c)^\vee)\simeq\widetilde{\Cal T}^{BS}\quad
\text{ and }\quad
\Lb(\Dbar^\r_{\ell,P},(\widehat\Sig_P^c)^\vee)\simeq\widetilde{\Cal T}^\r
$$
from (3.4.8.1); morphisms in $\PCp(A_P\times W_P(\RR))$ imply: $\Lb(W_P(\RR),
V_P(\RR))
\simeq W_P(\RR)$,
$\Lb(V_P(\RR),V_P(\RR))\simeq V_P(\RR)$, and $\Lb(\{1\},V_P(\RR))\simeq
V_P(\RR)$.
Using Lemma 1.1.2, we see we may take products and combine these to obtain:
$$
\Lb((\Dbar_{\ell,P}\times W_P(\RR)),((\widehat\Sig_P^c)^\vee
\times V_P(\RR)))\simeq \widetilde{\Cal T}^{BS}\times W_P(\RR);\tag 3.5.6.3
$$
$$\Lb((\Dbar_{\ell,P}\times V_P(\RR)),((\widehat\Sig_P^c)^\vee
\times
V_P(\RR)))\simeq \widetilde{\Cal T}^{BS}\times V_P(\RR);\tag 3.5.6.4
$$
$$\Lb(\Dbar^\r_{\ell,P}\times\{1\},((\widehat\Sig_P^c)^\vee\times
V_P(\RR)))\simeq \widetilde{\Cal T}^\r\times V_P(\RR).\tag 3.5.6.5
$$
With $W_P(\RR)$ and $V_P(\RR)$ being contractible, the horizontal morphisms
in (3.5.6.2) are,
up to homotopy, given by the first two in Proposition 3.4.8.  In particular,
the fibers in
(3.5.6.2) are contractible.

However, we must pass to the arithmetic quotients for (3.5.6.1).  For that we
can appeal to
the diagonality given in Proposition 1.3.7. Because of Proposition 3.5.5, i),
Corollary 3.4.2
applies
in determining the fibers; they are the same as in (3.5.6.2). With the aid of
Proposition 3.4.3,
we obtain from (3.5.6.3) that the fiber in the first row, over the point
represented by
$(q,w)\in (\widehat\Sigma_P^c)^\vee\times W_P(\RR)$, is
$\widehat\tau^{BS}\times [w+U_P]$ whenever $q\in (\tau^\vee)^\circ$.
Similarly, from (3.5.6.4) in the second row the fiber is
$\widehat\tau^{BS}\times [v]$, and from (3.5.6.5) in the third row of (3.5.6.1)
the fiber is $\widehat\tau^\r$.  (Here, $\widehat\tau^{BS}$ and $\widehat\tau^
\r$ are as in
(3.4.5).)  In all three cases, the fiber is contractible.
\sq
\enddemo
\demo{{\rm (3.5.6.6)} Remark} Over the interior of $\tau^\vee$,
the left-hand column in (3.5.6.1) is given by the canonical
projections
$$
\widehat\tau^{BS}\times\{w\}\to\widehat\tau^{BS}\times\{v\}\to\widehat\tau^\r
\qquad (v = w+U_P).
$$

%We recast part of the above in the language of (1.1) (cf.~Proposition 3.4.3):
%
%***\proclaim{Corollary 3.4.10} $\L$-basechange holds for $\Xtorexc$ with
%respect
%to the tower $\Xbar\to\Xbarexc\to\Xbarred$.
%\endproclaim
%\demo{Proof} On fibers over $X^*$ in the interior of $\tau^\vee\in
%\widehat
%\Sig_P^\vee$, the tower
%$$
%\L(\Xbar,\Xtorexc)\to\L(\Xbarexc,\Xtorexc)\to\L(\Xbarred,\Xtorexc)
%$$
%associated to $\Xbar\to\Xbarexc\to\Xbarred$ is given by (3.4.9).\sq
%\enddemo
%
%We use this to prove a result that contains Theorem 3.1.1 as one
%of the cases:

\proclaim{Corollary 3.5.7} For $Y$ a Borel-Serre space (i.e., $Y=\Xbar$,
$Y=\Xbarexc$,
and $Y=\Xbarred$), the mapping
$$
\L(Y,\Xtorexc)\to\Xtorexc
$$
has contractible fibers, so is a homotopy equivalence.
\endproclaim
\demo{Proof} Over $M_P\subset \d X^*$, the fiber of
$\L(Y,\Xtorexc)\to\Xtorexc$ is given by that of the corresponding
row of (3.5.6.1).  We have already seen that this is contractible
for all $P$ (Proposition 3.5.6).  We invoke the criterion from
[GT,\S 8]: a morphism of compact stratified spaces having
contractible fibers is a homotopy equivalence.\sq
\enddemo
By inverting the above homotopy equivalence in the second case, we obtain a
homotopy
class of mappings (analogous to $h:X^\t\to\Xbarred$ in [GT]):
$$
k:\Xtorexc\to\Xbarexc.\tag 3.5.8
$$
We can actually take this further by considering the other projection,
$$
\L(\Xbarexc,\Xtorexc)\to\Xbarexc\tag 3.5.9
$$
\proclaim{Theorem 3.5.10} The canonical mapping {\rm (3.5.9)} has contractible
fibers, so is a homotopy equivalence.
\endproclaim
\demo{Proof} The proof follows the same lines as that of Corollary
3.5.7. Over a point of $M_P\subset\d X^*$, the fiber of
$\L(\Xbarexc\!,\Xtorexc)\to\Xbarexc$ is given by
$$
\Lb(\Gam'\backslash(\Dbar_{\ell,P}\times V_P(\RR)),\Gam'\backslash((\widehat
\Sigma_P^c)^\vee
\times V_P(\RR)))\to\Gam'\backslash(\Dbar_{\ell,P}\times V_P(\RR)).
$$
Because $\Gam'$ acts diagonally on (3.5.6.4), and freely on
$\Dbar_{\ell,P}\times V_P(\RR)$, we again apply Corollary 3.4.2
here to determine the fiber. Thus, we may forget about $\Gam'$,
and then drop the common factor of $V_P(\RR)$.  This reduces us to
considering the fibers of $\Lb(\Dbar_{\ell,P}, (\widehat
\Sigma_P^c)^\vee)\to\Dbar_{\ell,P}$ for all $P$.  This was already
treated in Corollary 3.4.13: the fibers are contractible.\sq

%By $\Lb$-basechange with respect to the quotient mapping $\Dbar_{\ell,P}\to
%D_{\ell,P}^{Sa}$
%(***), it is the same to determine the fiber of
%$$
%\Lb(\widehat\Sigma_1,\widehat\Sigma_2)= \Lb(D_{\ell,P}^{Sa},(\widehat
%\Sigma_P^c)^\vee)\to D_{\ell,P}^{Sa}
%$$
%By Proposition 3.3.3, the fibers are all of the form $\tau^\vee$, which is
%contractible
%(Lemma 3.2.3), so we are done.\sq
%\enddemo
%
Putting Corollary 3.5.7 and Theorem 3.5.10 together, we obtain the following
fundamental
assertion, which is cited in [Z5:\S11].
\proclaim{Corollary 3.5.11} The spaces $\Xbarexc$ and $X^\te$ are
homotopy equivalent.\sq
\endproclaim
% \demo{{\rm (3.5.12)} Remark}  It is tempting to seek a stronger statement.
% For example,
% when $X$ is of $\QQ$-rank one, it is not hard to see that $\Xbarexc$ and
% $X^\te$ are
% actually {\it homeomorphic} as topological spaces, though not generally as
% compactifications of $X$.
% \enddemo

We finish this section with a further $\L$-basechange result.
\proclaim{Proposition 3.5.12} $\L$-basechange holds for $\Xbarred$
and $\Xbarexc$ with respect to the morphism $X^\t\to\Xtorexc$.
\endproclaim
\demo{Proof} To prove that the assertion holds at the boundary, we
must show that for all maximal parabolic $P$,
$$
\Lb(\d_P\Xbarexc,\d_P X^\t) = \Lb(\d_P\Xbarexc,\d_P
X^\te)\times_{\d_P X^\te}\d_P X^\t,\tag 3.5.12.1
$$
where $\d_P$ indicates the part of the boundaries mapping to
$M_P\subset \d X^*$ (i.e., $\Lb$-basechange holds there), and the
same with $\Xbarexc$ replaced by $\Xbarred$.  We are in the
$\Lb$-basechange variant of a situation from \S 1.  The torus
$T_P^c$ acts trivially on both $\d_P\Xbarexc$ and $\d_P\Xbarred$,
and $\d_P\Xtorexc$ is the $T_P^c$-quotient of $\d_P X^\t$. The
argument proving (3.5.12.1) follows the one in the proof of
Proposition 1.1.14 verbatim.\sq
\enddemo
From this, we obtain the following: \proclaim{Corollary 3.5.13} i)
The canonical mapping $\L(\Xbarexc,X^\t)\to X^\t$ has contractible
fibers.
\smallskip

ii) {\rm [GT]} The canonical mapping $\L(\Xbarred,X^\t)\to X^\t$ has
contractible fibers.
\endproclaim
\demo{Proof} From the definition (1.1.8), the $\L$-basechange
given by Proposition 3.5.12 implies that $\L(Y,X^\t)\to X^\t$ and
$\L(Y,\Xtorexc)\to\Xtorexc$ have the same fiber when either
$Y=\Xbarred$ or $Y=\Xbarexc$.  These fibers are contractible by
Corollary 3.5.7, so we are done.\sq
\enddemo
\medskip

\centerline{\bf 4. Canonical extension of homogeneous vector bundles.}
\medskip

In this section, we perform the toroidal construction of [AMRT] (see our (2.2))
on a homogeneous vector bundle $\CE_\Gam$ on $X=\Gam\backslash D$  (we drop
henceforth the subscript ``$\Gam$'').
This yields a holomorphic vector bundle $\CE^\t$ on
$X^\t$.  We show (4.4.5) that $\CE^\t$ is the canonical extension of $\CE$
to $X^\t$ in the sense of [Mu] (see [HZ1:\,3.2]) when $\CE$ is holomorphic.
It is then an
easy matter to descend $\CE^\t$ to a complex vector bundle
$\CE^\te$ on $X^\te$.
The treatment is similar in tone to the Borel-Serre construction
for the canonical extension $\Ebarred$ of $\CE$ to $\Xbarred$ given in
[Z4:\,1.10]; the latter also
gives immediately the canonical extension bundle $\Ebarexc$ on $\Xbarexc$.
\medskip

{\bf (4.1)} {\bf Standard notions.} Let $D=G/K$ as before.  Let
$E$ be a finite-dimensional complex vector space, and $\rho:K\to
\roman{GL}(E)$ a representation of $K$.  The action of $K$ on $E$
extends to one of its complexification $K(\CC)$, and thereby to
the (correct choice of) $\CC$-parabolic subgroup $\frak P$ having
$K(\CC)$ as Levi quotient. The so-called {\it compact dual} of $D$
is given as $\check D= G(\CC)/\frak P$, and it contains $D$ as an
open subset.  The $G$-homogeneous vector bundle
$\widetilde{\CE}=G\times^K E$ on $D$ descends to $\CE$ on $X$ by
taking the quotient by $\Gam$ on the left. It also extends to the
$G(\CC)$-homogeneous vector bundle $\check\CE=G(\CC)\times^{\frak
P} E$ on $\check D$, with the bundle projection given by
$$
\check\CE=G(\CC)\times^{\frak P} E\to G(\CC)\times^{\frak P}\{0\}=\check D.
\tag 4.1.1
$$

\demo {{\rm (4.1.2)} Remark} Taking $E=\{0\}$ gives $\widetilde\CE = D$ and
$\check\CE= \check D$. Thus the constructions of [AMRT] will
become a special case of ours.  For this reason, we shall cease to
talk about $D$ and $\check D$ unless that is needed to clarify the
discussion for the general homogeneous vector bundle.
\enddemo

Let $P$ be a maximal parabolic subgroup of $G$.  This determines
the open subset $\check D(P) = U_P(\CC)\!\cdot\!D$ of $\check D$
(cf.~(3.1.2.2)). It is convenient to allow the improper parabolic
$G$ here, and then $P\prec G$ for all maximal parabolic $P$ and
$\check D(G)=D$. One sees that $\check D(P)\subset\check D(Q)$
whenever $Q\prec P$, as $U_Q \supset U_P$. We point out that the
complement of $\check D(P)$ in $\check D(Q)$ has non-empty
interior.

Moving $K$ by the inverse Cayley transform for $P$ (see [HZ1:\,1.8])
determines
a basepoint for
$\check D(P)$ that is left fixed by $\GlP$.  One sees that $\check D(P)$ is
homogeneous
under the group $P'= \GhP\!\cdot\! W_P\!\cdot\!U_P(\CC)$.  As $P'\cap K
= K_{h,P}$,
where $K_{h,P} = K\cap \GhP$, one obtains the decomposition
$$
\check D(P)\simeq D_P\times V_P\times U_P(\CC)\tag 4.1.3
$$
(see (3.1.2.2); also (2.2.2)).  The action of $\GlP$ on $\check D(P)$ is
induced
by its adjoint action on $P$.  It preserves the factors in (4.1.3); in
particular, itsufficiently
small)
is trivial on the factor $D_P$.  This yields a projection of $D\subset\check
D(P)$ onto $D_P$.

Let $\check\CE(P)$ denote the restriction of $\check\CE$ to $\check D(P)$.
We can write
the bundle projection, the restriction of (4.1.1), as
$$
\check\CE(P)=P'\times^{K_{h,P}} E\to
P'\times^{K_{h,P}}\{0\}=\check D(P), \tag 4.1.4
$$
with $P'$ acting on the left.  We have
$\check\CE(P)\subset\check\CE(Q)$ whenever $Q\prec P$.
\medskip

{\bf (4.2)} {\bf Torus actions and torus embeddings, revisited.}
We form the quotient $\check\CE'_P = \Gam(P')\back\check \CE(P)$.
On $P'\times E$, the actions of $K_{h,P}$ (as for (4.1.4)) and the
$\CC$-torus $T_P=\Gam(U_P)\back U_P(\CC)$ (trivial on $E$)
commute, as $U_P$ is the center of $P'$.  Therefore, $T_P$ acts on
$\check \CE(P)'$. We get a commutative diagram:
$$\CD
\check\CE'_P @>\check\pi_2 >> \check\CE^A_P \\
@VVV   @VV{\pi_A} V \\
M'_P @>\pi_2 >> \CA_P\endCD\leqno (4.2.1)
$$
with the rows giving the $T_P$-quotients; it is the
$\Gam(P')$-quotient of
$$\CD
P'\times^{K_{h,P}}E @>>> (P'/T_P)\times^{K_{h,P}}E \\
@VVV   @VVV \\
P'\times^{K_{h,P}}\{0\} @>>> (P'/T_P)\times^{K_{h,P}}\{0\}\,.
\endCD\leqno (4.2.1.1)
$$

For any fan $\Sigma_P$ in $U_P$, let $T_{P,\Sigma_P}$ be the
corresponding torus embedding, as in (2.2).  Put
$$
\check\CE'_{P,\Sigma_P} =
\check\CE'_P\times^{T_P}T_{P,\Sigma_P}\,.\tag 4.2.2
$$
When $E=0$, this is the basic building block of the toroidal
construction in [AMRT] (see our (2.2.5.1)), and it plays the same
role here.  Moreover, the toroidal
construction (4.2.2) yields
$$
\check\pi_{2,\Sigma_P}:\check\CE'_{P,\Sigma_P} =
\check\CE'_P\times^{T_P}T_{P,\Sigma_P}\to\check\CE^A_P\,,\tag
4.2.3
$$
and this fits into a commutative diagram
$$\CD
\check\CE'_{P,\Sigma_P} @>\check\pi_{2,\Sigma_P}>> \check\CE^A_P \\
 @V VV @VV{\pi_A} V\\
 M'_{P,\Sigma_P}  @>\pi_{2,\Sigma_P}>> \CA_P\endCD\leqno (4.2.4)
$$
that extends (4.2.1). In terms of (4.2.1.1), this is
$$\CD
P'\times^{K_{h,P}}E\times^{T_P} T_{P,\Sigma_P} @>>> P'\times^{K_{h,P}}E\times^{
T_P}\{pt\} \\
 @V VV @VV  V\\
P'\times^{K_{h,P}}\{0\}\times^{T_P} T_{P,\Sigma_P} @>>> P'\times^{K_{h,P}}\{0\}
\times^{T_P}\{pt\}\endCD\leqno (4.2.4.1)
$$
In particular, (4.2.3) is a $T_{P,\Sigma_P}$-fibration. We will show that
(4.2.4),
so also (4.2.1), is Cartesian in Proposition 4.3.2 below.

Finally, we use the hypothesis that $\Sigma_P$ is $\Gam(\GlP)$-equivariant.
Then
$\Gam(\GlP)$ acts on $\check\CE'_{P,\Sigma_P}$:  for $\ell\in\Gam(\GlP)$,
$p'\in P'$,
and $t\in T_{P,\Sigma_P}$, $\ell\cdot(p',e,t) = (\ell p'\ell^\-,\ell\cdot e,
\ell\cdot t)$
in terms of (4.2.4.1).
We check this is compatible with the actions of $k\in K_{h,P}$ and $s\in T_P$:
$$\alignat 1
\ell\cdot(p's^\- k^\-, ke,st) & = (\ell p's^\- k^\-\ell^\-,\ell k\cdot e,\ell
\cdot st)\\
& = ((\ell p'\ell^\-(\ell s^\-\ell^\-) k^\-, k\ell\cdot e,(\ell s\ell^\-)\ell
\cdot t)\\
& \sim (\ell p'\ell^\-, \ell\cdot e,\ell\cdot t)
\endalignat
$$
Therefore, we can form $\Gam(\GlP)\back\check\CE'_{P,\Sigma_P}$, a partial
compactification
of $\Gam(P)\back\check\CE(P)$.
\medskip

{\bf (4.3)} {\bf Identification of canonical extensions.}  For a
topological vector bundle $\Cal V$, let $\Ct(\Cal V)$ denote the
sheaf of continuous sections of $\Cal V$; if $\Cal V$ is
holomorphic, we denote by $\CO(\Cal V)$ the sheaf of holomorphic
sections. We begin by stating some basic facts about the
unextended bundles:

\proclaim {Lemma 4.3.1} i) The $\GlP$-equivariant mapping
$\check\CE'_P\to M'_P$ in {\rm (4.2.1)} induces the mapping
$\check\CE^A_P\to\CA_P$
(by $T_P$-quotient), which is a vector bundle projection.
\smallskip

ii) $\check\CE'_P\simeq\pi_2^*\check\CE^A_P= M'_P\times_{\CA_P}\check\CE^A_P$.
\smallskip

iii) In terms of ii), the $T_P$-action on $\check\CE'_P$ is given by the
$T_P$-action
on $M'_P$ and the trivial action on $\check\CE^A_P$. In particular, $[(\pi_2)_*
\Ct(\check\CE'_P)]^{T_P} =\Ct(\check\CE^A_P)$.
\endproclaim
\demo{Proof} Statements i) and ii) are immediate from (4.2.1.1).  Note that we
can then
write $\check\pi_2$ as
$$
M'_P\times_{\CA_P}\check\CE^A_P\to \CA_P\times_{\CA_P}\check\CE^A_P;
$$
iii) is now evident.\sq
\enddemo

We continue by giving next the extension of Lemma 4.3.1 over $M'_{P,\Sigma_P}$.
Let $\pi_{2,\Sigma_P}: M'_{P,\Sigma_P}\to \CA_P$ be as in (4.2.4).

\proclaim{Proposition 4.3.2} i) Diagram {\rm (4.2.4)} is a pullback diagram.
In particular, the mapping
$\check
\CE'_{P,\Sig_P}\to M'_{P,\Sig_P}$ is a vector bundle projection, and $\check
\CE'_{P,\Sig_P}
\simeq\pi_{2,\Sig_P}^*\check\CE^A_P$.
\smallskip

ii) Moreover, $[(\pi_2)_*\Ct(\check\CE'_{P,\Sigma_P})]^{T_P} =\Ct(\check
\CE^A_P)$.
\endproclaim
\demo{Proof} Let $O$ be an open subset of $\CA_P$ over which both fibrations
$\check\pi_2$
and $b$ in (4.2.1) are trivial.  We can then trivialize (4.2.1.1) over $O$ as
$$\CD
E\times O\times T_P @>\check\pi_2 >> E\times O \\
@VVV   @VV \pi_A V \\
O\times T_P @>\pi_2 >> O\endCD\leqno (4.3.2.1)
$$
This trivialization extends when we replace $T_P$ by $T_{P,\Sig_P}$.  Thus,
the diagram
(4.2.4) is Cartesian.  All of the assertions follow.\sq
\enddemo
\demo{{\rm (4.3.2.2)} Remark}  As $T_P$ acts trivially on
$\check\CE^A_P$, we can view Proposition 4.3.2 i) as
$$
\check\CE'_{P,\Sigma_P} =
(M'_P\times_{\CA_P}\check\CE^A_P)\times^{T_P} T_{P,\Sigma_P}=
(M'_P\times^{T_P}T_{P,\Sigma_P})\times_{\CA_P}\check\CE^A_P=\pi_{2,\Sigma_P}^*
\check\CE^A_P.
$$
\enddemo
\proclaim{Corollary 4.3.3} The canonical extension of $\check\CE'_P$, from
$M'_P$ to
$M'_{P,\Sigma_P}$, is $\check\CE'_{P,\Sigma_P}$.
\endproclaim
\demo{Proof}  By definition (see [HZ1,\,3.2]), the canonical
extension of $\check\CE'_P$ is $\pi_{2, \Sigma_P}^*\check\CE^A_P$,
in view of Lemma 4.3.1 ii).  This coincides with $\check\CE'_{P,\Sigma_P}$ by
Proposition 4.3.2 i). \sq
\enddemo

To see that, as $P$ varies, the $\check\CE'_{P,\Sigma_P}$'s can be glued to
produce a vector
bundle on $X^\t$, we start by following [AMRT:\,III,\S 5].  Suppose that
$Q\prec P$,
so $\check\CE(P)\subset\check\CE(Q)$ (the complement has non-empty
interior, for
the same holds in $\check D(Q)$).  One can decompose $T_Q$ as the product of
$T_P$ and a
complementary torus---call it $T_{P,Q}$.  Assuming that $\Sigma_P\subseteq
U_P\cap
\Sigma_Q$, we get that
$$
\check\CE(Q)'_{\Sigma_P} = \check\CE(Q)'\times^{T_Q}T_{Q,\Sigma_P} = \check
\CE(Q)'\times^{T_Q}
(T_{P,Q}\times T_{P,\Sigma_P})\simeq\check\CE(Q)'\times^{T_P}T_{P,\Sigma_P}.
\tag 4.3.4
$$
From this, we see that $\check\CE(Q)'_{\Sigma_P}$ contains the
$\Gam(T_{P,Q} )$-quotient of $\check\CE(P)'_{\Sigma_P}$ as an open
set. Since $\GlP\subset\GlQ$, we obtain that
$\Gam(\GlP)\back\check\CE'_{P,\Sigma_P}\to \Gam(\GlQ)
\back\check\CE'_{Q,\Sigma_P}$ is a covering space over its image.
As $\check\CE'_{Q,\Sigma_P}$ is a dense open subset of
$\check\CE'_{Q,\Sigma_Q}$:
\proclaim{Proposition 4.3.5} $\Gam(\GlP)\back\check\CE'_{P,\Sigma_P}$ is a
covering
space over an open subset of $\Gam(\GlQ)\back\check\CE'_{Q,\Sigma_Q}$.\sq
\endproclaim

It is time to impose the usual {\it compatibility conditions} on the
collection of
fans $\Sigma=
\{\Sigma_P\}$: $\Sigma_P =\Sigma_Q\cap U_P$ whenever $Q\prec P$, and $\Sigma_{
\Int(\gamma) P}
=\Int(\gamma)\Sigma_P$ for all $\gamma\in\Gamma$.  We present next the
construction
of $\CE^\t$
(when $E=0$, it reduces to the main theorem of [AMRT,\,III,\S\S5,6], treated
in our (2.2)).

Only to put the construction of the vector bundle in ``usual''
patching format, we choose, for each $P$ (allowing, for
convenience, $S_G=D$), a {\it realm of reduction $S_P\subset D$
for} $P$. By this, we mean that $S_P$ is a sufficiently large
$\Gam(P)$-invariant open set on which $\Gam$-equivalence reduces
to $\Gam(P)$-equivalence.  If we select $S_P$ in a manner that is
compatible with the action of $\Gam$ by conjugation on the set of
parabolic subgroups, we need make only finitely many arbitrary
choices.  For instance, we can take $S_P$ to be the inverse image
in $D$ of a suitable deleted collar of the closed Borel-Serre face
$\overline{e'(P)}\subset\Xbar$. Let $\Cal B_P = \CE(S_P)$ (the
restriction of $\widetilde\CE$ to $S_P$). We define $\Cal
B'_{P,\Sigma}$ to be the interior of the closure of the image of
$\Cal B_P$ in $\Gam(\GlP)\back\check\CE'_{P,\Sigma_P}$. In
particular, $\widetilde\CE=\CE_{G,\Sigma_G}$, with $\Cal B_G =
\widetilde\CE$, $G_{\ell,G}=G$, and $\Sigma_G = \{0\}$. Let
$S'_{P,\Sigma}$ denote the same for the zero vector bundle, i.e.,
for $X$ itself, so
$$
\bigcup_P S'_{P,\Sig_P}=X^\t.\tag 4.3.6
$$

\proclaim{Proposition 4.3.7} With identifications induced from
Proposition {\rm 4.3.5}, $\bigcup_P \Cal B'_{P,\Sigma_P}$ is a
vector bundle over $X^\t=\bigcup_P S'_{P,\Sigma_P}$, which we
denote $\CE^\t$.
\endproclaim
\demo{Proof} We have, from Lemma 4.3.1 i) by restriction, that
each $\Cal B'_{P,\Sigma_P}$ is a vector bundle over
$S'_{P,\Sigma_P}$.  The restrictions of Proposition 4.3.2 and 4.3.5 imply that
they patch to
define a vector bundle over $X^\t$.\sq
\enddemo

Since the characterization of canonical extension can
be given in terms of the homogeneous bundle on the sets $M'_{P,
\Sigma_P}$ (for all $P$) [HZ1:\,3.2], we obtain immediately from
Proposition 4.3.3:

\proclaim{ Proposition 4.3.8} $\CE^\t$ is the
canonical extension (in the sense of {\rm [Mu]}) of $\CE$ from $X$
to $X^\t$.\sq
\endproclaim

We now define a (topological) vector bundle $\CE^\te$ on $X^\te$ such that
$\CE^\t$ is the
pullback of $\CE^\te$ under the quotient mapping $q:X^\t\to X^\te$
(from (2.2.18)).
We use the stratification of $X^\t\to X^\te$ to
do the same for $\CE^\t$. Note that for $P$ maximal
${}^<\!Z_P\subset S'_{P,\Sigma_P}$, and it is for ${}^<\!Z_P$ that
one takes the $T_P^c$-quotient to obtain the stratum
${}^<\!Z^\e_P$ of $X^\te$. For convenience, we write for $P=G$,
${}^<\!Z_G = X$, and use it for the interior in:
$$
X^\t = \bigsqcup_P\,{}^<\!Z_P\quad\text{ and }\quad X^\te = \bigsqcup_P\,
{}^<\!Z_P^\e.\tag 4.3.9.1
$$
We can also write
$$
\CE^\t = \bigsqcup_P \CE^\t({}^<\!Z_P)\tag 4.3.9.2
$$
though these do not give {\it open} covers.  We then put accordingly
$$
\CE^\te = \bigsqcup_P\,\CE^\te({}^<\!Z^\e_P)\to X^\te,\tag 4.3.9.3
$$
where $\CE^\te({}^<\!Z^\e_P)= \CE^\t({}^<\!Z_P)/T_P^c$, as a quotient of
$\CE^\t\to X^\t$.
Because the $T_P$ action is
given as
in Lemma 4.3.1 (iii), we see from (4.3.2.1) that
\proclaim{Proposition 4.3.10}
$\CE^\te$ is a vector bundle (i.e., locally trivial) over the space $X^\te$,
with
$q^*\CE^\te\simeq\CE^\t$.\sq
\endproclaim

% We can alternatively give the construction of $\CE^\te$ in terms of real
% torus embeddings
% (compare [HZ1:\,2.8]) and thereby avoid the use of non-open sets.  Put
% $(M'_P)^\nc = M'_P/T^c_P$, and let $g:(M'_P)^\nc\to\CA_P$ be the mapping
% induced from $\pi_2:M'_P\to\CA_P$.  Then $g$ is a principal
% fibration for
% the group $T^\nc_P = T_P/T^c_P$.  We also have the vector bundle $\check
% \CE^\nc(P)'=T^c_P\back\check\CE(P)'$, with pullback diagram
% $$\CD
% \check\CE^\nc(P)' @>\check g >> \check\CE^A_P\\
% @VVV   @VVV\\
% (M'_P)^\nc @> g >> \CA_P\endCD\tag 4.4.3
% $$
% such that the analogue of Lemma 4.3.1 holds.  Let (cf.~(4.2.1))
% $$
% \check\CE^\nc(P)'_{\Sigma_P} = \check\CE^\nc(P)'\times^{T^\nc_P}(T^\nc_P)_
% {\Sigma_P}.
% $$
% Then (4.4.3) extends to give (cf.~(4.2.4)):
% $$\CD
% \check\CE^\nc(P)'_{\Sigma_P} @>>> \check\CE^A_P\\
% @VVV   @VVV\\
% \left(M'_P\right)^\nc_{\Sigma_P} @> >> \CA_P\endCD\tag 4.4.4
% $$
% For $Q\prec P$, we have (cf.~(4.3.4)) that
%
% {\eightpoint
% $$
% \check\CE^\nc(Q)'_{\Sigma_P} = \check\CE^\nc(Q)'\times^{T^\nc_Q}(T^\nc_Q)_{
% \Sigma_P} = \check\CE^\nc(Q)'\times^{T^\nc_Q}
% (T^\nc_{P,Q}\times(T^\nc_P)_{\Sigma_P})\simeq\check\CE^\nc(Q)'\times^{
% T^\nc_P} (T^\nc_P)_{\Sigma_P}\tag 4.4.5
% $$}
% One finishes this construction of $\CE^\te$ by the method of (4.3).
\medskip

{\bf (4.4)} {\bf The Goresky-Tai conjecture.}  We recall the
statement of the conjecture (also given in [Z4:\,p.954]):
\proclaim{Conjecture A {\rm [GT:\,9.5]}} Let $h:X^\t \to\Xbarred$
be any of the continuous mappings constructed in {\rm [GT]}.  Then
the canonical extension $\CE^\t$ is topologically isomorphic to
the pullback $h^*\Ebarred$.
\endproclaim
We have been leading up to the following rather natural result:
\proclaim{Theorem 4.4.1} Suppose that for a mapping
$k:X^\te\to\Xbarexc$ as in {\rm (3.5.12)},
$\CE^\te\simeq k^*\Ebarexc$.  Then {\rm Conjecture A} is true.
\endproclaim
\demo{Proof} Since pullbacks of vector bundles are, up to
isomorphism, determined by the homotopy class of the morphism, it
is enough to show that there exist $h$ and $k$ in the homotopy
classes for which our hypothesis implies Conjecture A.

%It is convenient to simplify the notation. Let $A=\Xbarexc$,
%$B=\Xbarred$, $Y=X^\t$, and $Z=X^\te$. We have morphisms of
%compactifications of $X$: $A\to B$ and $Y\to Z$. We need to relate
%the homotopy inverses of $\L(Y,B)\to Y$ and $\L(Z,A)\to Z$, and
%the resulting homotopy classes of mappings $Y\to B$ and $Z\to A$.
%We have commutative diagrams
%$$\CD
%\L(Y,B) @>>> \L(Z,B)\times_Z Y\\
%@VVV @VVV\\
%Y @<<< Z\times_Z Y;\endCD\tag 4.4.1.1
%$$
%$$\CD
%\L(Z,A) @>>> \L(Z,B)\times_B A\\
%@VVV @VVV\\
%Z @<<< \L(Z,B).\endCD\tag 4.4.1.2
%$$
%By our $\L$-basechange results (Proposition 3.4.12 and Corollary
%3.4.10, resp.), the top rows of both diagrams are isomorphisms.
%Since $\L(Y,B)\to Y$ is a simple basechange from $\L(Z,B)\to Z$,
%we concern ourselves with the latter.  The diagram (4.4.1.2)
%displays the mapping $\L(Z,A)\to Z$ as the composite of
%$\L(Z,B)\to Z$ and one whose fibers are those of the mapping $A\to
%B$, and the latter has contractible fibers (cf. Proposition
%3.4.8).  Thus, we may assume that the homotopy inverse
%$Z\to\L(Z,A)$ is likewise decomposed.
%

We have the following picture:
$$\matrix
\CE^\te  & & & & & &               \Ebarexc\\
@VVV   & & & &                          @VVV\\
X^\te & @<\bullet << & \L(\Xbarred,X^\te) & @<\bullet << &
\L(\Xbarexc,X^\te) &
@>\bullet>k > & \Xbarexc\\
@AA q A   @AA * A  @AA * A  @VV \overline q V \\
X^\t & @<\bullet <h < & \L(\Xbarred,X^\t) & @<\bullet << & \L(\Xbarexc,X^\t) &
@>>> & \Xbarred\\
@AAA   & & & &  @AAA \\ \CE^\t & & & & & & \Ebarred
\endmatrix
$$
where (we appeal to Corollaries (3.5.7) and (3.5.13)) an arrow labelled with a
dot has a homotopy
inverse and one labelled with an asterisk is an $\L$-basechange.
When the homotopy inverses are taken, the composite in the upper
row gives $k$ and that of the lower row gives $h$. We need to
check that we can make choices so that the diagram commutes, for
then $\overline qkq=h$ and
$$
h^* \Ebarred = q^*k^*\overline q^*\Ebarred =q^*k^*\Ebarexc =
q^*\CE^\te=\CE^\t.
$$
For that, we can take homotopy inverses for both marked arrows in
the first row.  This implies (by $\L$-basechange) the same for the
second row, such that
$$\CD
X^\te @>>> \L(\Xbarred,X^\te) @>>> \L(\Xbarexc,X^\te) \\
@AAA  @AAA  @AAA  \\
X^\t @>>> \L(\Xbarred,X^\t) @>>> \L(\Xbarexc,X^\t) \endCD
$$
commutes.  This is almost what we need.  However,
$\L(\Xbarexc,X^\t)$ is not involved in the construction from [GT],
so we must check that invoking it is harmless. The diagram
$$\CD
\L(\Xbarexc,X^\t)  @> p >> \L(\Xbarred,X^\t)\\
@VV\alpha V  @VV\beta V  \\
\Xbarexc   @>\overline q >> \Xbarred\endCD
$$
commutes: $ \overline q \alpha=\beta p$.  Let $\eta$ be a homotopy inverse of
the projection
$p$.  Then
$\overline q\alpha\eta =\beta p \eta = \beta$, as desired.\sq
\enddemo

The natural setting for Conjecture A is really the excentric
compactifications.  Indeed, our investigation of these in
[HZ2:\S1] was guided by the feeling that the two had much in
common, yet they were not isomorphic as compactifications of $X$.
We emphasize this by formulating the excentric version of Conjecture A:
\proclaim{Conjecture A$'$} Let $k:X^\te\to\Xbarexc$ be any of the
continuous mappings
constructed as above.  Then the canonical extension $\CE^\te$ is
isomorphic to the pullback $k^*\Ebarexc$.
\endproclaim
\bigskip

{\bf (4.5)} {\bf Proof of the excentric Goresky-Tai conjecture.}
We abstract the set-up just a little. Let $Y_1$ and $Y_2$ be two
compactifications of $X$, with projections $p_1:\L(Y_1,Y_2)\to
Y_1$ and $p_2:\L(Y_1,Y_2)\to Y_2$.  Suppose that $p_2$ is a
homotopy equivalence, with homotopy inverse $i:Y_2\to\L(Y_1,Y_2)$.
When $X$ is a locally symmetric variety, we will be taking
$Y_1=\Xbarexc$ and $Y_2=X^\te$.

Suppose that for a vector bundle $E$ on $X$, $E$ has extensions $E_1\to Y_1$
and
$E_2\to Y_2$.  We assert:
\proclaim{Lemma 4.5.1}  In the above setting, put $k=p_1\circ i$.  Then
the
following are equivalent:
\smallskip

i) $k^*E_1\simeq E_2$,

\smallskip

ii) $p_1^*E_1\simeq p_2^*E_2$ (as vector bundles on $\L(Y_1,Y_2)$).
\endproclaim
\demo{Proof}  Write i) out as $i^*p_1^*E_1\simeq E_2$, and apply
$p_2^*$ to both sides. As $i\circ p_2$ is homotopic to the
identity, we have $p_2^*i^*$ is the identity on isomorphism
classes of vector bundles, giving ii). The other
direction is similar: apply $i^*$ to ii).\sq
\enddemo

Thus, we seek a method for verifying (4.5.1,\,ii) in the situation
of interest,  where we have $E_1 = \Ebarexc$ and $E_2 = \CE^\te$. We
continue, though, in the abstract setting.
\proclaim{Proposition 4.5.2} Let $Y_1,\,Y_2\in\PCp(X)$, $E\to X$ a vector
bundle, and
$E_1\to Y_1$,\,$E_2\to Y_2$ vector bundle extensions of $E$. Concerning the
diagram
$$\CD
E_1 @<\widetilde p_1 << \L(E_1,E_2) @>\widetilde p_2 >> E_2 \\
@V\varphi_1 VV  @VV\varphi V  @VV\varphi_2 V \\
Y_1 @< p_1 << \L(Y_1,Y_2) @> p_2 >> Y_2\,,\endCD
$$
the following are equivalent:
\smallskip

i) $\varphi$ is a vector bundle projection;
\smallskip

ii) $\L(E_1,E_2)\simeq p_1^*E_1\simeq p_2^*E_2$.
\endproclaim
\demo{Proof} We show that i) implies ii), the other direction being obvious.
By the universal property
of a pullback, there is a canonical morphism $\L(E_1,E_2)\to p_1^*E_1$ over
$\L(Y_1,Y_2)$.
By assumption, it is a morphism of total spaces of vector bundles of the same
rank (that of
$E$). To see it is an isomorphism, we may compute locally on $Y_1\times Y_2$
and thereby
assume that $E_1$ and $E_2$ are trivial. Thus, we write $E_1\simeq Y_1\times
V$ and
$E_2\simeq Y_2\times V'$, where $V=V'$ is a vector space.  Then:
$$\alignat1
E_1\times E_2 = (Y_1\times V)\times (Y_2\times V') & \simeq (Y_1\times Y_2)
\times
(V\times V'),\\
\L(E_1\times E_2)& \subset \L(Y_1,Y_2)\times(V\times V'),\\
p_1^*E_1 & = \L(Y_1,Y_2)\times(V\times\{0\}).\endalignat
$$
The morphism $\L(E_1,E_2)\to p_1^*E_1$ is, in these terms, induced
by the projection $\L(Y_1,Y_2)\times(V\times
V')\to\L(Y_1,Y_2)\times(V\times\{0\})$. We have a priori that this
contains as $E\to X$ as a subset, thus isomorphisms $\iota_x: V\to
V'$ for $x\in X$.  Our hypothesis implies the extension of
$\{\iota_x\}$, at least as mappings, to all of $\L(Y_1,Y_2)$.
Making the same argument for $p_2^*E_2$, we see that the extension
of $\iota_x$ is invertible and ii) holds.\sq
\enddemo
\demo{{\rm (4.5.2.1)} Remark} It is instructive to consider the
case $Y_1=Y_2$ with $E_1$ and $E_2$ {\it non}-isomorphic.
\enddemo
With the criterion given by Proposition 4.5.2, it remains to
verify: \proclaim{Proposition 4.5.3} The mapping
$\L(\Ebarexc,\CE^\te)\to\L(\Xbarexc, X^\te)$ is a vector bundle
projection.
\endproclaim
\demo{Proof} We must look simultaneously at the Borel-Serre and toroidal
parametrizations
of $\CE$, relative to $D_P$.  (In particular, we do not have to deal with the
$K_h$-equivariance.)
There is a $K_\ell$-equivariant mapping
$$
\Psi_P:\GlP\times A_P\times U_P\times V_P\times E\to U_P^\im\times U_P\times 
V_P
\times E,\tag 4.5.3.1
$$
given by $\Psi_P(g,u,v,e)=(gu'_0g^\-,u,v,e)$, where $u'_0$ is an
element of $C_P$ fixed by $K_\ell$. This ultimately induces the
pair of spaces: a $(T^c_P\times E)$-fibration over an arithmetic
quotient of $(\overline C_P\times V_P)$, and a
$(T_{P,\Sig_P})^\e\times E$-fibration over an arithmetic quotient
of $V_P$. The $\L$ in the case of $E=\{0\}$ (for $\Xbarexc$ and
$X^\te$) is covered by Proposition 3.5.6.  Since (4.5.3.1) is a
product with the identity mapping of $E$, we see that the action
of $\Gam(\GlP\cdot V_P)$ is diagonal for general $\CE$, and
$\L(\Ebarexc, \CE^\te)$ an $E$-fibration over $\L(\Xbarexc,
X^\te)$, as we wanted to show.\sq

\enddemo

% are characterized in terms of $U_P(\RR)$- and $T_P^c$-invariants along
% the respective Baily-Borel-type $P$-strata of $Y_1$ and $Y_2$. Write
% $$
% p_1^*E_1 = \{(e_1,(y_1,y_2)): e_1\in E_1(y_1)\}\qquad
% p_2^*E_2 = \{(e_2,(y_1,y_2)): e_2\in E_2(y_2)\},
% $$
% where $E_1(y_1)$ denotes the fiber of $E_1$ over $y_1$, and
% $E_2(y_2)$ the analogous.  We need a continuously-varying linear mapping
% $\lambda_{(y_1,y_2)}:E_1(y_1)\to E_2(y_2)$ whenever
% $(y_1,y_2)\in\L(Y_1,Y_2)$. Let $\{x_j\}$ be a sequence in $X$ that
% converges to $y_1\in Y_1$ and $y_2\in Y_2$.  Then $e_1$ lies on a
% unique $U_P$-invariant local section near $y_1$.  Take
% $e_2=\lambda_{(y_1,y_2)}(e_1)$ as the value of this section, viewed as
% $T_P^c$-invariant, in $E_2(y_2)$. Evidently, this does not
% depend on the choice of the sequence $\{x_j\}$.
% We can run the discussion in the opposite direction, starting from $E_2$.
% Since the two directions are mutually inverse, we obtain (4.5.1,\,ii),
% and Conjecture A$'$ is proved.
% \demo{{\rm (4.5.2)} Remark} We see no merit in invoking (3.5.6.4) in the
% above to parametrize $\L(Y_1,Y_2)$.

By Theorem 4.4.1, the conjecture of Goresky and Tai (Conjecture A) is likewise
proved.
\bigskip

{\eightpoint \centerline{\bf References}
\medskip

\noindent [AMRT] Ash, A., Mumford, D., Rapoport, M., Tai, Y.-S., Smooth
Compactification of Locally Symmetric Varieties, Math.~Sci.~Press,
Brookline, MA, 1975.
\smallskip

\noindent [BB] Baily, W., Borel, A., {\it Compactification of arithmetic
quotients
of bounded symmetric domains}.  Ann.~of Math.~{\bf 84} (1966), 442--528.
\smallskip

\noindent [B] Borel, A., {\it Introduction aux Groupes Arithm\'etiques}.
Hermann, 1969.

\smallskip
\noindent [BS] Borel, A., Serre, J.-P., {\it Corners and arithmetic
groups}.
Comm.~Math.~Helv.~{\bf 4} (1973), 436--491.
\smallskip

\noindent [GT] Goresky, M., Tai, Y.-S., {\it Toroidal and reductive
Borel-Serre
compactifications of locally symmetric spaces}.
Amer.~J.~Math.~{\bf 121} (1999), 1095--1151.
\smallskip

\noindent [HZ1] Harris, M., Zucker, S., {\it Boundary cohomology of Shimura
varieties, I: Coherent cohomology on toroidal compactifications.}
Ann.~Sci.~ENS {\bf 27} (1994), 249--344.
\smallskip

\noindent [HZ2] Harris, M., Zucker, S., {\it Boundary cohomology of Shimura
varieties, II: Hodge theory at the boundary}.
Invent.~Math.~{\bf 116} (1994), 243--307.
\smallskip

\noindent [J] Ji, L.,  {\it The greatest common quotient of Borel-Serre and
the toroidal compactifications},
Geom. funct.~anal.~{\bf 8} (1994), 978--1015.
\smallskip

\noindent [Mu] Mumford, D., {\it Hirzebruch's proportionality theorem in the
non-compact
case}. Invent.~Math.~{\bf 42} (1977), 239--272.
\smallskip

\noindent [N] Namikawa, Y., {\it Toroidal Compactification of Siegel Spaces}.
Lecture Notes in
Math.~{\bf 812}, Springer-Verlag, Berlin, Heidelberg, New York, 1980.
\smallskip

\noindent [S] Satake, I., {\it On compactifications of the quotient spaces
for arithmetically
defined discontinuous groups}. Ann.~of Math.~{\bf 72} (1960), 555--580.
\smallskip

\noindent [Z1]  Zucker, S.,  {\it $L_2$-cohomology of warped products and
arithmetic
groups.}  Invent.~Math.~{\bf 70} (1982), 169--218.
\smallskip

\noindent [Z2] Zucker, S., {\it $L_2$ cohomology and intersection homology
of locally
symmetric varieties, II.} Compositio Math.~{\bf 59} (1986), 339--398.
\smallskip

\noindent [Z3] Zucker, S., {\it Satake compactifications}.
Commentarii Math.~Helvetica {\bf 58} (1983), 312--343.  (Erratum
on p.~337 of {\it $L^p$-cohomology and Satake compactifications},
J.~Noguchi, T.~Ohsawa (eds.), {\it Prospects in Complex Geometry:
Proceedings, Katata/Kyoto 1989}. Springer LNM {\bf 1468} (1991),
317--339.)
\smallskip

\noindent [Z4] Zucker, S., {\it On the reductive Borel-Serre
compactification: $L^p$-cohomology of arithmetic groups (for large
$p$)}.  Amer.~J.~Math.~{\bf 123} (2001), 951--984.
\smallskip

\noindent [Z5] Zucker, S., {\it On the reductive Borel-Serre
compactification, III: Mixed Hodge structures.} Asian J.~Math.~{\bf 8}
(2004), 881--912.
\smallskip

\noindent [Z6] Zucker, S., {\it Excentric compactifications.} Quarterly
J.~Math.~(renamed Pure and App.~Math.~Quarterly) {\bf 1} (2005), 222--226.

\end